\documentclass{article}
\usepackage{amsmath}
\usepackage{amssymb}
\usepackage{graphicx}
\usepackage{algpseudocode}
\usepackage{natbib}
\usepackage{tikz}
\usepackage{endnotes}
\usepackage{amsthm}
\usepackage{pipeline-authoring}
\usepackage{informs4-compat}
\usepackage[margin=1in]{geometry}

\pipelineauthor{Haisheng Zhu}{Mitch Daniels School of Business, Purdue University}{zhu1214@purdue.edu}{}{}
\pipelineauthor{Taotao He}{Antai School of Management, Shanghai Jiao Tong University}{hetaotao@sjtu.edu.cn}{https://taotaoohe.github.io}{}
\pipelineauthor{Mohit Tawarmalani}{Mitch Daniels School of Business, Purdue University, Corresponding author.}{mtawarma@purdue.edu}{https://mohit.prof/}{}

\OneAndAHalfSpacedXI

\DeclareMathOperator{\softplus}{softplus}
\DeclareMathOperator{\proj}{proj}

\DeclareMathOperator{\conv}{conv}
\DeclareMathOperator{\conc}{conc}

\DeclareMathOperator{\vertex}{vert}

\DeclareMathOperator{\R}{\mathbb{R}}

\DeclareMathOperator{\for}{for}

\DeclareMathOperator{\z}{\boldsymbol{z}}
\DeclareMathOperator{\x}{\boldsymbol{x}}
\DeclareMathOperator{\y}{\boldsymbol{y}}
\DeclareMathOperator{\X}{\mathcal{X}}
\DeclareMathOperator{\blambda}{\boldsymbol{\lambda}}
\DeclareMathOperator{\cost}{\chi}
\DeclareMathOperator{\s}{\boldsymbol{s}}

\usepackage{xcolor}
\usepackage{svg}
\usepackage{pgfplots}
\usepgfplotslibrary{groupplots}
\pgfplotsset{compat=1.18}
\usepackage{booktabs}
\usepackage{comment}
\usepackage{bbm}
\usepackage{tikz, tikz-3dplot}
\usepackage{pgfplots}
\usepackage{tkz-fct}
\usetikzlibrary{positioning,calc}
\usepackage{etoolbox}
\usepackage{xcolor}
\usepackage[default]{mycolor}
\makeatletter
\ifdefined\tikzexternalrealjob
  \providecommand{\href}[2]{#2}
  \providecommand{\url}[1]{\texttt{#1}}

  \providecommand{\hypersetup}[1]{}
  
  \providecommand{\pdfbookmark}[3][0]{}
  
  \let\latexpipeline@origusepackage\usepackage
  \def\usepackage{\@ifnextchar[{\latexpipeline@usepkg@opt}{\latexpipeline@usepkg@noopt}}  \def\latexpipeline@usepkg@opt[#1]#2{    \def\latexpipeline@usepkg@tmp{#2}    \ifx\latexpipeline@usepkg@tmp\latexpipeline@hyperrefname
      \relax
    \else
      \latexpipeline@origusepackage[#1]{#2}    \fi
  }  \def\latexpipeline@usepkg@noopt#1{    \def\latexpipeline@usepkg@tmp{#1}    \ifx\latexpipeline@usepkg@tmp\latexpipeline@hyperrefname
      \relax
    \else
      \latexpipeline@origusepackage{#1}    \fi
  }  \def\latexpipeline@hyperrefname{hyperref}\fi
\ifdefined\tikzexternalrealjob\else\usepackage[colorlinks=true,breaklinks=true,bookmarks=true,urlcolor=links,citecolor=links,linkcolor=links,bookmarksopen=false,draft=false]{hyperref}\def\tikzsetnextfilename#1{}\fi
\makeatother
\usepackage{empheq}

\usepackage{multirow}
\usepackage{algorithm}
\usepackage{algpseudocode}
\usepackage{amsmath,amsfonts}
\usepackage{caption}
\usepackage{subcaption}
\usepackage{xfrac}
\usetikzlibrary{arrows.meta,calc,positioning}

\definecolor{masterblue}{HTML}{1F4E79}
\definecolor{candidateblue}{HTML}{DCEAF7}
\definecolor{blockgreen}{HTML}{E7F2EC}
\definecolor{cutorange}{HTML}{C65D21}
\definecolor{textgray}{HTML}{25313C}

\usepackage{endnotes}
\let\footnote=\endnote

\definecolor{links}{RGB}{204,36,29}

\usepackage{natbib}
 \bibpunct[, ]{(}{)}{,}{a}{}{,}    \def\newblock{\ } 
\TheoremsNumberedThrough     \ECRepeatTheorems

\EquationsNumberedThrough

\begin{document}

\RUNAUTHOR{Double Blind Review}

\RUNTITLE{Nonlinear discrete functions and retail revenue management}

\pipelinetitle{Optimal Pricing under Logit-Based Choice Model}

\pipelineabstract{Pricing under the finite-mixture logit (FML) consumer choice model is a foundational revenue management problem that links consumer preferences and their heterogeneous price sensitivity to firm-level profitability. Its relevance is increasing with availability of large-scale consumer data and because empirical evidence demonstrates the prevalence of consumer heterogeneity. However, we show that this problem is NP-hard. Existing methods are either restrictive (assume common price sensitivity across segments), heuristic (based on local search), or do not scale to practical problem sizes. We develop the first scalable optimization algorithm for pricing over a finite price grid without restrictive assumptions on price sensitivity. Our approach bridges pricing and convexification by developing (i) a mixed-integer linear program (MILP) reformulation via a unary price encoding, whose linear programming (LP) relaxation dominates factorable and fractional programming based relaxations, (ii) a hierarchy of relaxations based on a multipartite transportation problem, which yields an ideal (exact convex hull) LP formulation for the pricing problem that is exponential only in the number of products, and (iii) a Benders-style decomposition algorithm where the master problem coordinates a common price vector across independent per-segment subproblems using a Euclidean-deepest-cut separation oracle. Our formulation uses scaling techniques to embed the convex hull of product prices with segment-level purchase probabilities. Numerical experiments show that global optimization more than doubles revenue relative to local search for over half the instances. For baseline instances with 20 segments, 200 products, and 5 price levels, our method certifies optimality within 0.01\% gap, with a median runtime of about 2 minutes, while Gurobi leaves a 30\% gap after an hour. The method scales to problems with 10,000 product-segment combinations and solves real-data Expedia pricing instances on which the solver typically fails to close the gap.
}

\pipelinekeywords{pricing optimization, multinomial logit model, convexification, fractional programming, Benders decomposition} \HISTORY{}

\maketitle 

\section{Introduction}
Pricing optimization is a central problem in revenue management. Determining optimal prices is challenging because customers often have heterogeneous preferences \citep{ettl2020data,meyer2024interpretable,lu2020modelling} and, in a multi-product setting, demand for each product depends jointly on the prices of all products through substitution effects. Given the growing availability of customer-level purchase data and the expansion of online retailing, there is an increasing need for scalable pricing algorithms that can handle a large number of products and fine customer segmentation.

To formulate the pricing problem, various choice models have been developed to estimate the purchase probabilities of customers. Among them, the multinomial logit (\({\rm MNL}\)) model is the most prominent, which is derived from the random utility maximization (RUM) framework~(\cite{mcfadden1972conditional,falmagne1978representation,block1959random}). 
However, since \({\rm MNL}\) assumes that customers share the same utility, it fails to capture the heterogeneity in purchasing behavior among customers. To address this limitation, the finite-mixture logit $(\rm FML)$ model segments the customer population into finitely many classes. Customers within each segment follow an MNL model with segment-specific parameters. Several methods have been developed to learn $\rm FML$ from data (\cite{train2009discrete,hu2022learning}). 

A useful theoretical property of $\rm FML$ is that it can approximate an arbitrary RUM model with any precision~\citep{mcfadden2000mixed}. Moreover, it is in practical use. Based on the $\rm FML$ model, \cite{kannan2009practice} develop and implement pricing strategies for print editions, PDF editions, and print--PDF bundles at the National Academies Press, whereas \cite{lewis2005research} derived tenure-dependent subscription prices for a newspaper. \cite{li2019product} developed product-line pricing tools for Intel, and \citet{arslan2022data} optimized season-ticket and single-game prices for collegiate football. 

Although much of the existing literature focuses on continuous pricing optimization, we study pricing optimization under the ${\rm FML}$ model over finite price grids. This setting reflects the restriction of prices that is consistent with real-world pricing practice and the established literature. Such restrictions arise from marketing considerations, including consumers' responses to particular price endings (\cite{schindler1996increased,thomas2005penny}), as well as from the managerial convenience of implementing and controlling a finite price list~\citep{bitran2003overview,harsha2019practical}. 
This setting is applied in several real world pricing practice: 
\citet{ferreira2016analytics} optimize prices for Rue La La over finite price grids with prices ending in \(4.90\) or \(9.90\);
\cite{caro2012clearance} selects clearance prices from finite sets of “commercial prices” for Zara; 
\cite{harsha2019practical} develop a retail pricing system that uses product-specific ladders containing 20 admissible price points with `magic-number' endings.
To address this practical need, we develop a formulation and solution method for discrete pricing for large-scale instances under the ${\rm FML}$ model. Our methods can be used for continuous pricing by refining the discretization iteratively.

\subsection{Literature Review}
The computational challenges of pricing optimization under logit-based choice models were recognized by~\cite{hanson1996optimizing}, who showed that the profit function can be nonconcave in prices and that local search may terminate at suboptimal solutions. 
The literature has since focused on developing efficient methods for pricing optimization under various assumptions on logit-based choice models and/or problem settings.
A first line of work studies pricing under the (two-stage) nested logit model, where the first stage selects either the outside option or one of the nests and the second stage selects the product within the chosen nest. \citet{li2011pricing} assume equal price sensitivity within each nest to show that the profit is concave with respect to market shares, reducing the problem to a convex optimization along with a single-variable root-finding problem. \citet{gallego2014multiproduct} relax the assumption of common price sensitivity within each nest and develop a solution approach. \citet{davis2017pricing} study discrete pricing under the two-stage nested logit model with quality consistency constraints that impose price orderings within or across nests. They develop an exact solution method based on linear programming whose size is polynomial in the number of products, nests, and admissible prices.

Pricing under the FML model is substantially harder. \citet{li2019product} show that the concavity of profit in market shares established for MNL does not generally extend to FML. They reformulate the problem using choice probabilities of a single reference segment and express the profit as a weighted sum of quasiconcave functions. They develop a multistart gradient method that converges to stationary points but does not provide a global optimality certificate. Efforts have also been made to obtain global optimality guarantees. \citet{marandi2023exact} develop an exact algorithm for continuous pricing under the FML choice model that allows price sensitivity to vary across products and consumer segments. Their approach relies on spatial branch-and-bound that combines perspective reformulation and McCormick inequalities to obtain bounds. However, the computational cost of the algorithm increases substantially with number of products, limiting its use for large-scale problems. A scalable approach to FML pricing has also been developed, albeit under a restrictive assumption. \citet{van2022price} develop a $(1-\varepsilon)$-approximation algorithm for continuous pricing under the FML choice model allowing for product-specific but segment-independent price sensitivities. For a fixed number of segments, the running time is polynomial in the number of products but exponential in the number of segments. Their method conditions on the vector of segment no-purchase probabilities and reduces the problem to a concave maximization in reparametrized prices. Their analysis also shows that the common-sensitivity assumption is essential to their method. Specifically, when price sensitivity varies across segments, to enforce same price for all segments requires nonlinear coupling constraints that destroy the concave reformulation (their Remark 2). Because \cite{van2022price} solve the continuous problem and address  segment sizes smaller than studied here, we do not benchmark against their method; instead, Section \ref{exp1} quantifies the revenue cost of the homogeneity assumption.

To the best of our knowledge, ours is the first scalable global optimization algorithm over discrete price grids for pricing under FML without restrictive assumptions on price sensitivity.

\subsection{Contributions}
Let $n$ be the number of products, $d$ be the number of admissible price levels for each product, and $m$ be the number of segments. In this paper, we make the following contributions:
\begin{enumerate}
    \item We develop a new scalable algorithm for discrete pricing under the FML choice model that allows utility intercepts and price sensitivities to vary across both products and customer segments. This algorithm relies on an MILP formulation whose LP relaxation dominates the relaxations obtained through factorable programming (as used in generic state-of-the-art solvers)     and from fractional programming literature     with the dominance shown to be strict for various instances.
    \item We exploit the structure of the pricing problem to accelerate the algorithm. First, we develop a linear programming relaxation with $O(m(nd)^2)$ variables and $O(mn^2d)$ balance constraints. Then, we write the problem using discretization variables that lie in unary simplices and recognize that the monotonic decrease of attraction values with price implies that the coefficients in the relevant sums of bilinear terms share a common sign. Exploiting this and the domain structure, we project out the auxiliary variables associated with bilinear terms. The resulting convexification subproblem is then derived using $O((nd)^2)$ comparisons for each segment, replacing an exponential set of inequalities with $O(nd)$ constraints for each segment. 
    \item We develop an LP relaxation hierarchy based on a multipartite transportation interpretation of the pricing problem. Discrete price combinations correspond to paths in a layered graph, and higher levels capture interactions among price levels for larger subsets of products. Each level provides an upper bound on optimal profit, with higher orders strengthening the relaxation and the final hierarchy recovering an exact convex-hull description. The final hierarchy requires $O(md^{n})$ variables and constraints, reducing the $O(m2^{n(d-1)})$ variables and constraints in a typical RLT hierarchy.    \item We develop a decomposition algorithm that exploits the structure of the formulation and is of independent interest for fractional programming. The base relaxation allows price discrimination, relaxing the requirement that prices do not vary across customer segments. Then, it iteratively solves price coordination subproblems, one per customer segment, to derive separating inequalities. The master problem contains $nd$ shared pricing variables rather than a monolithic formulation with $(m+1)nd$ variables. The algorithm scales to practical sized pricing problems with many products and customer segments.
    \item We conduct numerical experiments to examine the revenue consequences of model misspecification and local search. The results show that ignoring heterogeneity in price sensitivity can lead to suboptimal pricing decisions and revenue losses. They also demonstrate that local search can converge to suboptimal solutions, highlighting the importance of both global optimization and leveraging heterogeneity in price sensitivity. We further show that the method certifies a $0.01\%$ optimality gap on baseline (20 segments, 200 products, 5 price levels) instances with a median runtime of approximately 2 minutes, where Gurobi leaves 30\% gap after an hour of computation. The method scales to 10,000 product-segment combinations and solves real-data Expedia pricing instances on which the solver typically fails to close the gap.
\end{enumerate}

\section{Price Optimization under FML}
We consider a pricing optimization problem for a set of products $N := \{1,2,\ldots,n\}$ and an outside option indexed by $0$, where the price of product $i$ is restricted to a finite set of admissible price levels $x_i\in \{p_{i1},p_{i2},\cdots,p_{id_i}\}$.
We assume that the customer population is divided into \(M=\{1,2,\cdots,m\}\) segments. Customers within the same segment share the same utility. Specifically, for a customer from segment \(c\), his utility of product \(i\) is given by
\[
\mathcal{U}_{ci}=a_{ci}-b_{ci}x_i+\epsilon_{ci},
\]
where \(a_{ci}\) represents the intrinsic preference for product \(i\), and \(b_{ci}>0\) denotes the price sensitivity, and $\epsilon_{ci}$ captures the unobserved error term that is independently and identically distributed according to a standard Gumbel distribution. Furthermore, the utility of outside option is normalized to $\epsilon_{c0}$.

Then, under the random utility maximization framework, the probability that a customer from segment \(c\in M\) purchases product \(i\in N\) under prices \(\x=(x_1,\ldots,x_n)\) is
\[
\pi_{ci}( \x):=\mathbb P\left\{\operatorname*{arg\,max}_{j\in N\cup\{0\}}\mathcal U_{cj}=i\right\}=\frac{\exp(a_{ci}-b_{ci}x_i)}{1+\sum_{j\in N}\exp(a_{cj}-b_{cj} x_j)}.
\]
For each segment $c$ and  item $i$, the term $\exp(a_{ci}-b_{ci}x_i)$ represents the attraction value of item $i$ at price $x_i$ to customers in segment $c$. The outside option is assigned a normalized attraction value of $1$, which accounts for the constant term in the denominator.

Let the share of the segment $c$ be $\omega_c$, which is non-negative so that $\sum_{c\in M}\omega_c = 1$. Let $\cost_i$ be the cost of product $i$. Under a pricing decision  $\x$, the expected profit across all customer segments  is:
\[
\Pi(\x) = \sum_{c\in M}\omega_c\sum_{i\in N}(x_i-\cost_i) \cdot \pi_{ci}(\x).
\]
Accordingly, the pricing problem can be formulated as the following discrete nonlinear program:
\begin{equation}
\label{price optimization}
\tag{\textsc{DPO}}
    \max_{\x} \, \Pi(\x)\quad \text{s.t. } x_i\in \{p_{i1},p_{i2},\cdots,p_{id_i}\}\text{ for }i\in N \text{ and } \x \in \X,
\end{equation}
where $0<p_{i1}< p_{i2} < \ldots < p_{id_i}$ defines the price ladder for item $i$, and $\X$ denotes the set of feasible price vectors satisfying any additional joint constraints on product prices. 
\begin{proposition}
\label{NP-Hard}
The  pricing problem in \eqref{price optimization} is NP-hard even with two price levels, two customer segments, and no additional joint constraints on product prices. 
\end{proposition}
The proof of Proposition~\ref{NP-Hard} shows that coordinating prices to be low or high for $n$ products across two customer segments is NP-Hard by reducing the partition problem to such an instance of \eqref{price optimization}.

\section{Solution Approaches}\label{section:solution}
We will focus on deriving MIP formulations for  the \textit{unconstrained} discrete pricing problem. Additional side constraints that are linear in both continuous and integer decision variables can be included simply. In Section~\ref{section:MIP}, we first present an LP formulation for the personalized pricing problem and introduce a family of \textit{pricing coordination inequalities} to obtain an MIP formulation for uniform pricing~\eqref{price optimization}. In Section~\ref{theoretical-properties}, we investigate the theoretical properties of our MIP formulation. Finally, in Section~\ref{section:LP}, we generalize the  pricing coordination inequalities to capture higher-order interactions, establishing a hierarchy of LP relaxations whose $n^{\text{th}}$ level yields the optimal value of~\eqref{price optimization}. 

\subsection{An exact MIP formulation: from personalized pricing to uniform pricing}\label{section:MIP}
Our solution approach consists of two main steps. First, in Theorem~\ref{them:personalized}, we establish an exact linear programming (LP) formulation for the personalized pricing problem. We then show in Theorem~\ref{them:exactnessMIP} that the uniform pricing problem~\eqref{price optimization} can be formulated as a mixed-integer linear program (MIP) by introducing additional constraints that coordinate pricing decisions across customer segments.
Instead of imposing the discrete requirement on the pricing variables, we represent them using the unary binarization~\citep{markowitz1957solution}.  Specifically, for each item $i \in N$, we introduce a vector of variables $\z_i = (z_{i1},z_{i2} \cdots, z_{id_i-1})$, and impose the following constraints:
\begin{equation}
    \label{eq:unary}
    \begin{aligned}
            &x_i =  \sum_{k=1}^{d_i} (p_{ik}-p_{ik-1})z_{ik-1}\\
            &1\geq z_{i1}\geq z_{i2}\geq \cdots\geq z_{i,d_i-1} \geq 0 \\
            &\z_i\in \{0,1\}^{d_i-1},
    \end{aligned}
    \tag{\textsc{Unary}}
\end{equation}
where, for notational convenience, we let $p_{i0}=0$  and $z_{i0} = 1$ for $i \in N$. Under this representation, a unary vector of the form
\[
\z_i=(\overbrace{1,\ldots,1}^{k},0,\ldots,0)  
\]
corresponds to the price $x_i=p_{ik+1}$, for $k \in \{0, \ldots, d_i-1\}$. The ordering constraint defines a simplex, which we denote by $\Delta^{d_i-1}$ and refer to as the \textit{unary simplex}. Together with the binary restrictions, it restricts $\z_i$ to  the vertices of this simplex.  The Cartesian product of these simplices $\Delta:=\prod_{i \in N}\Delta^{d_i-1}$  will be referred to as the \textit{unary simplotope}.
Beyond encoding the discrete pricing decisions, the unary variables also yield a lifted linear representation of the attraction functions. Specifically, for $c \in M$, $i \in N$, and $k \in \{1, \ldots, d_i\}$, we define attraction values and incremental attraction values as follows:
\[
v^c_{ik}:=\exp(a_{ci}-b_{ci}p_{ik}) \quad \text{ and } \quad \delta^c_{ik} := v^c_{ik} - v^c_{ik-1},
\]
where, for notational convenience, we let $v^c_{i0} = 0$. Since the price sensitivity parameter $b_{ci}$ is positive, the incremental attraction values $\delta^c_{i2}, \ldots, \delta^c_{id_i}$ are negative, reflecting the decrease in attraction as the price increases. Then, for each segment $c$ and product $i$, the attraction function  can be expressed as an affine function of the unary variables:
\begin{equation}\label{eq:attraction_lift}
\exp(a_{ci}-b_{ci}x_i) = \sum_{k \in [d_i]}\delta^c_{ik}z_{ik-1} =:\ell^c_i(\z_i).\tag{\textsc{Lift}}
\end{equation}
Next, we introduce normalized variables to represent the choice probabilities of each customer segment. Specifically, for each segment $c\in M$, define
\begin{equation}\label{eq:CC}
    \rho^c := \frac{1}{1+\sum_{j \in N}\exp(a_{cj}-b_{cj} x_j)} \quad \text{ and } \quad \y^c_i = \frac{\z_i}{1+\sum_{j \in N}\exp(a_{cj}-b_{cj} x_j)}\quad \for i \in N.
\end{equation}
Here, the variable $\rho^c$ represents the probability that a customer from segment $c$ chooses the outside option, while the term $\delta^c_{ik}y^c_{ik-1}$ captures the incremental contribution of the $k^{\text{th}}$ price level
to the choice probability of item $i$ for segment $c$. The normalization of the choice probabilities, together with  \eqref{eq:attraction_lift}, therefore yields that for each $c \in M$
\begin{equation}\label{eq:mip-1}
    \rho^c + \sum_{i \in N} \Bigl( \delta^c_{i1} \cdot \rho^c+ \sum_{k = 2}^{d_i}\delta^c_{ik} \cdot y^c_{ik-1} \Bigr) = 1.
\end{equation}
Furthermore, for each $c \in M$, the unary ordering constraints imply the following inequalities:
\begin{equation} \label{eq:mip-2}
    \rho^c\geq y^c_{i1}\geq\cdots
    \geq y^c_{i,d_i-1}\geq0 \quad \for\ i \in N.
\end{equation}
The uniform pricing problem~\eqref{price optimization} requires all customer segments to share a common pricing vector, which couples their choice probability variables $(\rho^c,\y^c)$ through the same unary variables $\z$. If we remove this common pricing requirement and allow each segment to have its own pricing decisions, the resulting optimization problem decomposes across segments. In the normalized variables introduced above, this relaxation is represented by the following linear program
\[
\Pi^*_{\text{LP}}:=\max_{\rho,\y} \biggl\{\sum_{c \in M} \omega_c \cdot \sum_{i\in N} \biggl(\gamma^c_{i1} \cdot \rho^c +\sum_{k=2}^{d_i}\gamma^c_{ik} \cdot y^c_{ik-1}\biggr) \biggm|~\eqref{eq:mip-1} \text{ and }~\eqref{eq:mip-2}  \biggr\},
\]
where the objective coefficients are defined as
\begin{equation*}
\begin{aligned}
&\gamma^c_{i1}:=(p_{i1}-\cost_i)v^c_{i1} \text{ and } \gamma^c_{ik}:=(p_{ik} - \cost_i)v^c_{ik}-(p_{i,k-1}-\cost_i)v^c_{ik-1} \quad \for  k \in \{2,\ldots, d_i\}.
\end{aligned}
\end{equation*}
The objective expresses each segment's expected profit in terms of the normalized unary variables. In particular, $\gamma^c_{i1} \cdot \rho^c$ represents the contribution of the lowest
admissible price level to the expected profit, while
$\gamma^c_{ik} \cdot y^c_{ik-1}$  captures the incremental contribution associated with increasing the price of item $i$ from $p_{i,k-1}$ to $p_{ik}$.
We next establish that the optimal value of this LP coincides with
the optimal profit achievable under personalized pricing, where each customer segment $c$ is offered its individual pricing vector $\x^c$. Specifically, the personalized pricing problem is given by
\begin{equation*}
\label{eq:personalized}
\begin{aligned}
\Pi^*_{\text{personalized}}
:=\max\quad
&\sum_{c\in M}\omega_c\sum_{i\in N}(x^c_i-\cost_i) \cdot \pi^c_{i}(\x^c)\\
\text{s.t.}\quad
&x_i^c\in\{p_{i1},p_{i2},\ldots,p_{id_i}\},
\quad c\in M,\ i\in N.
\end{aligned}
\end{equation*}
An optimal personalized pricing decision can be recovered from an optimal extreme point solution of the LP, as we detail in the proof of Theorem~\ref{them:personalized}.
\begin{theorem}\label{them:personalized}
$\Pi^*_{\text{LP}} = \Pi^*_{\text{personalized}}$.
\end{theorem}

The second step is to coordinate pricing decisions across customer segments using additional linear inequalities. These inequalities provide a complete polyhedral characterization of the nonlinear relationship between the pricing decision of an individual item and the choice probabilities of an individual customer segment. More generally, in Theorem~\ref{them:indhull}, we show that, for each segment $c$ and item $i$, they characterize the following \textit{individual hull}:
\begin{equation}\label{eq:hull}
\conv\Bigl\{\bigl(\phi_i(\z_i), \rho^c,\y^c \bigr) \Bigm| (\z,\rho^c,\y^c) \text{ satisfies }\eqref{eq:CC}, (x_j,\z_j)  \text{ satisfies }\eqref{eq:unary} \for j \in N  \Bigr\},\tag{\textsc{IndHull}}
\end{equation}
where $\phi_i(\z_i)$ is a non-decreasing linear response function in the unary variable $\z_i$,  with  $\sum_{k=1}^{d_i} (p_{ik} - p_{ik-1}) z_{ik-1}$ as a special case that maps the unary encoding of item $i$ to its corresponding price level. 

We show how these inequalities can be used to coordinate pricing decisions across customer segments, thereby yielding an optimal uniform pricing decision. We defer the investigation of their theoretical properties to Section~\ref{theoretical-properties}.
To derive the coordinating inequalities, we begin with the following nonlinear representation of the unary variables. Let $D^c(\z):=1+\sum_{j\in N}\ell^c_j(\z_j)$ denote the total (unnormalized) attraction of segment $c$. By~\eqref{eq:attraction_lift} and~\eqref{eq:CC}, we have $\rho^c = 1/D^c(\z)$ on the set defined by these constraints, and hence, for each $c\in M$, $i\in N$, and price level $k\in [d_i-1]$,
$z_{ik} = z_{ik} \cdot D^c(\z) \cdot \rho^c$.
The coordinating inequalities are then derived from a two-step relaxation of this nonlinear equality. First, we relax the product of the first two factors using McCormick envelopes~\citep{mccormick1976computability}, together with the identity $z_{ik}z_{ir} = z_{i\max\{k,r\}}$ for unary variables. We then homogenize the resulting inequalities using the positive variable $\rho^c$. This yields the following coordinating inequalities: for each segment $c \in M$, item $i \in N$, and price level $k \in [d_i-1]$,
\begin{subequations}\label{eq:coordinating}
        \begin{align}
        z_{ik} &\leq \biggl(1 + \sum_{j \in N}\delta^c_{j1} \biggr) y^c_{ik}
        + \sum_{r=1}^{d_i-1}\delta^c_{ir+1}y^c_{i,\max\{k,r\} }
        + \sum_{j \neq i} \sum_{r =1}^{d_j-1} \delta^c_{jr+1}\max\{0,y^c_{ik}+ y^c_{jr}- \rho^c\} \label{eq:coordinating-1} \\
        z_{ik} &\geq \biggl(1 + \sum_{j \in N}\delta^c_{j1} \biggr) y^c_{ik} + \sum_{r=1}^{d_i-1}\delta^c_{ir+1}y^c_{i,\max\{k,r\} } + \sum_{j \neq i} \sum_{r =1}^{d_j-1} \delta^c_{jr+1}\min\{y^c_{ik}, y^c_{jr}\}.  \label{eq:coordinating-2}
    \end{align}
\end{subequations}
The two inequalities differ only in their final terms, which capture the cross-item effects. Specifically, the final term in~\eqref{eq:coordinating-1} is obtained by homogenizing the lower bound $z_{ik}z_{jr} \geq \max \{0,z_{ik} + z_{jr}-1\}$, whereas the corresponding term in~\eqref{eq:coordinating-2} is obtained by homogenizing the upper bound $z_{ik}z_{jr} \leq \min \{z_{ik},z_{jr}\}$.
Incorporating the pricing coordinating inequalities and the unary binarization into the LP for the personalized pricing problem, we obtain an exact formulation of~\eqref{price optimization}:
\begin{equation}\label{eq:MIP}
\begin{aligned}
\max \quad & \sum_{c \in M} \omega_c \cdot \sum_{i\in N} \biggl(\gamma^c_{i1} \cdot \rho^c +\sum_{k=2}^{d_i}\gamma^c_{ik} \cdot y^c_{ik-1}\biggr)      \\
\text{s.t.} \quad & \eqref{eq:unary},~\eqref{eq:mip-1}, ~\eqref{eq:mip-2} \text{ and } \eqref{eq:coordinating} \for c\in M, i \in N,  k \in [d_i-1].
\end{aligned}\tag{\textsc{PriceConv}}
\end{equation}
As we show in Section~\ref{theoretical-properties}, the coordinating inequalities~\eqref{eq:coordinating} provide an explicit description of~\eqref{eq:hull}.
\begin{theorem}\label{them:exactnessMIP}
    The MIP formulation~\eqref{eq:MIP} solves~\eqref{price optimization} exactly.
\end{theorem}

The MIP formulation \eqref{eq:MIP} differs from typical formulations of \eqref{price optimization}. First, the discretization~\eqref{eq:unary} converts the utility function and $(x_i-\chi_i)\exp(a_{ci}-b_{ci}x_i)$ to linear functions (see \eqref{eq:attraction_lift} and $\gamma^c_{i1}z_{i1} + \sum_{k=2}^{d_i} \gamma^c_{ik}z^c_{ik-1}$, respectively). It is useful to compare with \citet{van2022price}, which converts the problem to a convex optimization problem once the non-purchase probability $\rho^c$ for all customer-segments is fixed. With discretization, this simplification is implicit for \eqref{eq:MIP} since $y^c_{ik}=\rho^cz_{ik}$, and the price coordination problem is an LP. Second, after discretization, the problem reduces to a 0-1 linear fractional programming problem for which the relaxation technique of \cite{csen2018conic} is one of the prominent approaches. Instead, we introduce \eqref{eq:mip-2} and \eqref{eq:coordinating} by scaling the simplotope containing the unary variables and relaxing the bilinear expression obtained by multiplying $z_{ik}$ with $D^c(\z)\rho^c$. Third, in \eqref{eq:coordinating}, we relax the product $\rho^cz_{ik}z_{jr}$ using scaled McCormick constraints, $\min\{y^c_{ik},y^c_{jr}\}\ge \rho^cz_{ik}z_{jr}\ge \max\{0,y^c_{ik}+y^c_{jr}-\rho^c\}$. The literature typically uses products of price level indicators \citep{cohen2021promotion,abdolhamidi2025tactical}. Instead, we relax the unary variable products because they have favorable convexification properties. The next section analyzes \eqref{price optimization} theoretically, developing insights into our choices.

\subsection{Theoretical Properties}\label{theoretical-properties}
This subsection investigates theoretical properties of \eqref{eq:MIP}. More specifically, in Section~\ref{section:hull}, we establish that the pricing coordinating inequalities provide an explicit description of the individual hull defined in~\eqref{eq:hull}. In Section~\ref{section:full_binarization}, we motivate our choice of unary binarization by examining the limitations of the prevalent full binarization approach in the pricing literature.  Last, in Section~\ref{section:twoformulations}, we present theoretical comparisons between our formulation and two alternative formulations.
\subsubsection{Characterizing the individual hull}\label{section:hull}
Deriving the individual hull is challenging for two reasons: the underlying choice probabilities involve multiple ratios whose denominator couples pricing decisions across items, and the underlying domain is the Cartesian product of unary simplices rather than a box. Surprisingly, the structure of the problem, particularly that of the unary variables, allows us to characterize the individual hull by reducing its convexification to that of pairwise bilinear terms over unit squares. The proof relies on three decomposition results established in~\ref{app:indhull}. We provide a brief sketch of the major ideas. First, Lemma~\ref{lemma:projective} reformulates the description of~\eqref{eq:hull} as a problem of \textit{simultaneously convexifying} the cross-effect between the unary pricing decision of item $i$ and the total attraction to the segment $c$:
\[
\conv\Bigl\{\bigl(\z, \phi_i(\z_i) \cdot D^c(\z) \bigr)   \Bigm| \z_j \in \Delta^{d_j-1} \cap\{0,1\}^{d_j-1} \for j \in N \Bigr\},
\]
where recall that  $\phi_i(\z_i)$ is  a non-decreasing affine function and  $D^c(\z)$ denotes the total attraction to segment $c$. Next, Lemma~\ref{lemma:common_simplex} shows that this simultaneous convexification decomposes into individual convexification problems for each item pair. Specifically, for each $j\in N$, recall that $\ell^c_j (\z_j)$ denotes the attraction value of $j^{\text{th}}$ item to a customer in the $c^{\text{th}}$ segment. We consider
\[
\conv\Bigl\{(\z_i, \z_j, \phi_i(\z_i) \cdot \ell^c_j(\z_j) )   \Bigm| \z_i \in \Delta^{d_i-1} \cap\{0,1\}^{d_i-1}, \z_j \in \Delta^{d_j-1} \cap\{0,1\}^{d_j-1} \Bigr\}, 
\]
which involves $(d_i-1)(d_j-1)$ bilinear terms over two unary simplices. These two steps use the  projective-convexification and decomposition machinery, exploiting multilinear envelope results for fractional programming. The resulting formulation lends itself to a decomposition algorithm developed later in Section~\ref{sec:decomposition}, with pricing being the first illustration. Finally, Lemma~\ref{lemma:common_triangulation} shows that this pairwise convexification can be reduced to convexifying individual bilinear terms over the unit square:
\[
\conv\bigl\{(z_{ik},z_{jr},z_{ik}z_{jr}) \bigm| z_{ik} \in \{0,1\},z_{jr} \in \{0,1\} \bigr\},
\]
a new observation exploiting the structure of \eqref{price optimization}.
Each such convex hull admits an exact description through four inequalities~\citep{mccormick1976computability}.
\begin{theorem}\label{them:indhull}
    For each segment $c$ and item $i$, a point $(\mu_i,\rho^c,\y^c)$ satisfies~\eqref{eq:hull} if and only if there exists a unary vector $\z_i$ such that 
    \begin{equation}\label{eq:indhull_rep}
    \begin{aligned}
        \mu_i = \phi_i(\z_i)  \quad \text{ and } \quad  (z_{ik}, \rho^c,\y^c)\text{ satisfies}~\eqref{eq:mip-1},~\eqref{eq:mip-2},\text{ and }~\eqref{eq:coordinating} \for k \in \{1, \ldots, d_i-1\}. 
    \end{aligned}
    \end{equation}
\end{theorem}

\begin{remark}\label{rmk:simult}
 The proof of Theorem~\ref{them:indhull} shows that the hulls of multiple response functions,  
\begin{equation}\label{eq:indhull_rep_simul}
\begin{aligned}
        \mu_{ij} \le \phi_{ij}(\z_i)  \quad \text{ and } \quad  (z_{ik}, \rho^c,\y^c)\text{ satisfies}~\eqref{eq:mip-1},~\eqref{eq:mip-2},\text{ and }~\eqref{eq:coordinating} \for k \in \{1, \ldots, d_i-1\},
\end{aligned}
\end{equation}
can be derived simultaneously, as the underestimators share the same triangulation of $\prod_{t=1}^n\Delta^{d_{t}-1}$.
\end{remark}
Observe that \eqref{eq:indhull_rep} does not include linearizations of the bilinear variables $z_{ik}z_{jr}$, thus reducing the variables from linearizing $\rho^cz_{ik}z_{jr}$ for all $c\in M$, $i,j\in N$, $k\in \{1,\ldots,d_i\}$, $r\in \{1,\ldots,d_j\}$, which are $O(mn^2d^2)$, to $O(mnd)$. In the following $\phi_i(\z_i)$ will serve as an expression that coordinates the price across customer segments. Specifically, $\phi_i(\z_i)$ could be the expression defining $x_i$ in \eqref{eq:unary}, in which case Theorem~\ref{them:indhull} gives the convex hull of the price of product $i$ with the probability of purchase for segment $c$ at each price level.

\subsubsection{Unary versus full binarization}\label{section:full_binarization}
The \textit{full binarization}, also known as  \textit{one-hot encoding}, is a prevalent binarization approach in the pricing literature~\citep{cohen2017impact,ke2019approximate,harsha2019practical,cohen2021promotion,guan2026randomized}.  Specifically, for each item $i\in N$, we introduce a vector of binary variables  $\blambda_i = (\lambda_{i1} , \ldots, \lambda_{id_i} )$, and impose the following constraints:
\begin{equation}\label{eq:full}
    x_i = \sum_{k=1}^{d_i}p_{ik} \lambda_{ik}, \ \sum_{k=1}^{d_i} \lambda_{ik} = 1, \text{ and } \blambda_i \in \{0,1\}^{d_i}. \tag{\textsc{Full}}
\end{equation}
Under this encoding, $\lambda_{ik}=1$ if and only if $x_i=p_{ik}$. Accordingly, the one-hot encoding counterpart of the individual hull can be defined as follows: for each segment $c\in M$ and item $i\in N$,
\begin{equation}\label{eq:indihull-full}
  \begin{aligned}
       \conv\biggl\{( x_i, \rho^c, \s^c) \biggm| &\rho^c = \frac{1}{1+\sum_{j \in N}\exp(a_{cj}-b_{cj} x_j)},\  \s^c = \blambda \cdot \rho^c,\ \\
       & \qquad (x_j,\blambda_j) \text{ satisfies }\eqref{eq:full} \for j \in N \biggr\},
  \end{aligned} 
\end{equation}
where $\blambda$ denotes the concatenation of the one-hot encoding vectors $\blambda_j$ over all items $j \in N$.

The decomposition approach developed in Section~\ref{section:hull} does not extend directly to the one-hot encoding. Although the first two steps of the decomposition, established in Lemmas~\ref{lemma:projective} and~\ref{lemma:common_simplex} can be generalized to this setting, the last step, established in Lemma~\ref{lemma:common_triangulation}, does not generalize. In particular, the reduction to individual bilinear terms relies on the ordered structure of the unary encoding, which is absent under one-hot encoding. The next example illustrates that termwise relaxations are insufficient to capture the relationships with one-hot encoding but the unary encoding allows these relationships to be captured.

\newcommand{\zprod}{Z}
\newcommand{\lprod}{\Gamma}

\begin{example}[One-hot versus unary binarization]\label{ex:onehot-vs-unary}
Consider two items, each with three admissible price levels, and let
$\blambda_1=(\lambda_{11},\lambda_{12},\lambda_{13})$ and
$\blambda_2=(\lambda_{21},\lambda_{22},\lambda_{23})$,
$\blambda_i\in\vertex(\Sigma^3)$, where $\Sigma^3:=\left\{
\blambda_j\in\mathbb R_+^{3}\;\middle|\;
\sum_{k=1}^{3}\lambda_{jk}=1 \right\}$, denote the one-hot (full-binarization)
price variables. We focus on the products of two linear functions of
these variables,
\begin{align*}
    \mu_1 &\le 2\lambda_{13}\lambda_{21}+\lambda_{13}\lambda_{22}
        = \lambda_{13}\bigl(2\lambda_{21}+\lambda_{22}\bigr),\\[2pt]
    \mu_2 &\le 2\lambda_{12}\lambda_{21}+\lambda_{12}\lambda_{22}
        +2\lambda_{13}\lambda_{21}+\lambda_{13}\lambda_{22}
        = (\lambda_{12}+\lambda_{13})\bigl(2\lambda_{21}+\lambda_{22}\bigr),
\end{align*}
which mirror the pricing structure. The factor $2\lambda_{21}+\lambda_{22}$ takes the value $2,1,0$ on increasing price levels, modeling the decreasing attraction value for the item $2$ in the pricing application. The coefficients of $(\lambda_{11},\lambda_{12},\lambda_{13})$ in $\lambda_{13}$ and $\lambda_{12}+\lambda_{13}$, take the values
$(0,0,1)$ and $(0,1,1)$, i.e.\ they are \emph{non-decreasing}. This monotonicity (equivalently, the common sign of the coefficients in the unary representation) is used in Lemma~\ref{lemma:common_triangulation} and Remark~\ref{rmk:simult} to derive the concave envelope for $\mu_1$ and $\mu_2$ simultaneously using termwise relaxation in the unary encoding space. The current example provides a concrete setting that illustrates that result. We will also show that if the direction of constraints on $\mu_1$ and $\mu_2$ is switched simultaneously, a similar argument continues to hold. In the following, let $\lprod_{kt}$ (resp. $\zprod_{kt}$) denote a variable representing $\lambda_{1k}\lambda_{2t}$ (resp. $z_{1k}z_{2t}$) for $k,t\in\{1,2,3\}$ (resp. $k,t\in \{1,2\}$).

\paragraph{Set 1 (one-hot, RLT/balance constraints).}
Introduce $\lprod_{kt}\ge 0$ and impose the balance (RLT) constraints
\begin{equation}\label{eq:transport}
    \sum_{t=1}^3\lprod_{kt} = \lambda_{1k},\quad k=1,2,3,
    \qquad
    \sum_{k=1}^3\lprod_{kt} = \lambda_{2t},\quad t=1,2,3,
\end{equation}
together with $\sum_k\lambda_{1k} = 1$, $\sum_t\lambda_{2t} = 1$,
$\lprod_{kt} \ge 0$, and the defining inequalities
\begin{align*}
    \mu_1 \le 2\lprod_{31} + \lprod_{32},
    \qquad
    \mu_2 \le 2\lprod_{31} + \lprod_{32} + 2\lprod_{21} + \lprod_{22}.
\end{align*}
Projecting out $\boldsymbol\lprod$ yields $13$ inequalities including
\begin{align*}
    & \mu_2 + \lambda_{22} + 2\lambda_{23} \le 2.\\
\end{align*}

\paragraph{Set 2 (one-hot, termwise McCormick).}
Alternatively, relax each product by its McCormick envelope,
\begin{align*}
    \max\{0,\, \lambda_{1k} + \lambda_{2t} - 1\}
       \le \lprod_{kt}
       \le \min\{\lambda_{1k},\, \lambda_{2t}\},
       \qquad k,t = 1,2,3,
\end{align*}
together with $\boldsymbol\lambda_i\in \Sigma^3$ and the same linking inequalities
as in Set~1. This relaxation is \emph{not} exact.
Specifically, consider the point
\begin{align*}
    \boldsymbol\lprod =
    \begin{bmatrix}
        \tfrac13 & 0 & 0\\
        \tfrac13 & 0 & 0\\
        \tfrac13 & \tfrac13 & 0
    \end{bmatrix},
    \quad
    \blambda_1 = \Bigl(\tfrac13, \tfrac13, \tfrac13\Bigr),
    \quad
    \blambda_2 = \Bigl(\tfrac13, \tfrac23, 0\Bigr),
    \quad
    \mu_1 = 1,
    \quad
    \mu_2 = \tfrac53.
\end{align*}
By direct verification, this point satisfies McCormick inequalities. 
Yet the row sums of $\boldsymbol\lprod$ are $(1/3,1/3,2/3)$ and its
column sums are $(1,1/3,0)$, which do not match $\blambda_1$ and
$\blambda_2$, and hence violate \eqref{eq:transport}. 
Moreover, the point violates the valid inequality
$\mu_2+\lambda_{22}+2\lambda_{23}\le 2$ of Set~1, since the left-hand
side equals $5/3 + 2/3 = 7/3$. Thus, the termwise
McCormick relaxation in the one-hot variables does not describe the
convex hull.

\paragraph{Set 3 (unary, termwise McCormick).}
Now encode the same price choices with the unary variables of
Section~\ref{section:hull}. For $i\in\{1,2\}$ and $k\in\{2,3\}$,
\begin{align*}
    z_{i,k-1} = \sum_{j=k}^{3} \lambda_{ij},
\end{align*}
which satisfy the unary (nested) constraints
\begin{align*}
    0 \le z_{12} \le z_{11} \le 1,
    \qquad
    0 \le z_{22} \le z_{21} \le 1.
\end{align*}
In these variables, the two products become
\begin{align*}
    \mu_1 \le 2z_{12} - \zprod_{21} - \zprod_{22},
    \qquad
    \mu_2 \le 2z_{11} - \zprod_{11} - \zprod_{12}.
\end{align*}
Then, $\zprod_{ab}$, which represents $z_{1a}z_{2b}$, is relaxed by
its McCormick envelope
\begin{align*}
    \max\{0,\, z_{1a} + z_{2b} - 1\}
       \le \zprod_{ab}
       \le \min\{z_{1a},\, z_{2b}\},
       \qquad a,b \in \{1,2\}.
\end{align*}
Projecting the resulting $(\z,\boldsymbol\zprod,\boldsymbol\mu)$-set back to the $(\boldsymbol\lambda,\boldsymbol\mu)$-space
recovers \emph{exactly} the description in Set~1, which is a direct
instance of Lemma~\ref{lemma:common_triangulation}.
For instance, combining the lower McCormick bounds
\begin{align*}
    z_{11} + z_{21} - 1 \le \zprod_{11},
    \qquad
    z_{11} + z_{22} - 1 \le \zprod_{12},
\end{align*}
with $\mu_2 \le 2z_{11} - \zprod_{11} - \zprod_{12}$ gives
$\mu_2 + z_{21} + z_{22} \le 2$; substituting
$z_{21} = \lambda_{22} + \lambda_{23}$ and $z_{22} = \lambda_{23}$
yields:
\begin{align*}
    \mu_2 + \lambda_{22} + 2\lambda_{23} \le 2,
\end{align*}
an inequality not implied by Set~2. The remaining inequalities follow analogously.

\medskip
\noindent
Section~\ref{ex:reverse} shows that when the inequalities defining $\mu_1$ and $\mu_2$ are switched from $\le$ to $\ge$, the conclusions are similar. This example motivates the discussion in this section through three observations.
First, Set~2 shows that the last step of the decomposition in Section~\ref{section:hull} fails under one-hot encoding. Second, Set~3 shows
that the termwise McCormick is \emph{exact} when the variables are unary and
the two linear factors have consistent signs, the content of Lemma~\ref{lemma:common_triangulation}. Third, the exactness has direct computational consequence. Separation with the one-hot convex hull using Set~1 requires solving an LP with $O(d_1d_2)$ variables and $O(d_1+d_2)$ balance constraints. On the other hand, the unary formulation reduces the separation to
$O(d_1d_2)$ comparisons and avoids the introduction of $\boldsymbol{\zprod}$ variables. \Halmos
\end{example}

As shown in the case of Set~2 in Example~\ref{ex:onehot-vs-unary}, independently convexifying each bilinear term $\lambda_{ik}\lambda_{jt}$, for $i\ne j$, using McCormick inequalities is insufficient to capture the coupling between these bilinear terms. If the relaxation technique of \citet{he2025convexification} is used with one-hot encoding, the natural set to convexify would be \eqref{eq:indihull-full}. The convex hull of \eqref{eq:indihull-full} can be described in the space of variables that represent these bilinear terms by imposing marginal consistency requirements such as \eqref{eq:transport} in Example~\ref{ex:onehot-vs-unary}. More specifically, consider a segment $c \in M$ and item $i\in N$. For each $j\in N\backslash\{i\}$ define a matrix of variables $\boldsymbol{\Gamma}^c_j \in \R^{d_i \times d_j}$ that represents the normalized pairwise effects $\Gamma^c_{jkt} := \rho^c\lambda_{ik}\lambda_{jt}$, where $k\in \{1,\ldots,d_i\}$ and $t\in \{1,\ldots,d_j\}$ and, for each $j\in N$, let $s^c_{jt} := \rho^c\lambda_{jt}$ for $t\in \{1,\ldots,d_j\}$. Then, $(\blambda_i,\rho^c,\s^c)$ satisfies~\eqref{eq:indihull-full} if and only if:
\begin{subequations}\label{eq:onehot-convexhull}
\begin{align}
     &x_i = \sum_{k=1}^{d_i} p_{ik} \lambda_{ik}\label{eq:onehot-definex}\\ 
     &\lambda_{ik} = (1+v^c_{ik}) \cdot s^c_{ik} + \sum_{j \neq i} \sum_{t=1}^{d_j} v^c_{jt} \cdot \Gamma^c_{jkt}  \label{eq:onehot-defn} &&\for k\in [d_i]\\
     &1 = \rho^c + \sum_j\sum_{t=1}^{d_j} v^c_{jt}s^c_{jt} \quad \text{and}\quad \rho^c\ge 0\label{eq:onehot-normalize}\\
     & \sum_{k=1}^{d_i} s^c_{ik} = \rho^c \quad \text{and} \quad 0 \leq \s^c_{j} && \for j \in N \label{eq:onehot-beginflow}\\
     &\sum_{\bar{t}=1}^{d_j} \Gamma^c_{jk\bar{t}} = s^c_{ik} \text{ and } \sum_{\bar{k}=1}^{d_i} \Gamma^c_{j\bar{k}t} = s^c_{jt} \quad &&\for j \in N \setminus \{i\}, k\in [d_i], \text{ and } t \in [d_j] \label{eq:onehot-balance}\\
     &\Gamma^c_{jkt} \ge 0 &&\for j\in N\setminus \{i\}, k\in [d_i], \text{ and } t\in [d_j].\label{eq:onehot-nonnegative}
\end{align}
\end{subequations}

\begin{proposition}\label{prop:onehot-convexhull}
\eqref{eq:onehot-convexhull} provides a lifted representation of 
\begin{equation*}
  \begin{aligned}
       \conv\biggl\{(x_i,\blambda_i, \rho^c, \s^c) \biggm| &\rho^c = \frac{1}{1+\sum_{j \in N}\exp(a_{cj}-b_{cj} x_j)},\  \s^c_j = \rho^c\blambda_j,\ \\
       & \qquad (x_j,\blambda_j) \text{ satisfies }\eqref{eq:full} \for j\in N\biggr\}.
  \end{aligned} 
\end{equation*}
In particular, its projection to $(x_i,\rho^c,\s^c)$ is the one-hot individual hull in~\eqref{eq:indihull-full}. 
\end{proposition}

This extended formulation introduces $\sum_{j\in N\setminus\{i\}}d_i d_j$ additional $\boldsymbol{\Gamma}^c$ variables to represent the normalized pairwise effects, whereas the unary-based description (Set 3 in Example~\ref{ex:onehot-vs-unary}) does not need such additional variables. Consequently, the one-hot approach increases the size of the formulation substantially. In Section~\ref{sec:example_projection}, for several example settings, we show that separating in unary space is significantly faster.

\subsubsection{Comparisons with two alternative formulations}\label{section:twoformulations}
The tightness of our formulation is the result of several steps in our construction. First, in Section~\ref{section:full_binarization}, we showed the role unary binarization plays in deriving a tight formulation. The second source of tightness is the use of inequalities \eqref{eq:coordinating} to bound $z_{ik}$. Recall that these inequalities aggregate termwise envelopes for bilinear terms in $z_{ik}D^c(\z)$ and then scale the estimator by $\rho^c$. In this subsection, we focus mostly on this second source. To isolate this benefit, we will formulate the problem for competing approaches using unary binarization. 

The first formulation is derived using factorable programming (FP), a framework adopted by several state-of-the-art general-purpose MINLP solvers~\citep{tawarmalani2004global,bolusani2024scip}. Assume pricing decisions and attraction value functions are binarized using unary variables $\z$. 

FP introduces a variable $\mu^c$ to present the profit associated with $c^{\text{th}}$ segment as follows:
\[
 \mu^c=\sum_{j\in N}\frac{\gamma^c_{j1}+\sum_{k=2}^{d_j}\gamma^c_{jk}z_{jk-1}}{D^c(\z)}.
\]
Introducing $d^c$ to represent $D^c(\z)$, this is equivalent to the bilinear equality:
\[
\mu^c d^c=\sum_{j\in N}\left(\gamma^c_{j1}+\sum_{k=2}^{d_j}\gamma^c_{jk}z_{jk-1}\right).
\]
Given $\mu^c\in [\underline \mu^c,\overline \mu^c]$ and $d^c\in [\underline d^c,\overline d^c]$,
the McCormick envelope of $\mu^c d^c$ then constructs a relaxation:
\begin{equation}\label{FPMCC}
        \begin{aligned}
        \underline\mu^c d^c+\overline d^c\mu^c-\underline\mu^c\overline d^c\geq
        \sum_{j\in N}\gamma^c_{j1}+\sum_{j\in N}\sum_{k=2}^{d_j}\gamma^c_{jk}z_{jk-1}\ge \overline\mu^c d^c+\overline d^c\mu^c-\overline\mu^c\overline d^c,c\in M,\\
        \overline\mu^c d^c+\underline d^c\mu^c-\overline\mu^c\underline d^c\geq
        \sum_{j\in N}\gamma^c_{j1}+\sum_{j\in N}\sum_{k=2}^{d_j}\gamma^c_{jk}z_{jk-1}\ge \underline\mu^c d^c+\underline d^c\mu^c-\underline\mu^c\underline d^c,c\in M.\\
    \end{aligned}
\end{equation}
We consider a basic factorable relaxation in which the bounds used to construct the McCormick inequalities are derived from the component-wise restrictions $0\le z_{jk}\le1$. Specifically, let $\underline\mu^c$ and $\overline\mu^c$ denote valid lower and upper bounds on $\mu^c=\sum_{j\in N}\frac{\gamma^c_{j1}+\sum_{k=2}^{d_j}\gamma^c_{jk}z_{jk-1}}{D^c(\z)}$ over $\z_j\in[0,1]^{d_j-1}$ for all $j\in N$. 
Since $\exp(a_{cj}-b_{cj}x_j)$ is decreasing in $x_j$ the bounds on $d^c=D^c(\z)$ over this box are
\[
\underline d^c:=1+\sum_{j\in N}\exp(a_{cj}-b_{cj}p_{jd_j}),\qquad
\overline d^c:=1+\sum_{j\in N}\exp(a_{cj}-b_{cj}p_{j1}),\qquad c\in M.
\]
These bounds are attained when all $z_{jk}=1$ and all $z_{jk}=0$, respectively.
The relaxation is then:
\begin{equation}
   \begin{aligned}\label{fp}
\max_{\mu,\bold d,\z}\quad & \sum_{c\in M}\omega_c\mu^c\\
\text{s.t.}\quad & d^c=D^c(\z), && c\in M,\\
& \underline\mu^c\le\mu^c\le\overline\mu^c,\quad \underline d^c\le d^c\le\overline d^c, && c\in M,\notag\\
& \z_j\in\Delta^{d_j-1}\cap \{0,1\}^{d_j-1}, && j\in N,\notag\\
&\eqref{FPMCC}.
\end{aligned}
    \tag{\textsc{FactProg}}
\end{equation}
\begin{proposition}\label{prop:PD-dominates-FP}
    The LP relaxation of~\eqref{eq:MIP} yields a tighter upper bound for~\eqref{price optimization} than that of~\eqref{fp}.  
\end{proposition}
\begin{example}\label{ex:dominateFactProg}
We show that the dominance can be strict and also compare with the root relaxation from SCIP. Consider pricing with two products and two customer segments (assuming $\cost_j = 0$ for $j=1,2$):
\[
\max_{\substack{x_1\in\{3/2,2\}\\x_2\in\{2,5\}}}
\sum_{c=1}^2\omega_c
\frac{\sum_{j=1}^2x_j\exp(a_{cj}-b_{cj}x_j)}
{1+\sum_{j=1}^2\exp(a_{cj}-b_{cj}x_j)},
\]
where
$
(\omega_1,\omega_2)=\left(\frac25,\frac35\right)$,
$a=
\begin{pmatrix}
3\log(\frac{3}{2})&\frac23\log 3\\
3\log 3&\frac23\log 6
\end{pmatrix}$,
and
$b=
\begin{pmatrix}
2\log(\frac{3}{2})&\frac13\log 3\\
2\log 3&\frac13\log 6
\end{pmatrix}$.
Let \(z_j=1\) indicate that the higher price is selected for product \(j\). Substituting, the problem reduces to:
\[
\max_{z_1,z_2\in\{0,1\}}
\eta(z_1,z_2) := \frac{\frac75-\frac1{15}z_1-\frac2{15}z_2}
{3-\frac13z_1-\frac23z_2}
+
\frac{\frac{21}{10}-\frac12z_1-\frac7{10}z_2}
{3-\frac23z_1-\frac56z_2}.
\]
Figure~\ref{scip_pd_projection} compares the projections of SCIP's root-LP relaxation with default cuts enabled, the relaxation~\eqref{fp}, and the continuous relaxation of~\eqref{eq:MIP} onto the $(z_1,z_2,\eta)$-space. Their maximum objective values are approximately $1.290741$, $104/85$, and $6/5$, respectively. For this instance, the relaxation of~\eqref{fp} yields a tighter upper bound than SCIP's root LP, while the relaxation of~\eqref{eq:MIP} yields a strictly tighter bound than the relaxation of~\eqref{fp}, demonstrating that the dominance established in Proposition~\ref{prop:PD-dominates-FP} can be strict. \Halmos
\begin{figure}
    \centering
    \includegraphics[scale=0.4]{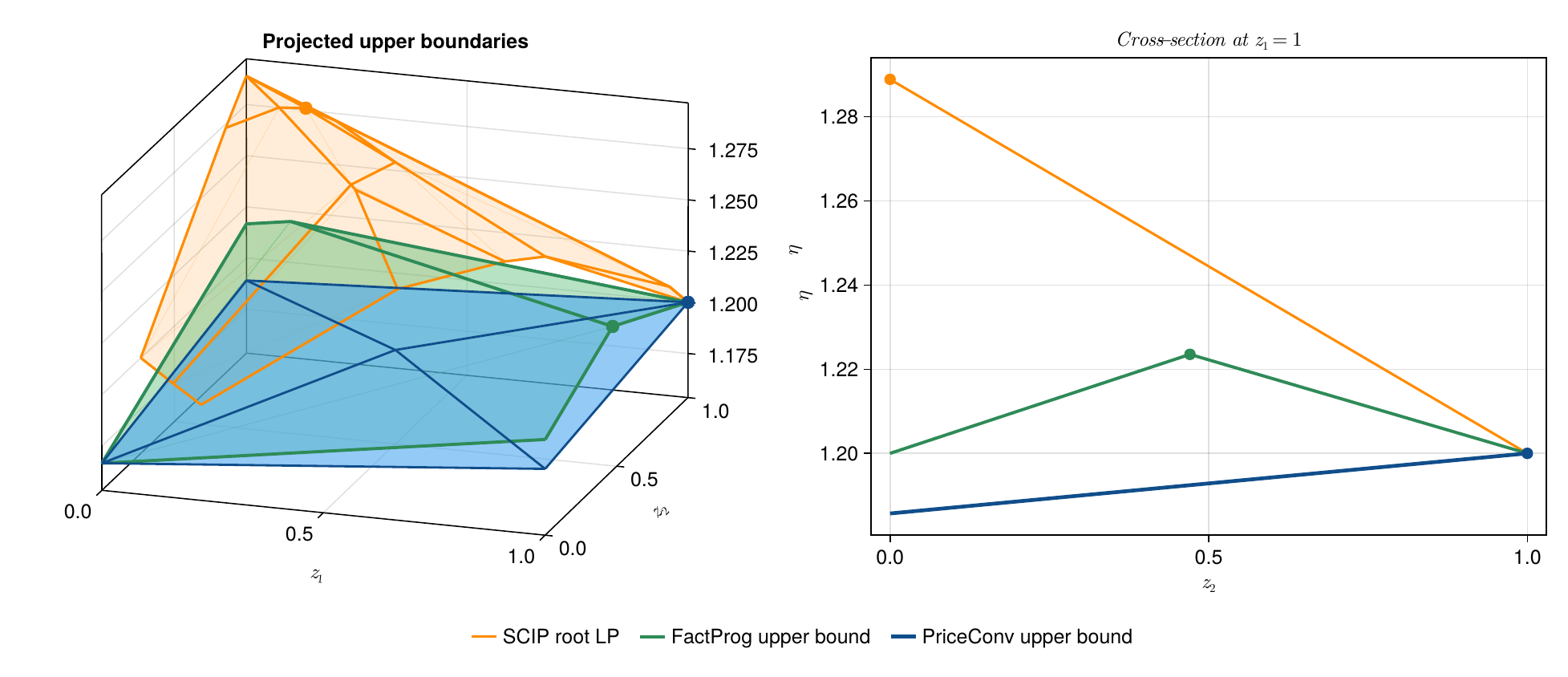}
    \caption{Comparison of the projected upper boundaries of the SCIP root-LP relaxation(orange), the relaxation of~\eqref{fp} (green), and the relaxation of~\eqref{eq:MIP} (blue) in the \((z_1,z_2,\eta)\)-space. The left panel displays the complete upper-boundary surfaces, while the right panel shows the cross-section at \(z_1=1\).}
    \label{scip_pd_projection}
\end{figure}
\end{example}
To relax fractional programs, the literature has used McCormick envelopes to replace the linking equalities \(y_{jr}^c=\rho^cz_{jr}\)~\citep{mehmanchi2019fractional,tawarmalani2002global}.
Let $\underline \rho^c$ and $\overline \rho^c$ be the lower and upper bounds of $1/D^c(\z)$ over $z\in \vertex(\Delta)$ respectively. Since $\exp(a_{cj}-b_{cj}x_j)$ is decreasing in $x_j$, we have
\[
\underline \rho^c=\frac{1}{1+\sum_{j\in N}\exp(a_{cj}-b_{cj}p_{j1})},\qquad
\overline \rho^c=\frac{1}{1+\sum_{j\in N}\exp(a_{cj}-b_{cj}p_{jd_j})},\qquad c\in M.
\]
The McCormick inequalities for the product \(z_{jk}\rho^c\) are
\[
\underline\rho^c z_{jk}\leq y^c_{jk}\leq\overline\rho^c z_{jk},\quad
\rho^c-\overline\rho^c(1-z_{jk})\leq y^c_{jk}\leq\rho^c-\underline\rho^c(1-z_{jk}),\qquad
c\in M,\ j\in N,\ k\in[d_j-1]\label{mcc}
    \tag{\textsc{McC}}
\]
We refer to the formulation obtained by replacing constraints~\eqref{eq:coordinating-1} and~\eqref{eq:coordinating-2} with constraints~\eqref{mcc} as the \textsc{FracProg} formulation, which is formally written as follows:
\begin{equation}
\begin{aligned}\label{fracprog}
\max \quad
& \sum_{c \in M} \omega_c \cdot \sum_{i\in N} \biggl(\gamma^c_{i1} \cdot \rho^c +\sum_{k=2}^{d_i}\gamma^c_{ik} \cdot y^c_{ik-1}\biggr) \\
\text{\rm s.t.}\quad
&\eqref{eq:mip-1},\eqref{eq:mip-2},\eqref{eq:unary},\eqref{mcc}.
\end{aligned} \tag{\textsc{FracProg}}
\end{equation}
Constraints~\eqref{eq:coordinating-1} and~\eqref{eq:coordinating-2} are strictly stronger than the McCormick inequalities~\eqref{mcc}. The following result formalizes the dominance relation between the two relaxations.
\begin{proposition}
\label{prop:PD-dominates-MCPD}
The LP relaxation of~\eqref{eq:MIP} is tighter than that of~\eqref{fracprog}.
\end{proposition}
Proposition~\ref{prop:PD-dominates-MCPD} follows from Theorem~\ref{them:indhull} or the proof of Theorem~\ref{them:exactnessMIP}. To see the first derivation, observe that \eqref{mcc} is implied in the convex hull of $(z_{jk},\rho^c, y^c_{jk})$. In fact, we will show an even stronger result in Proposition~\ref{cpd better than mcc}. Specifically, we will retain only a subset of coordinating inequalities~\eqref{eq:coordinating} and show that the resulting formulation dominates \eqref{fracprog}.

\begin{example}
Consider the setting of Example~\ref{ex:dominateFactProg}:
\[
\max_{z_1,z_2\in\{0,1\}}
\frac{\frac75-\frac1{15}z_1-\frac2{15}z_2}
{3-\frac13z_1-\frac23z_2}
+
\frac{\frac{21}{10}-\frac12z_1-\frac7{10}z_2}
{3-\frac23z_1-\frac56z_2}.
\]
An optimal solution of the relaxation of~\eqref{fracprog} is
\[
z_1=z_2=\frac23,\qquad
\rho_1=\rho_2=\frac49,\qquad
(y_1^1,y_2^1)=\left(\frac13,\frac13\right),\qquad
(y_1^2,y_2^2)=\left(\frac29,\frac29\right).
\]
However, this solution violates the upper coordinating constraint of~\eqref{eq:MIP} associated with the second fractional term and \(z_2\), which requires
\[
z_2\leq\frac{13}{6}y_2^2-\frac23\max\{0,y_1^2+y_2^2-\rho^2\}.
\]
Hence, the optimal solution of the~\eqref{fracprog} relaxation can be separated from the relaxation of the feasible region of~\eqref{eq:MIP}. Indeed the relaxation of~\eqref{eq:MIP} yields an upper bound of $6/5$ which is tighter than the upper bound of $11/9$ from the relaxation of~\eqref{fracprog}. The feasible regions of the two relaxations projected to $(z_1,z_2,\eta)$ are compared in Figure~\ref{relaxation_projection}. \Halmos
\begin{figure}
    \centering
    \includegraphics[scale=0.4]{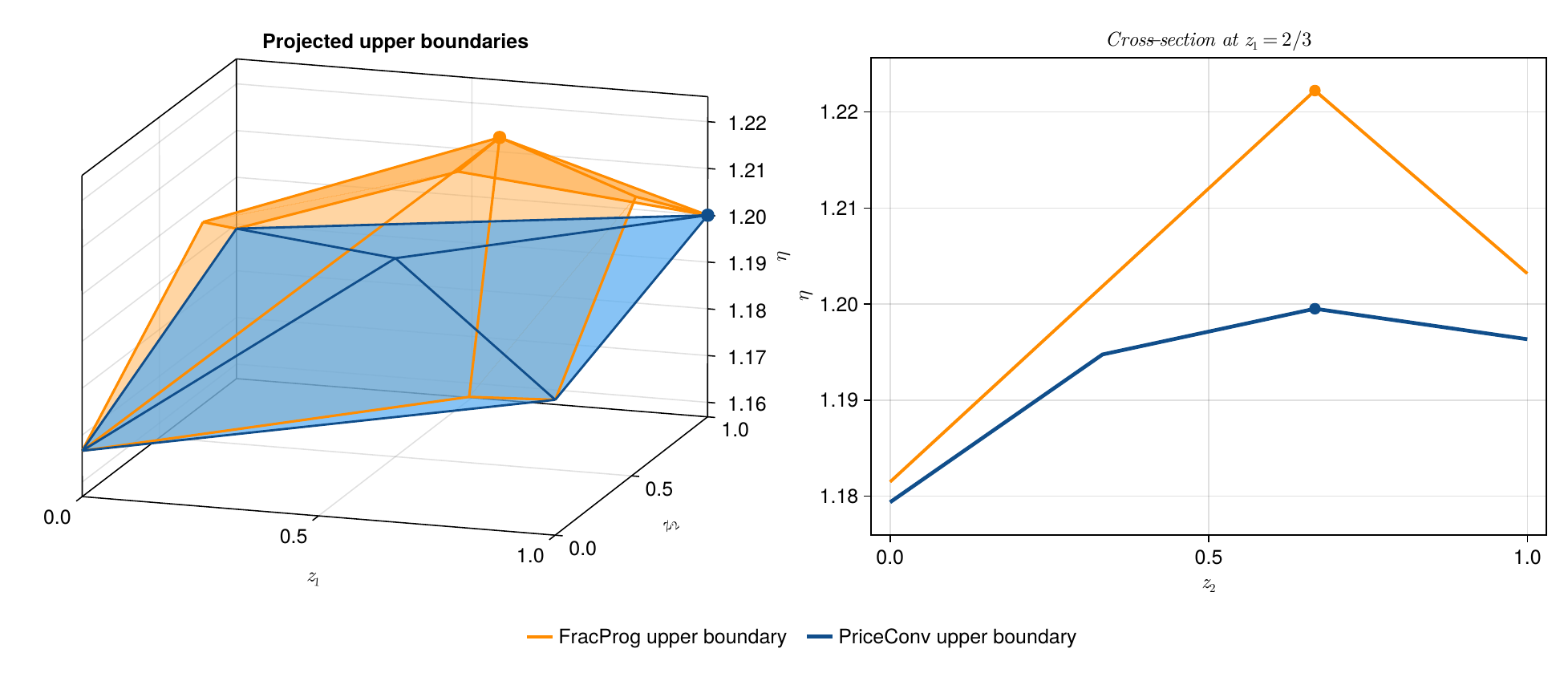}
    \caption{Comparison of the projected upper boundaries of the relaxation of~\eqref{fracprog} (orange), and~\eqref{eq:MIP} (blue) in the $(z_1,z_2,\eta)$-space.
    The left panel displays the complete upper-boundary surfaces, while the right panel shows the cross-section at \(z_1=2/3\).}
    \label{relaxation_projection}
\end{figure}
\end{example}

\subsection{An exact LP formulation via higher order pricing coordination}\label{section:LP}
In Section~\ref{section:MIP}, we derive the normalization equality in~\eqref{eq:mip-1} and the coordinating inequalities in~\eqref{eq:coordinating} through a two-step relaxation  of the following nonlinear relations:
\[
    1 = 1 \cdot D^c(\z) \cdot \rho^c \quad \text{ and } \quad
    z_{ik} = z_{ik} \cdot D^c(\z) \cdot \rho^c,
\]
respectively. Accordingly, we refer to~\eqref{eq:mip-1} as the zero-order pricing coordination and to~\eqref{eq:coordinating} as the first-order pricing coordination. This section extends such an approach to higher-order coordination involving products of multiple unary variables. 

To define these higher-order interactions, we first introduce some notation. Let $\mathcal{I} := \bigl\{ (i,k) \bigm| i \in N,\ k \in [d_i-1] \bigr\}$ denote the index set of the unary variables in $\z$. We define $\mathcal{A}$ as the family of all subsets of $\mathcal{I}$ that contain at most one pair associated with each product, that is,
\[
\mathcal{A} : = \Bigl\{ \alpha \subseteq \mathcal{I} \Bigm| \bigl\vert \alpha \cap (\{i\} \times [d_i-1] ) \bigr\vert \leq 1 \, \for i \in N \Bigr\}. 
\]
In other words, no two distinct pairs from $\mathcal{A}$ share the same first component.  For each $\alpha\in\mathcal{A}$, we define the corresponding monomial in the unary variables and its normalization with respect to customer segment $c$ as
\begin{equation}\label{eq:high-order}
    Z_\alpha = \prod_{(i,k) \in \alpha} z_{ik} \quad  \text{ and } \quad W^c_\alpha = \frac{\prod_{(i,k) \in \alpha}z_{ik}}{D^c(\z)},
\end{equation}
where $Z_{\emptyset} := 1$ and $W^c_{\emptyset} := \rho^c$. These definitions yield the following higher-order nonlinear identity
\begin{equation}\label{eq:high-order-identity}
Z_\alpha = \prod_{(i,k) \in \alpha} z_{ik} \cdot D^c(\z) \cdot \rho^c.     
\end{equation}
Note that $\vert \alpha \vert$ represents the degree of a monomial $Z_\alpha$. 

Next, we present a two-step relaxation for these nonlinear identities. The first step tackles the monomials in the unary variables. For an integer $\tau$ less than the number of items $n$, we define the set of monomials with degree up to $\tau$ as follows:
\[
\mathcal{M}^{\tau} := \biggl\{ \boldsymbol{Z} \biggm| \z_i \in \Delta^{d_i-1} \cap \{0,1\}^{d_i-1} \for i \in N,\  Z_\alpha = \prod_{(i,k) \in \alpha}z_{ik} \for \alpha \in  \mathcal{A} \text{ s.t. }  \vert \alpha \vert \leq \tau \biggr\}.
\]
Given a polyhedral relaxation of $\mathcal{M}^{\tau+1}$, we homogenize the resulting linear system and obtain a polyhedral relaxation of the nonlinear identities defined as in~\eqref{eq:high-order-identity}, for monomials of degree at most $\tau$.

To formally describe this relaxation, we introduce notation for the monomials arising from the product $Z_\alpha \cdot D^c(\z)$. Specifically, for each $\alpha\in\mathcal A$ and $(i,k)\in\mathcal I$, define
\[
\alpha\vee(i,k):=
\begin{cases}
\alpha \cup \bigl\{(i,k) \bigr\}
& \text{if } \alpha \cap (\{i\}\times[d_i-1])=\emptyset \\
\alpha \setminus\bigl\{(i,r) \bigr\}\cup \bigl\{(i,\max\{k,r\}) \bigr\} 
& \text{if there exists } r \text{ s.t. } (i,r) \in \alpha.
\end{cases}
\]
This operation accounts for the unary identity $z_{ik}z_{is}=z_{i,\max\{k,s\}}$ when $\alpha$ already contains a unary variable associated with item $i$. Now, for each segment $c\in M$ and degree $\tau$, we define the relaxation of the nonlinear identities involving monomials of degree at most $\tau$ as the set of points $(\boldsymbol Z,\boldsymbol W^c)$ satisfying
\begin{subequations}\label{eq:highorder}
\begin{align}
    & \boldsymbol{W}^c\in W^c_{\emptyset} \cdot \conv(\mathcal{M}^{\tau +1 }) \label{eq:highorderScaleF}\\
    &Z_\alpha = \biggl(1+\sum_{j\in N}\delta^c_{j1}\biggr)W^c_\alpha
+\sum_{j\in N}\sum_{r=2}^{d_j}
\delta^c_{jr}W^c_{\alpha\vee(j,r-1)} \quad \for  \alpha\in\mathcal A \text{ s.t. } \vert \alpha \vert \leq \tau \label{eq:highorderNormalize} \\
& Z_{\emptyset} = 1 \text{ and } W^c_{\emptyset}  \geq 0. \label{eq:highorderNonnegative}
\end{align}
\end{subequations}
The formulations~\eqref{eq:highorder} and \eqref{eq:onehot-convexhull} have structural similarities and are built from the same three pieces. First, a cone tying the normalized product variables to the zero-order variable ($W^c_\emptyset$ in \eqref{eq:highorder} and $\rho^c$ in \eqref{eq:onehot-convexhull}): \eqref{eq:highorderScaleF} places $\boldsymbol{W}^c$ in a cone obtained by scaling $\conv(\mathcal{M}^{\tau+1})$ with $W^c_\emptyset$, just as \eqref{eq:onehot-beginflow}--\eqref{eq:onehot-nonnegative} define a cone in $(\rho^c, \s^c, \boldsymbol{\Gamma}^c)$. Second, a normalizing plane fixes the zero-order variable: \eqref{eq:highorderNormalize} at $\alpha = \emptyset$, and \eqref{eq:onehot-normalize}. Third, an affine linking expresses the remaining variables using the normalized ones: \eqref{eq:highorderNormalize} for $Z_\alpha$, and \eqref{eq:onehot-defn} for $\lambda_{ik}$.

The formulations differ in two respects. First,  \eqref{eq:highorder} is based on  the unary binarization while \eqref{eq:onehot-convexhull} is written in the full binarization. Second, \eqref{eq:highorder} allows arbitrary interactions of order up to $n$, whereas \eqref{eq:onehot-convexhull} captures only pairwise interactions. We address the first difference by the linear transformation $L$ defined as in~\eqref{eq:Ltrans}: the flow variables $\Gamma^c_\alpha$ equal the normalized one-hot products, $\Gamma^c_\alpha = \rho^c\prod_{(j,l_j)\in \alpha}\lambda_{jl_j}$, with the lower-order sums in \eqref{eq:onehot-convexhull} obtained by summing over the slices of the unchosen items; hence, for two products, the flow variables coincide with $\boldsymbol{s}^c$ and $\boldsymbol{\Gamma}^c$ variables of \eqref{eq:onehot-convexhull}. The $L$ defined in~\eqref{eq:Ltrans} provides a bijection between $\boldsymbol{\Gamma}^c$ and $\boldsymbol{W}^c$ variables. Therefore, \eqref{eq:highorder} and \eqref{eq:onehot-convexhull} arise from the same template in two linearly related bases. The second difference is addressed by extending the individual hull characterization in Proposition~\ref{prop:onehot-convexhull} to a multi-segment, multi-product convex hull characterization in Theorem~\ref{hierarchy}.

To derive explicit polyhedral relaxations of constraint~\eqref{eq:highorderScaleF}, one-hot variables provide a more natural representation, although, as demonstrated in Section~\ref{section:full_binarization}, unary variables offer an advantage in describing the individual hull. In particular, polyhedral descriptions based on one-hot variables admit an interpretation in terms of network flows. To make this connection explicit, we introduce a directed multi-partite graph $G$ with $n$ intermediate layers. Each layer contains $d_i$ nodes.  Every node in layer $i$ is connected to every node in layer $i+1$, while no edges connect nonconsecutive layers. The source is connected to every node in the first layer, and every node in the final layer is connected to the sink. Formally, $G= ( V, E)$ such that
\[
V_i:=\{i\}\times[d_i],\qquad
 V:=\{s,t\}\cup
\bigcup_{i\in N}  V_{i},
\]
and
\[
 E:=
\bigl(\{s\}\times  V_{1}\bigr)\cup
\bigcup_{i=1}^{n-1}
\bigl( V_{i}\times  V_{i+1}\bigr)\cup
\bigl( V_{n}\times\{t\}\bigr).
\]
An example with three layers is illustrated in Figure~\ref{fig:multipartite-transport}.  Each node $(i,k)$ in layer $i$  corresponds to the one-hot variable $\lambda_{ik}$ of item $i$ at a certain price level $k$. Each $s$--$t$ path, specified by the sequence of intermediate nodes $(1,k_1), \ldots, (n,k_n)$, corresponds to a degree $n$ interaction $ \prod_{i \in N} \lambda_{ik_i}$. Moreover, the same path is associated with its lower-degree interactions $\prod_{i \in I} \lambda_{ik_i}$ for $I \subseteq N$, obtained by restricting the path to a subset of its intermediate layers. Let $\bar{\mathcal A}$ denote the collection of all $s$--$t$ paths and their subpaths, that is
\[
\bar{\mathcal{A}} := \Bigl\{\bigl\{(i,k_i)\bigr\}_{i \in I} \Bigm| I \subseteq N,\ k_i \in [d_i] \for i \in I \Bigr\}.
\]
For each customer segment $c \in M$, we associate a flow polytope with the graph $G$, with total flow $W^c_{\emptyset}$, which equals $\rho^c$ and represents the probability that a customer in segment $c$ chooses the outside option.  For each partial path $\alpha \in \bar{\mathcal{A}}$, let $\Gamma^c_\alpha$ denote the aggregate flow over all $s$--$t$ paths containing the nodes specified by $\alpha$. The resulting polytope is characterized by nonnegativity, a prescribed total flow, and flow aggregation equations.
More specifically, for each segment $c$ and each interaction order $\tau \geq 2 $, we impose  flow conservation equations for partial paths of length strictly less than  $\tau$, yielding the following polytope:  
\begin{equation}\label{eq:multipartitePolytope}
    \begin{aligned}
        P^\tau(W^c_\emptyset):=
    \biggl\{ &
    \boldsymbol{\Gamma}^c \in \mathbb R_+^{| \bar{\mathcal A}|}  \biggm| 
    \Gamma^c_{\emptyset} = W^c_\emptyset, \\
& \quad \Gamma_\alpha^c=\sum_{\ell=1}^{d_i}\Gamma_{\alpha\cup\{(i,\ell)\}}^c 
    \for \alpha\in\mathcal A,\ 
     i\in N\text{ s.t. } \vert \alpha \vert < \tau \text{ and } \alpha\cap(\{i\}\times[d_i])=\emptyset
    \biggr\}.
    \end{aligned}
\end{equation}
Note that $P^n(W^c_\emptyset) = P^{n+1}(W^c_\emptyset)$. We next introduce an invertible linear transformation $L$ that maps the flow variables $\boldsymbol{\Gamma}$ to the unary monomial variables $\boldsymbol{W}$. Namely, for each  $\alpha = \{(i,k_i): i \in I\} \in \mathcal{A}$, define
\begin{equation}\label{eq:Ltrans}
    W_\alpha = \sum_{\bar{\alpha} \in \{(i,\ell_i) \mid i \in I,k_i< \ell_i \leq d_i \}} \Gamma_{\bar{\alpha}}.
        \end{equation}
Conversely, for each
$\bar{\alpha}=\{(i,\ell_i):i\in I\}\in\bar{\mathcal A}$,
the inverse transformation is given by
\begin{equation}\label{eq:Linverse}
\begin{aligned}
\Gamma_{\bar{\alpha}}
=\sum_{S\subseteq I}(-1)^{|S|}
W_{\{(i,\ell_i):i\in S\}\cup
\{(i,\ell_i-1):i\in I\setminus S\}},
\end{aligned}
\end{equation}
where, for convenience,  $W_{\{(i,0)\}\cup\alpha}:=W_\alpha$ and $W_{\{(i,d_i)\}\cup\alpha}:=0$.

\begin{figure}[t]
\tikzsetnextfilename{priceopt-2}
\centering
\tikzsetnextfilename{pricelogit-1}
\begin{tikzpicture}[x=1cm, y=1cm,
  edge/.style={thin, black!80},
  faceA/.style={fill=black!5},
  faceB/.style={fill=black!18},
  faceC/.style={fill=black!14},
  faceD/.style={fill=black!9},
  dot/.style={circle, fill=black, draw=white, line width=0.4pt, minimum size=3pt, inner sep=0pt},
  arr/.style={->, thick, black!75},
  lbl/.style={font=\tiny},
  ttl/.style={font=\tiny, align=center},
  sub/.style={font=\tiny, align=center, black!70},
  box4/.style={draw=black!40, fill=black!5, rounded corners=3pt, font=\tiny, align=center, inner sep=8pt},
]
\begin{scope}[shift={(0.2,1.0)}]
\node (O1) at (0,0) {};
\node (A1) at (1,0) {};
\node (B1) at (0.3,0.6) {};
\node (C1) at (1.3,1.6) {};
\fill[faceA] (O1) -- (A1) -- (C1) -- cycle;
\fill[faceB] (O1) -- (A1) -- (B1) -- cycle;
\fill[faceC] (O1) -- (B1) -- (C1) -- cycle;
\fill[faceD] (A1) -- (B1) -- (C1) -- cycle;
\draw[edge] (O1) -- (A1);
\draw[edge] (A1) -- (C1);
\draw[edge] (C1) -- (B1);
\draw[edge] (B1) -- (O1);
\draw[edge] (O1) -- (C1);
\draw[edge, dashed] (A1) -- (B1);
\node[dot] at (0,0) {};
\node[dot] at (1,0) {};
\node[dot] at (0.3,0.6) {};
\node[dot] at (1.3,1.6) {};
\end{scope}
\node[ttl, text width=2.6cm] at (0.85,3.15) {Scaled multipartite partial flow polytope};
\node[sub, text width=2.6cm] at (0.85,0.55) {$P^{\tau+1}( W^c_\emptyset)$};
\node[sub, text width=2.6cm] at (0.85,0.3) {one-hot variables $\boldsymbol{\Gamma}^c$};

\begin{scope}[shift={(3.0,1.0)}]
\node (O2) at (0,0) {};
\node (A2) at (1.1,0.25) {};
\node (B2) at (0.3,0.6) {};
\node (C2) at (1.3,1.6) {};
\fill[faceA] (O2) -- (A2) -- (C2) -- cycle;
\fill[faceB] (O2) -- (A2) -- (B2) -- cycle;
\fill[faceC] (O2) -- (B2) -- (C2) -- cycle;
\fill[faceD] (A2) -- (B2) -- (C2) -- cycle;
\draw[edge] (O2) -- (A2);
\draw[edge] (A2) -- (C2);
\draw[edge] (C2) -- (B2);
\draw[edge] (B2) -- (O2);
\draw[edge] (O2) -- (C2);
\draw[edge, dashed] (A2) -- (B2);
\node[dot] at (0,0) {};
\node[dot] at (A2) {};
\node[dot] at (0.3,0.6) {};
\node[dot] at (1.3,1.6) {};
\end{scope}
\node[ttl, text width=4cm] at (3.9,3.15) {Linearly transformed polytope};
\node[sub, text width=2.6cm] at (3.9,0.55) {Relax $W^c_{\emptyset} \cdot \mathcal{M}^{\tau+1}$};
\node[sub, text width=2.6cm] at (3.9,0.3) {unary variables $\boldsymbol{W}$};

\begin{scope}[shift={(7.6,1.3)}]
\node (O3) at (0,0) {};
\node (A3) at (0.5,-0.2) {};
\node (B3) at (0.3,0.6) {};
\node (C3) at (0.8,1.267) {};
\fill[faceA] (O3) -- (A3) -- (C3) -- cycle;
\fill[faceB] (O3) -- (A3) -- (B3) -- cycle;
\fill[faceC] (O3) -- (B3) -- (C3) -- cycle;
\fill[faceD] (A3) -- (B3) -- (C3) -- cycle;
\draw[edge] (O3) -- (A3);
\draw[edge] (A3) -- (C3);
\draw[edge] (C3) -- (B3);
\draw[edge] (B3) -- (O3);
\draw[edge] (O3) -- (C3);
\draw[edge, dashed] (A3) -- (B3);
\node[dot] at (0,0) {};
\node[dot] at (A3) {};
\node[dot] at (0.3,0.6) {};
\node[dot] at (0.8,1.267) {};
\end{scope}
\node[ttl, text width=2.6cm] at (7.9,3.15) {Projective transformation};
\node[sub, text width=2.6cm] at (7.9,0.55) {Relax the high order interactions~\eqref{eq:high-order}};
\node[box4, text width=3cm] at (12.6,1.6) {Theorem \ref{hierarchy}: exact convex hull (ideal) formulation at full order};
\draw[arr] (1.5,1.6) -- (3,1.6);
\node[lbl, text width=2cm, align=center] at (2.45,2.20) {mapping via $L$};
\node[lbl, text width=2cm, align=center] at (2.45,1.95) {defined as in \eqref{eq:Ltrans}};
\draw[arr] (4.4,1.6) -- (7.6,1.6);
\node[lbl] at (5.95,2.65) {Projective normalization};
\node[lbl] at (5.95,2.35) {$1 = (1+\sum_{j\in N} \delta^c_{j1})W^c_\emptyset$};
\node[lbl] at (5.95,1.95) {$+ \sum_{j\in N}\sum_{r=2}^{d_j} \delta^c_{jr}W^c_{(j,r-1)}$};
\draw[arr] (8.5,1.6) -- (10.6,1.6);
\node[lbl] at (9.5,1.95) {full order (Theorem \ref{hierarchy})};
\node[draw=black!25, fill=black!4, rounded corners=3pt, align=left, font=\tiny, text width=14cm, inner sep=6pt] at (6.8,-1)
{\textbf{Pricing interpretation.}\quad
$\Gamma^c_\emptyset = W^c_\emptyset = \rho^c$: no-purchase probability (total flow used for scaling multipartite flow polytope).\\[0.5em]
$\Gamma^c_\alpha$: scaled mass on price combinations that extend the partial assignment $\alpha$ (a full path $=$ a complete price vector).\\[0.5em]
$W^c_\alpha$: normalized mass on price combinations that exceed the unary threshold $\alpha$, related to the flows by the invertible $L$.};
\end{tikzpicture}\caption{Theorem~\ref{hierarchy} illustration. (1) The polytope $P^\tau(W^c_\emptyset)$
is the scaled multipartite partial flow polytope in \eqref{eq:multipartitePolytope}. (2) The linear transformation $L$ in~\eqref{eq:Ltrans} maps $P^{\tau+1}( W^c_\emptyset)$ to $W^c_\emptyset \conv(\mathcal{M}^{\tau+1}  )$. (3) The projective map rescaling yields a lifted convex relaxation of the interactions in~\eqref{eq:high-order} with an order of $\tau$. (4) At level $\tau =n$, the relaxation hierarchy terminates with the optimal value for \eqref{price optimization}, as shown in Theorem~\ref{hierarchy}.}
\label{fig:thm4-engine}
\end{figure}
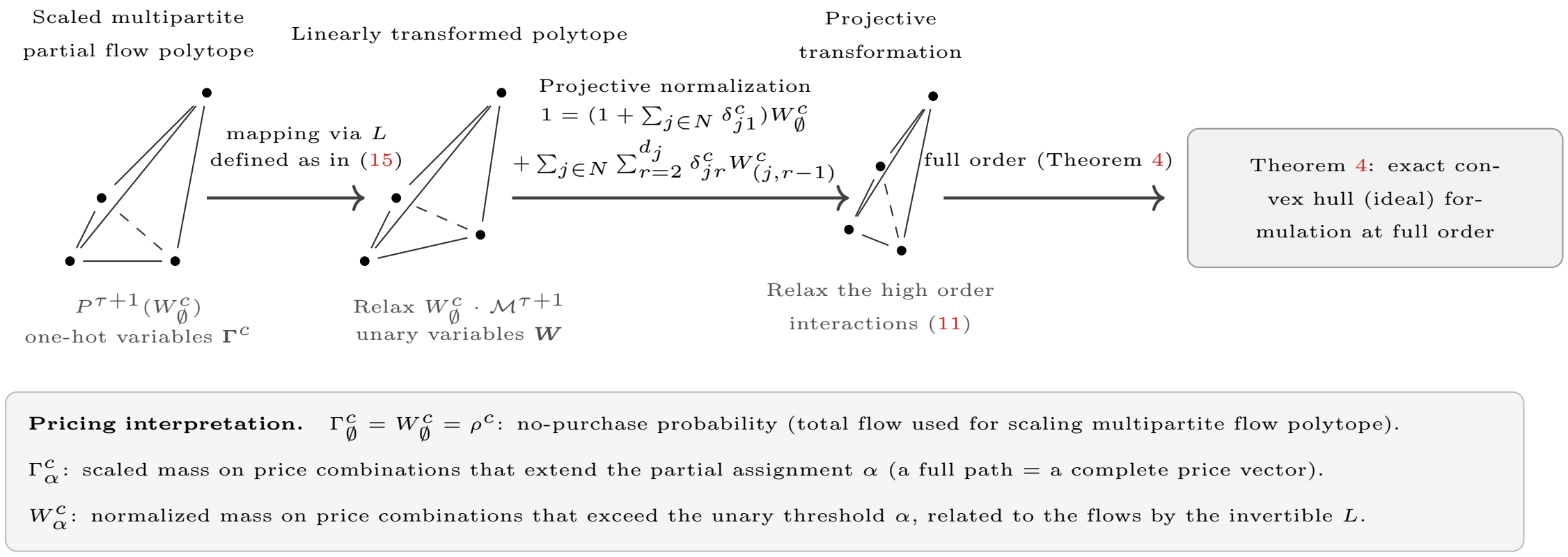

With the flow polytopes $P^\tau(W^c_{\emptyset})$ and the linear transformation $L$, we obtain an explicit polyhedral description for the constraint~\eqref{eq:highorderScaleF}, obtaining a hierarchy of LP formulations. In Figure~\ref{fig:thm4-engine}, we visualize the whole construction of our relaxation hierarchy. 
\begin{figure}[t]
\centering
\tikzsetnextfilename{priceopt-1}
\begin{tikzpicture}[
    >=Stealth,
    vertex/.style={circle,draw,minimum size=8mm,inner sep=1pt},
    terminal/.style={circle,draw,minimum size=9mm},
    arc/.style={->,thin,gray!65}
]

\node[terminal] (source) at (0,0) {$s$};
\node[vertex] (v11) at (2.5,1.6) {$(1,1)$};
\node[vertex] (v12) at (2.5,0.4) {$(1,2)$};
\node at (2.5,-0.45) {$\vdots$};
\node[vertex] (v1d) at (2.5,-1.6) {$(1,d_1)$};

\node[vertex] (v21) at (5.5,1.6) {$(2,1)$};
\node[vertex] (v22) at (5.5,0.4) {$(2,2)$};
\node at (5.5,-0.45) {$\vdots$};
\node[vertex] (v2d) at (5.5,-1.6) {$(2,d_2)$};

\node[vertex] (v31) at (8.5,1.6) {$(3,1)$};
\node[vertex] (v32) at (8.5,0.4) {$(3,2)$};
\node at (8.5,-0.45) {$\vdots$};
\node[vertex] (v3d) at (8.5,-1.6) {$(3,d_3)$};

\node[terminal] (sink) at (11,0) {$t$};

\foreach \v in {v11,v12,v1d}
    \draw[arc] (source) -- (\v);

\foreach \u in {v11,v12,v1d}
    \foreach \v in {v21,v22,v2d}
        \draw[arc] (\u) -- (\v);

\foreach \u in {v21,v22,v2d}
    \foreach \v in {v31,v32,v3d}
        \draw[arc] (\u) -- (\v);

\foreach \u in {v31,v32,v3d}
    \draw[arc] (\u) -- (sink);
    \foreach \u/\v in {source/v11, v11/v21, v21/v31, v31/sink}
    \draw[arc,red,thick] (\u) -- (\v);

\node[above] at (2.5,2.15) {Layer \(1\)};
\node[above] at (5.5,2.15) {Layer \(2\)};
\node[above] at (8.5,2.15) {Layer \(3\)};
\end{tikzpicture}
\caption{Layered transportation graph, illustrated for three products.}
\label{fig:multipartite-transport}
\end{figure}
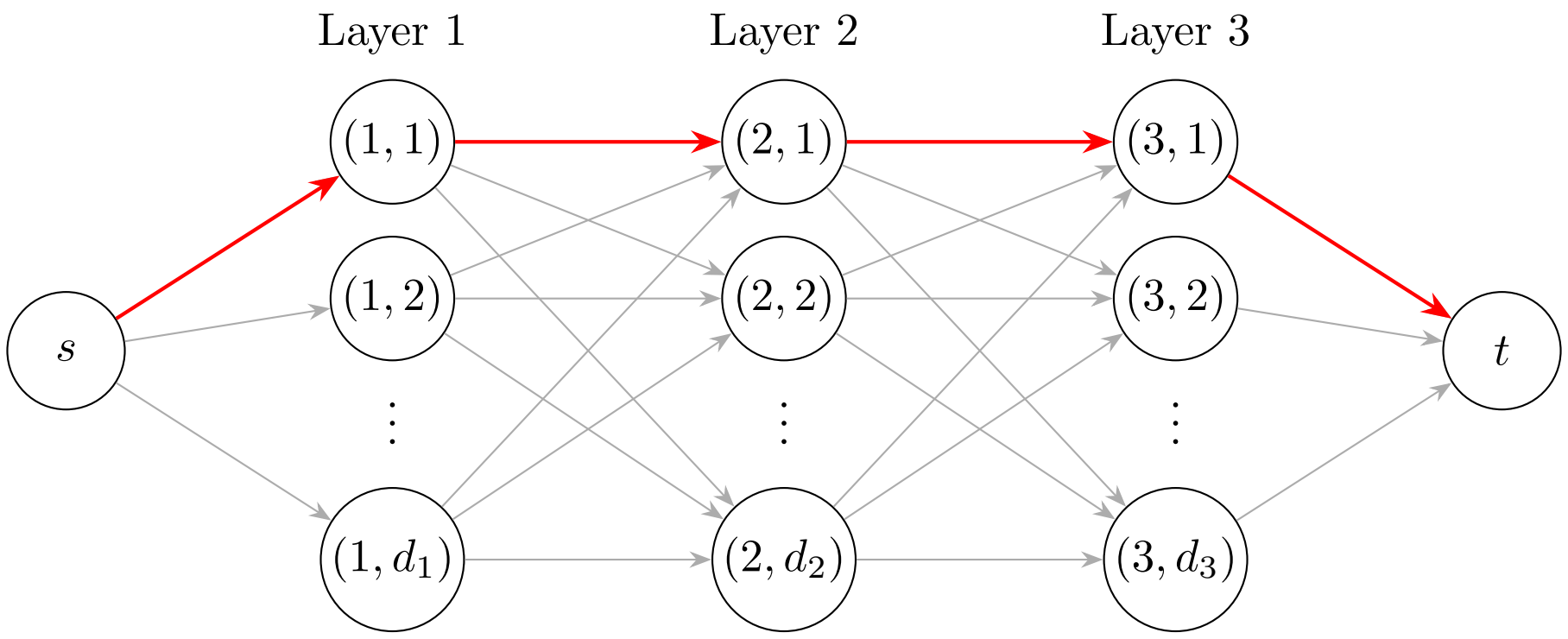

\begin{theorem}\label{hierarchy}
For $\tau \geq 2$, consider 
\begin{equation*}
\label{k-order hierarchy}
\begin{aligned}
R^*_\tau:=\max_{\boldsymbol{Z},\{\boldsymbol{W}^c\}_{c\in M}}\quad
&\sum_{c\in M}\omega_c\left(W^c_{\emptyset}\sum_{j\in N}\gamma^c_{j1}+\sum_{j\in N}\sum_{k=2}^{d_j}\gamma^c_{jk}W_{\{(j,k-1)\}}^c\right)\\
\text{\rm s.t.}\quad
&L^{-1}(\boldsymbol{W}^c) \in P^{\tau+1}(W^c_{\emptyset}), ~\eqref{eq:highorderNormalize}, \text{ and } ~\eqref{eq:highorderNonnegative} 
\quad \for  c\in M.
\end{aligned}
\end{equation*}
Let $R^*$ denote the optimal profit of~\eqref{price optimization}. 
Then, $R_2^*\geq\cdots\geq R_n^* = R^*$. 
\end{theorem}

At full order, all partial-flow variables can also be eliminated using the $d^n$ full-path flows. Consequently, the hierarchy admits an equivalent formulation with $O(md^n)$ variables and constraints.
In contrast, a direct application of the generic RLT hierarchy proposed in~\cite{he2025convexification} introduces moment variables indexed by all subsets of these indicators. Without exploiting the structure of unary simplex, its full order formulation has \(O(m2^{n(d-1)})\) variables and \(O(m2^{n(d-1)})\) constraints.

\section{Large-scale Solution via Decomposition}\label{sec:decomposition}
The MIP formulation~\eqref{eq:MIP} is impractical to solve directly: each \eqref{eq:coordinating} inequality resolves $O(nd)$ max/min terms inducing exponentially many linear inequalities. The following relaxation retains only a small subset of the pricing coordinating inequalities in \eqref{eq:coordinating}:
\begin{subequations}\label{eq:compactpriceconv}
\begin{align}
\max \quad & \sum_{c \in M} \omega_c \cdot \sum_{i\in N} \biggl(\gamma^c_{i1} \cdot \rho^c +\sum_{k=2}^{d_i}\gamma^c_{ik} \cdot y^c_{ik-1}\biggr)\label{CPD}\tag{\textsc Compact} \\
\text{s.t.}\quad &\eqref{eq:mip-1},\eqref{eq:mip-2},\eqref{eq:unary},\notag\\
&z_{ik} \leq \biggl(1 + \sum_{j \in N}\delta^c_{j1} \biggr) y^c_{ik}
+ \sum_{r=1}^{d_i-1}\delta^c_{ir+1}y^c_{i,\max\{k,r\} }
+ \sum_{j \neq i} \sum_{r =1}^{d_j-1} \delta^c_{jr+1}\left(y^c_{ik}+ y^c_{jr}- \rho^c\right) \label{eq:coordinating-1-1}\\
&z_{ik} \leq \biggl(1 + \sum_{j \in N}\delta^c_{j1} \biggr) y^c_{ik} + \sum_{r=1}^{d_i-1}\delta^c_{ir+1}y^c_{i,\max\{k,r\} }\label{eq:coordinating-1-2}\\
&z_{ik} \geq \biggl(1 + \sum_{j \in N}\delta^c_{j1} \biggr) y^c_{ik} + \sum_{r=1}^{d_i-1}\delta^c_{ir+1}y^c_{i,\max\{k,r\} } + \sum_{j \neq i} \sum_{r =1}^{d_j-1} \delta^c_{jr+1}y^c_{ik}.  \label{eq:coordinating-2-1}\\
&z_{ik} \geq \biggl(1 + \sum_{j \in N}\delta^c_{j1} \biggr) y^c_{ik} + \sum_{r=1}^{d_i-1}\delta^c_{ir+1}y^c_{i,\max\{k,r\} } + \sum_{j \neq i} \sum_{r =1}^{d_j-1} \delta^c_{jr+1}y^c_{jr}.  \label{eq:coordinating-2-2}\\
&\qquad\qquad\qquad\qquad\qquad\qquad\qquad\quad\for c\in M, i \in N,  k \in [d_i-1].\notag
\end{align}
\end{subequations}
The formulation \eqref{CPD} is exact and, as we show next, dominates~\eqref{fracprog}.
\begin{proposition}\label{cpd better than mcc}
    Let $R_{\mathrm{CP}}$ and $R_{\mathrm{McC}}$ be the continuous relaxations of~\eqref{CPD} and~\eqref{fracprog}, respectively, obtained by relaxing the integrality constraints on $\z$.
    If $\delta^c_{jr}\leq 0$ for every $c\in M$,
    $j\in N$, $r\in\{2,\cdots,d_j\}$, and $1+\sum_{j\in N}\sum_{r\in [d_j]}\delta^c_{jr} > 0$ for every $c\in M$,
    then $R_{\mathrm{CP}}\subseteq R_{\mathrm{McC}}$. 
\end{proposition}
Our numerical experiments in Section~\ref{sec:neceofdec} show that \eqref{CPD} outperforms \eqref{fracprog}  for \eqref{price optimization}. However, as the instances grow larger or restrictions are imposed on price changes, \eqref{CPD} can struggle to solve the problem to optimality. We will now develop a tailored solution algorithm exploiting the structure of~\eqref{eq:MIP} that improves computational scalability.
\subsection{Decomposition Architecture}
We observe that the MILP formulation~\eqref{eq:MIP} lends itself to
decomposition techniques. Here, we develop Price Decomposition (PD), a
Benders-style, Euclidean-deepest-cut decomposition algorithm for optimizing a sum of fractional functions coupled by a common variable, with per-fraction subproblems that are linear programs; pricing is the first illustration of this scheme. To go beyond personalized pricing, the constraint \eqref{eq:coordinating} links the segment problems, requiring that $z_{ik}$ match across segments. Removing these constraints and $\z$ decomposes the formulation into independent segment blocks, resulting in the pricing LP that allows price discrimination (Theorem~\ref{them:personalized}). The algorithm progressively coordinates segment-specific price decisions by aggregating constraints \eqref{eq:mip-1}, \eqref{eq:mip-2}, \eqref{eq:coordinating}, and $\mu^c = \sum_i(\gamma_{ci1}\rho^c + \sum_k
\gamma_{cjk}\,y_{cjk})$ so that $\mu^c$ is represented correctly at binary
points.

To link the master and subproblems, for each segment \(c\in M\), we introduce \(\mu^c\) to represent the profit generated from the $c^{\text{th}}$ segment. Then, the linear constraints associated with segment \(c\) can be encapsulated into  the following polyderal set:
\begin{equation}
\label{subproblem}
\tag{\textsc{SegBlock}}
\mathcal P_c:=\left\{(\z,\mu^c,\y^c,\rho^c)\;\middle|\;
\begin{aligned}
&\mu^c=\sum_{i\in N} \biggl(\gamma^c_{i1} \cdot \rho^c +\sum_{k=2}^{d_i}\gamma^c_{ik} \cdot y^c_{ik-1}\biggr)\\
&\z_i \in \Delta^{d_i -1},~\eqref{eq:mip-1}, ~\eqref{eq:mip-2} \text{ and } \eqref{eq:coordinating} \for i \in N \text{ and }  k \in [d_i-1]
\end{aligned}
\right\}.
\end{equation}
Accordingly, formulation~\eqref{eq:MIP} can be expressed as a master problem involving only the unary variables and segment profit variables:
\begin{equation}\label{eq:masterBenders}
\max_{\z,\boldsymbol{\mu}}\left\{\sum_{c\in M}\omega_c\mu^c\;\middle|\;
(\z,\mu^c)\in\proj_{(\z,\mu^c)}(\mathcal P_c) \, \for c\in M \text{ and }
\z_i \in \{0,1\}^{d_i-1} \for i \in N\right\},
\end{equation}
where $\proj_{(\z,\mu^c)}(\mathcal P_c)$  projects the segment-specific constraints onto the  unary variables and profit variable.

The family of constraints obtained from eliminating the segment choice probability variable $(\rho^c,\y^c)$ is too large to explicitly enumerate. Thus, to solve the master problem, we iteratively generate the constraints by solving subproblems. More specifically, 
given a candidate master solution \((\bar \z,\bar{\boldsymbol{\mu}})\), each segment independently checks whether \((\bar \z,\bar\mu^c)\in\proj_{(\z,\mu^c)}(\mathcal P_c)\) by solving an LP with its segment-specific variables \((\rho^c,\y^c)\). If the candidate does not belong to \(\proj_{(\z,\mu^c)}(\mathcal P_c)\), the dual of this LP generates a valid inequality involving only \(\z\) and \(\mu^c\) that separates \((\bar \z,\bar\mu^c)\) from \(\proj_{(\z,\mu^c)}(\mathcal P_c)\). These inequalities are added to the master problem. Consequently, the master problem retains only the unary variables  and segment-level profits, while all other segment-specific variables and constraints remain in independent subproblems.
Since the segment-specific variables are disjoint across segments, \((\z,\boldsymbol\mu)\) admits a feasible extension to the full formulation if and only if \((\z,\mu^c)\in\proj_{(\z,\mu^c)}(\mathcal P_c)\) for every \(c\in M\). 

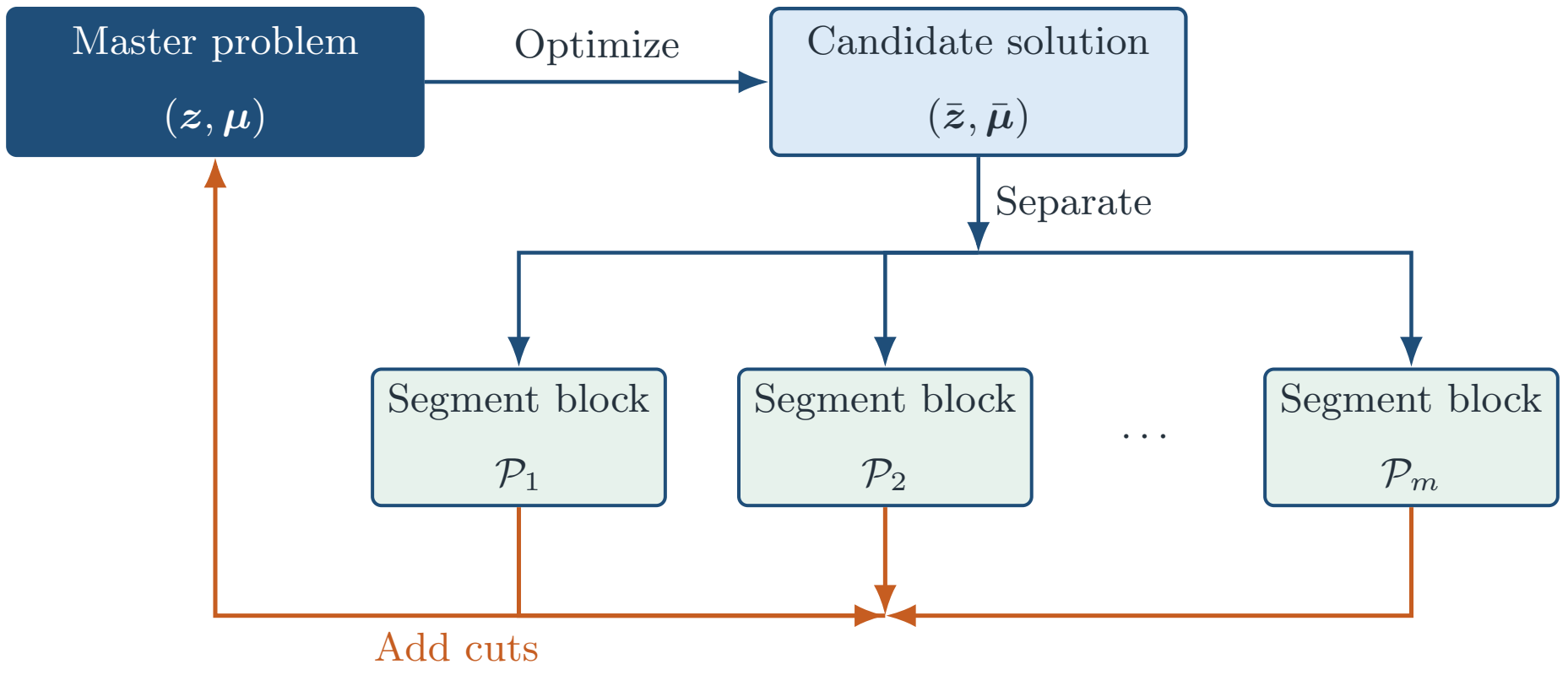
\begin{figure}
    \centering
\tikzsetnextfilename{priceopt-3}
\begin{tikzpicture}[
    font=\small,
    >={Latex[length=2.4mm]},
    master/.style={
        draw=masterblue,
        fill=masterblue,
        text=white,
        rounded corners=2pt,
        minimum width=3.25cm,
        minimum height=1.15cm,
        align=center,
        line width=0.8pt
    },
    candidate/.style={
        draw=masterblue,
        fill=candidateblue,
        text=textgray,
        rounded corners=2pt,
        minimum width=3.25cm,
        minimum height=1.15cm,
        align=center,
        line width=0.8pt
    },
    block/.style={
        draw=masterblue,
        fill=blockgreen,
        text=textgray,
        rounded corners=2pt,
        minimum width=2.25cm,
        minimum height=1.05cm,
        align=center,
        line width=0.75pt
    },
    arrow/.style={->, line width=0.85pt, draw=masterblue},
    cutarrow/.style={->, line width=0.95pt, draw=cutorange}
]

\node[master] (master) {Master problem\\$(\z,\boldsymbol\mu)$};
\node[candidate, right=2.7cm of master] (candidate)
    {Candidate solution\\$(\bar{\z},\bar{\boldsymbol{\mu}})$};

\draw[arrow] (master) -- node[above, text=textgray] {Optimize} (candidate);

\coordinate (dispatch) at ($(candidate.south)+(0,-0.75)$);
\node[block, below=1.65cm of candidate, xshift=-3.6cm] (p1)
    {Segment block\\$\mathcal P_1$};
\node[block, right=0.55cm of p1] (p2)
    {Segment block\\$\mathcal P_2$};
\node[right=0.56cm of p2, text=textgray] (dots) {$\cdots$};
\node[block, right=0.56cm of dots] (pm)
    {Segment block\\$\mathcal P_m$};

\draw[arrow] (candidate.south) -- node[right, text=textgray, align=left]
    {Separate} (dispatch);
\draw[arrow] (dispatch) -| (p1.north);
\draw[arrow] (dispatch) -| (p2.north);
\draw[arrow] (dispatch) -| (pm.north);

\coordinate (collector) at ($(p2.south)+(0,-0.85)$);
\draw[cutarrow] (p1.south) |- (collector);
\draw[cutarrow] (p2.south) -- (collector);
\draw[cutarrow] (pm.south) |- (collector);
\draw[cutarrow] (collector) -| node[below, pos=0.32, text=cutorange]
    {Add cuts} (master.south);
\end{tikzpicture}
\caption{Decomposition pipeline. The master problem selects a common price vector \(\bar \z\) and segment-level profits \(\bar{\boldsymbol{\mu}}\). For each customer segment \(c\), the corresponding segment block evaluates the profit achievable under the candidate common prices \(\bar{\z}\). If the master problem overestimates this profit, that is, if \(\bar\mu^c\) exceeds the value permitted by the segment's demand model, the segment block generates a cut that excludes the current candidate. Otherwise, the candidate is accepted for that segment. The procedure terminates when, for every segment \(c\), the profit \(\mu^c\) proposed by the master is attainable under the common price vector \(\z\).}
\label{fig:pipeline}
\end{figure}

Because the separations give outer approximations of $\proj_{(\z,\mu^c)}(\mathcal P_c)$, the optimal objective value of \eqref{eq:masterBenders} provides an upper bound on \eqref{price optimization}. Conversely, evaluating the candidate price vector \(\bar \z\) directly using~\eqref{price optimization} provides a feasible objective value and hence a lower bound. The master problem is iteratively reoptimized, and the procedure continues until the gap between the upper and lower bounds is within the prescribed tolerance. This procedure is visualized in Figure~\ref{fig:pipeline}, and is formally described in Algorithm~\ref{alg:pipeline}. We will develop a separation oracle to separate \((\bar \z,\bar\mu^c)\) from \(\proj_{(\z,\mu^c)}(\mathcal P_c)\) in the next subsection.
\begin{algorithm}[t]
\caption{Price Decomposition (PD)}
\label{alg:pipeline}
\begin{algorithmic}[1]
\Require Initial solution $(\bar \z,\bar{\boldsymbol{\mu}})$; Initial feasible region $\mathcal C = \{(\z,\boldsymbol\mu)\mid \z\in\vertex(\prod_{j\in N}\Delta^{d_j-1}),\mu^c\in[\underline\mu^c,\overline \mu^c]\}$; Tolerance \(\epsilon>0\);
\Repeat
    \ForAll{\(c\in M\)}
        \State Separate \((\bar \z,\bar\mu^c)\) from
        \(\proj_{(\z,\mu^c)}(\mathcal P_c)\)
        \State Add a cut in the space of $(\z,\mu^c)$ to $\mathcal C$ if \((\bar \z,\bar\mu^c)\) is not feasible to $\proj_{(\z,\mu^c)}(\mathcal P_c)$
    \EndFor
    \State Solve the master problem $\max_{(\z,\boldsymbol\mu)\in\mathcal C}\sum_{c\in M}\omega_c\mu^c$ 
    \State $(\bar \z,\bar{\boldsymbol{\mu}})\gets$ the solution of the master problem
    \State \(UB\gets\sum_{c\in M}\omega_c\bar\mu^c\)
    \State \(LB\gets\) Evaluation of \(\bar \z\) in~\eqref{price optimization}
\Until{\(UB-LB\leq\epsilon\)}
\State \Return the best feasible price vector found
\end{algorithmic}
\end{algorithm}
\subsection{Separation Oracle}
To separate $\bar \z \in\vertex(\prod_{j\in N}\Delta^{d_j-1})$, we do not include all the inequalities in~\eqref{eq:coordinating}. Instead, define
$$
m_{jr,\ell t}^{\bar z}(\y^c):=
\begin{cases}
y_{jr}^c,&\bar z_{jr}\leq\bar z_{\ell t},\\
y_{\ell t}^c,&\bar z_{jr}>\bar z_{\ell t},
\end{cases}
\qquad
h_{jr,\ell t}^{\bar z}(\y^c,\rho^c):=
\begin{cases}
0,&\bar z_{jr}+\bar z_{\ell t}\leq1,\\
y_{jr}^c+y_{\ell t}^c-\rho^c,&\bar z_{jr}+\bar z_{\ell t}>1.
\end{cases}
$$
Consider a relaxation of $\mathcal P_c$, denoted as $\mathcal P_c(\bar \z)$, which restricts the constraints in~\eqref{eq:coordinating} to the following:
\begin{subequations}
\begin{align}
    &z_{jr}\geq
\left(1+\sum_{\ell\in N}\delta^c_{\ell1}\right)y_{jr}^c
+\sum_{t=1}^{d_j-1}\delta^c_{jt+1}y_{j,\max\{r,t\}}^c
+\sum_{\substack{\ell\neq j\\t\in[d_\ell-1]}}
\delta^c_{\ell t+1}m_{jr,\ell t}^{\bar z}(\y^c),\label{eq:zjrLowerDecomp}\\
&z_{jr}\leq
\left(1+\sum_{\ell\in N}\delta^c_{\ell1}\right)y_{jr}^c
+\sum_{t=1}^{d_j-1}\delta^c_{jt+1}y_{j,\max\{r,t\}}^c
+\sum_{\substack{\ell\neq j\\t\in[d_\ell-1]}}
\delta^c_{\ell t+1}h_{jr,\ell t}^{\bar z}(\y^c,\rho^c)\label{eq:zjrUpperDecomp}.
\end{align}
\end{subequations}
It is sufficient to separate $(\bar \z,\bar{\boldsymbol{\mu}})$ from $\mathcal P_c(\bar \z)$.
\begin{proposition}\label{necessarycut}
For every \(\bar \z\in\prod_{j\in N}\Delta^{d_j-1}\), \(\mathcal P_c\subseteq\mathcal P_c(\bar \z)\) and, for every \(\bar \mu^c\in\mathbb R\),
if $\bar \z\in \vertex(\prod_{j\in N}\Delta^{d_j-1})$ then \((\bar \z,\bar \mu^c)\in\proj_{(\z,\mu^c)}(\mathcal P_c)\) if and only if \((\bar \z,\bar \mu^c)\in\proj_{(\z,\mu^c)}(\mathcal P_c(\bar \z))\).
\end{proposition}

\begin{figure}[t]
\centering
\tikzsetnextfilename{priceopt-4}
\begin{tikzpicture}[scale=0.7,>=stealth,
   every node/.style={font=\small}]

\coordinate (xstar) at (0,0);
\coordinate (xbar) at (0,2.20);
\coordinate (leftfoot) at (-1.005,0.653);
\coordinate (rightfoot) at (1.005,0.653);

\begin{scope}
\clip (-4.1,-3.65) rectangle (4.1,2.75);
\fill[gray!20] (-4.2,-4.2) -- (4.2,-4.2) -- (xstar) -- cycle;
\end{scope}
\draw[black!75,very thick]
(-3.45,-3.45) -- (xstar) -- (3.45,-3.45);

\draw[blue!70!black,dashed,thick]
(-3.9,2.535) -- (3.9,-2.535);
\draw[gray!75!black,dashed,thick]
(-3.9,-2.535) -- (3.9,2.535);
\draw[orange!90!red,very thick]
(-4.0,0) -- (4.0,0);

\draw[blue!70!black,thick]
(xbar) -- (leftfoot);
\draw[gray!75!black,thick]
(xbar) -- (rightfoot);
\draw[orange!90!red,very thick,->]
(xbar) -- (xstar);

\filldraw[fill=orange!90!red,draw=black]
(xstar) circle (3.2pt);
\filldraw[fill=black,draw=black]
(xbar) circle (3.6pt);

\node[above=4pt] at (xbar)
{$(\bar z,\bar\mu^c)$};

\node[align=center] at (0,-2.45)
{$\proj_{(\z,\mu^c)}
\!\left(\mathcal P_c(\bar z)\right)$};

\node[orange!90!red,anchor=west]
at (2.70,-0.18) {deepest cut};

\node[blue!70!black,rotate=-33,anchor=west]
at (-3.05,2.18) {alternative cut};

\node[gray!75!black,rotate=33,anchor=west]
at (1.50,1.42) {alternative cut};

\end{tikzpicture}
\caption{Euclidean-deepest cut selection for a maximization problem.}
\label{fig:deepest-segment-cut}
\end{figure}
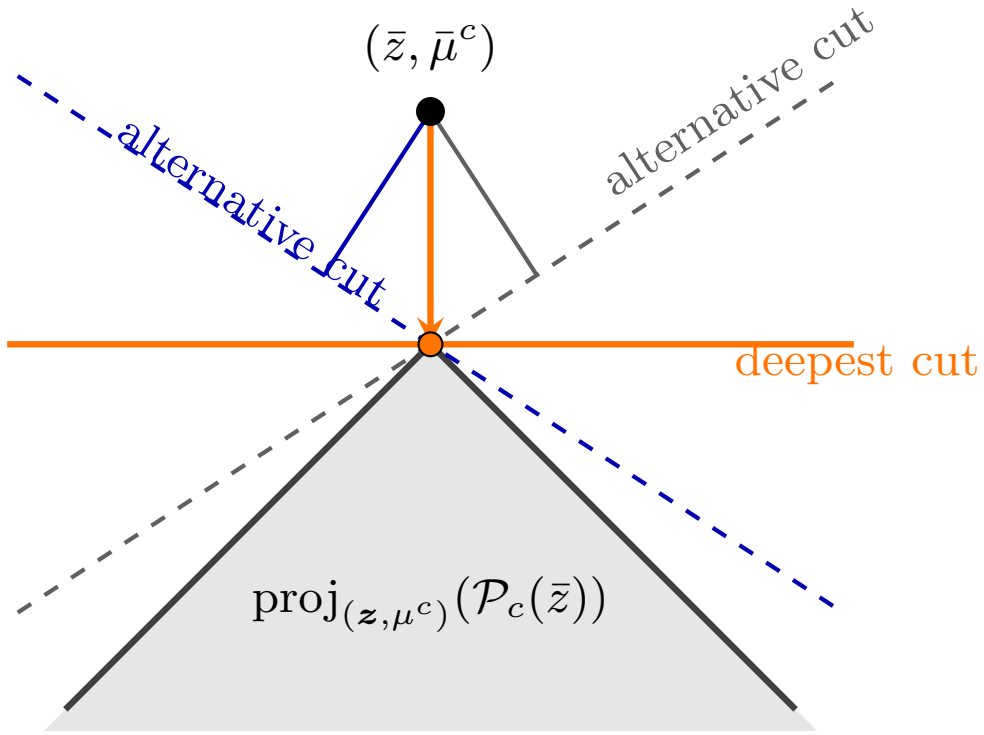

We separate $(\bar \z,\bar \mu^c)$ from $\mathcal P_c(\bar z)$ by solving a cut-generating program \eqref{eq:CGP}. Specifically, \eqref{eq:CGP} constructs the cut as a linear combination of the constraints defining \(\mathcal P_c(\bar{\z})\). Its multiplier ensures that the coefficients of the segment-specific variables \((\rho^c,\y^c)\) vanish, leaving an inequality involving only the master variables \((\z,\mu^c)\).
For notational simplicity, let \(\boldsymbol{t}_1:=(\z,\mu^c)\), \(\boldsymbol{t}_2:=(\y^c,\rho^c)\), and \(\bar{\boldsymbol{t}}_1:=(\bar \z,\bar\mu^c)\). We abbreviate \(\mathcal P_c(\bar \z)\) by
$
\mathcal P_c(\bar \z):=\left\{(\boldsymbol{t}_1,\boldsymbol{t}_2)\mid A\boldsymbol{t}_1+B\boldsymbol{t}_2\leq  d\right\}
$.
To generate a cut separating \(\bar{\boldsymbol{t}}_1\), we solve the following CGP:
\begin{equation}\label{eq:CGP}
\tag{\textsc{CGP}}
\begin{aligned}
\max_{\boldsymbol{u}}\quad &\boldsymbol{u}^\top(A\bar{\boldsymbol{t}}_1 -d)\\
\text{\rm s.t.}\quad &\boldsymbol{u}^\top B=0, \lVert \boldsymbol{u}^\top A\rVert_2\leq1, \boldsymbol{u}\geq0.
\end{aligned}
\end{equation}
If the optimal value of~\eqref{eq:CGP} is positive, 
let \(\boldsymbol{u}^*\) be an optimal solution and define
$
\boldsymbol\alpha^*:=A^\top \boldsymbol{u}^*, \boldsymbol\beta^*:=d^\top \boldsymbol{u}^*.
$
Then $\boldsymbol\alpha^{*\top}\boldsymbol{t}_1\leq \boldsymbol\beta^*$ separates $\bar{\boldsymbol{t}}_1$. Here, \eqref{eq:CGP} finds the deepest cut where the depth of a cut is measured by the Euclidean distance between the separated point and the cut hyperplane as shown in Figure~\ref{fig:deepest-segment-cut}. This is based on similar constructions in the literature~\citep{boyd1994fenchel,tawarmalani2002convexification,hosseini2025deepest}. Our algorithm includes several implementation details that we discuss in \ref{sec:implementation}.

\section{Numerical Results}
This section evaluates the computational performance and practical relevance of our method. We first show why price-sensitivity heterogeneity across customer segments and global, rather than local, optimality matter, then compare our method with \texttt{Gurobi} on synthetic and real-world instances. Throughout, we assume the cost of each product $\cost_i = 0$ for $i \in N$, and synthetic instances use the segment-specific utility function
\begin{equation}\label{eq:utility}
\mathcal U_{cj}=a_{cj}-b_{cj}x_j,
\end{equation}
where \(a_{cj}\in[-7,7]\) and \(b_{cj}\in[0.001,0.1]\). These ranges follow \cite{van2022price}, except that we expand \(b_{cj}\) from \([0.001,0.01]\) to \([0.001,0.1]\) to capture a broader range of price sensitivities. We derive bounds \(\underline x_j\) and \(\overline x_j\) such that \([\underline x_j,\overline x_j]\) contains an optimal price for product \(j\) (Section~\ref{bounds of prices}).

The code for optimization was written in \texttt{Julia} v1.11~(\cite{Julia-2017}), using \texttt{JuMP} v1.32.1~(\cite{lubin2023jump}), \texttt{Ipopt} v1.16.0~(\cite{wachter2006implementation}) (nonlinear solver), and \texttt{Gurobi} v13.0.3~(\cite{gurobi}) (MILP, subproblem and master problem solver).
The code for estimation was written in \texttt{Python} v3.12.
Estimation of FML model and optimization in Section~\ref{exp1} and~\ref{exp2} were conducted on a MacBook Pro equipped with an Apple M4 Pro processor (12-core CPU) and 24 GB of RAM.
The optimization problems in Sections~\ref{syn} and~\ref{exp4} were solved on a community cluster with Dell compute nodes, each equipped with two 96-core AMD EPYC Genoa processors and 384 GB of RAM. \texttt{Gurobi}’s branch-and-bound used 12 threads, while all other computations used a single thread.
\subsection{Impact of Ignoring Price Sensitivity Heterogeneity across Segments}\label{exp1}

The existing scalable algorithm \citep{van2022price} assumes that $b_{ci}$ is independent of $c$. This is practically restrictive, since consumer segments exhibit different price sensitivities, as supported by empirical studies in the literature. For example, FML models have been used to study the effect of price on the choice of consumer in the travel and tourism industry \citep{smith2007travelocity,cho2023estimating,namin2020improving}, where this sensitivity has been documented to vary across consumer segments. Specifically, \citet{ettl2020data} use airline transactions to find that the price elasticity of the first type of consumer is much lower than that of the second type. They suggest business-related products for the first type while discounts on group-related products for the latter. Similar evidence is also presented in \citet{meyer2024interpretable} and \citet{lu2020modelling}. We show in Section~\ref{sec:homogeneityRestrictive} that assuming that $b_{ci}$ is independent of $c$ also implies that arbitrary RUM models can no longer be uniformly approximated, unlike general FML~\citep{mcfadden2000mixed}.

We show that estimating models assuming this independence produces suboptimal revenue. For this, we generate synthetic choice data from a FML model with \(M=5\) customer segments and \(N=3\) products. The segment weights are \(\omega_c=1/M\). 
For estimation, each customer's segment is assumed to be known.

For each instance, $a_{cj}\sim U[-7,7]$ and $b_{cj}\sim U[0.05-\delta,0.05+\delta]$. We consider $\delta \in \{0.01,0.02,0.03,0.04\}$, where larger values admit more heterogeneity in price sensitivity. The price $x_j\sim U[\underline x_j,\overline x_j]$, and the  product is chosen according to the true MNL distribution. Each instance is trained using \(10^5\) observations. We compare two estimated models: the first allows segment-specific price sensitivity,
\[
    \text{Correctly specified model:}\qquad
    U_{cj}=\alpha_{cj}-(0.001+\softplus(\beta_{cj}))x_j,
\]
while the second restricts price sensitivity to be the same across customer segments for each product,
\[
    \text{Misspecified model:}\qquad
    U_{cj}=\alpha_{cj}-(0.001+\softplus(\beta_j))x_j .
\]
Both models are estimated by maximum likelihood. For each product $j$, we construct $10$ equally spaced admissible price levels in $[\underline x_j,\overline x_j]$ and solve the resulting discrete pricing problem globally.

\begin{figure}[htbp]
\centering
\includesvg[scale=0.3]{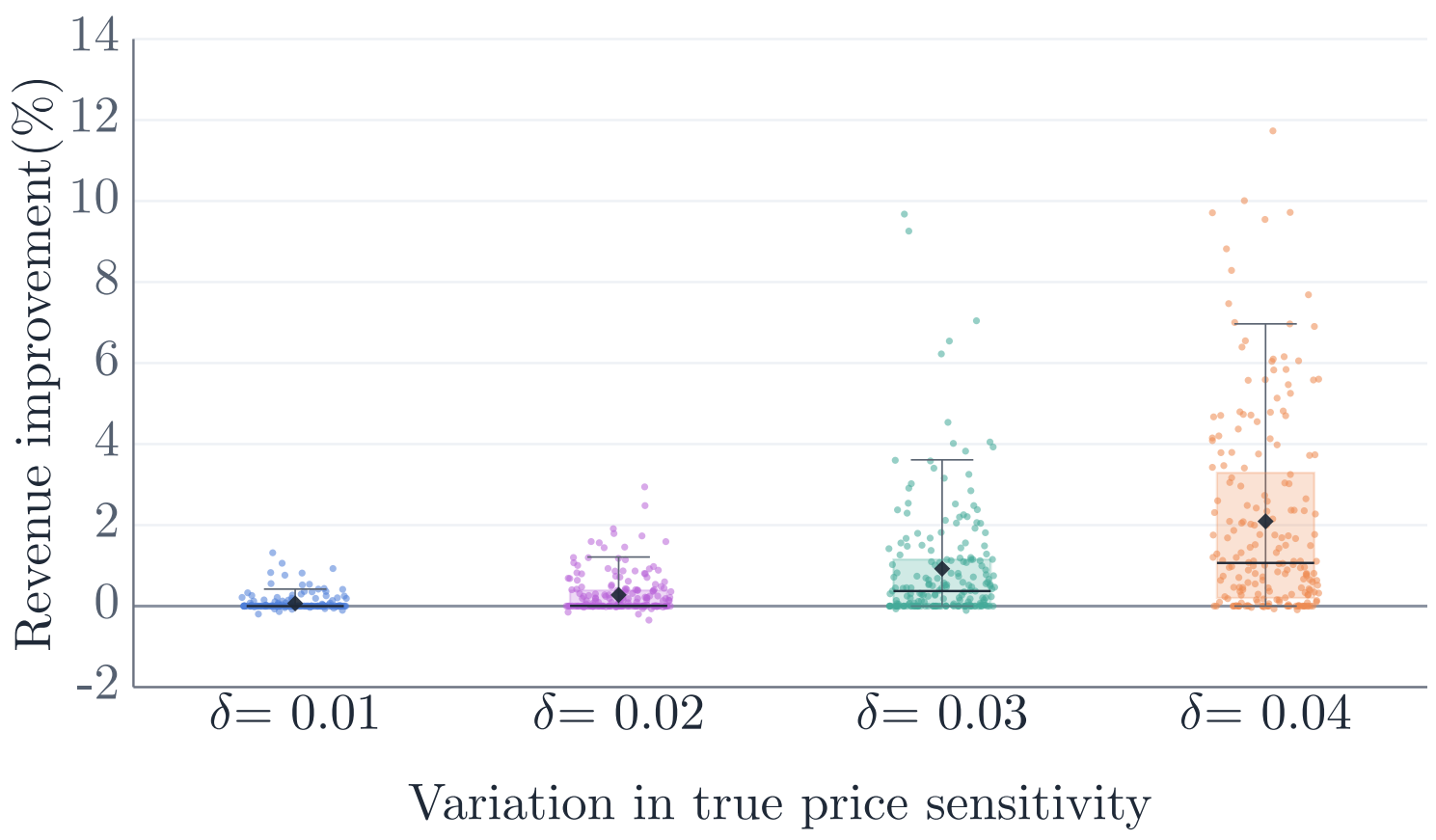}
\caption{Each point corresponds to an instance. The shaded box shows the interquartile range of revenue improvements, the short bars show the 5th and 95th percentiles, the long bar marks the median, and the black diamond marks the mean.}
\label{fig:price-sensitivity-box}
\end{figure}

We evaluate both pricing decisions under the true model and compute the relative improvement in revenue from the correctly specified model relative to that from the misspecified model. Figure~\ref{fig:price-sensitivity-box} shows that the revenue loss of the misspecified model increases with $\delta$ (price sensitivity varies more across segments).
 
\subsection{Revenue Impact of Relying on Local Optima}\label{exp2}
We next show that the revenue loss from local optimization can be substantial. For this experiment, we compare the revenue obtained from a local-search benchmark with the globally optimal revenue.
We generate $100$ synthetic instances with
$a_{cj} \sim \mathit{U}(-7,7)$, and $b_{cj} \sim \mathit{U}(0.001,0.1)$.
The segment weights \(\omega_c\) are sampled independently from \(U[0,1]\) and then normalized so that they sum to one. We show that even for small instances local optimality leads to revenue loss. Specifically, we set the number of segments to \(20\) and the number of products to \(20\). For each product, we uniformly select \(5\) breakpoints over the interval defined by the theoretical lower and upper bounds given by Proposition~\ref{prop:price-bounds}.

The local-search benchmark solves a continuous relaxation of the pricing problem over the theoretical price bounds using IPOPT. The resulting prices are then rounded to the nearest admissible grid points to obtain a discrete solution.
\begin{figure}[htbp]
\centering
\includesvg[scale=0.5]{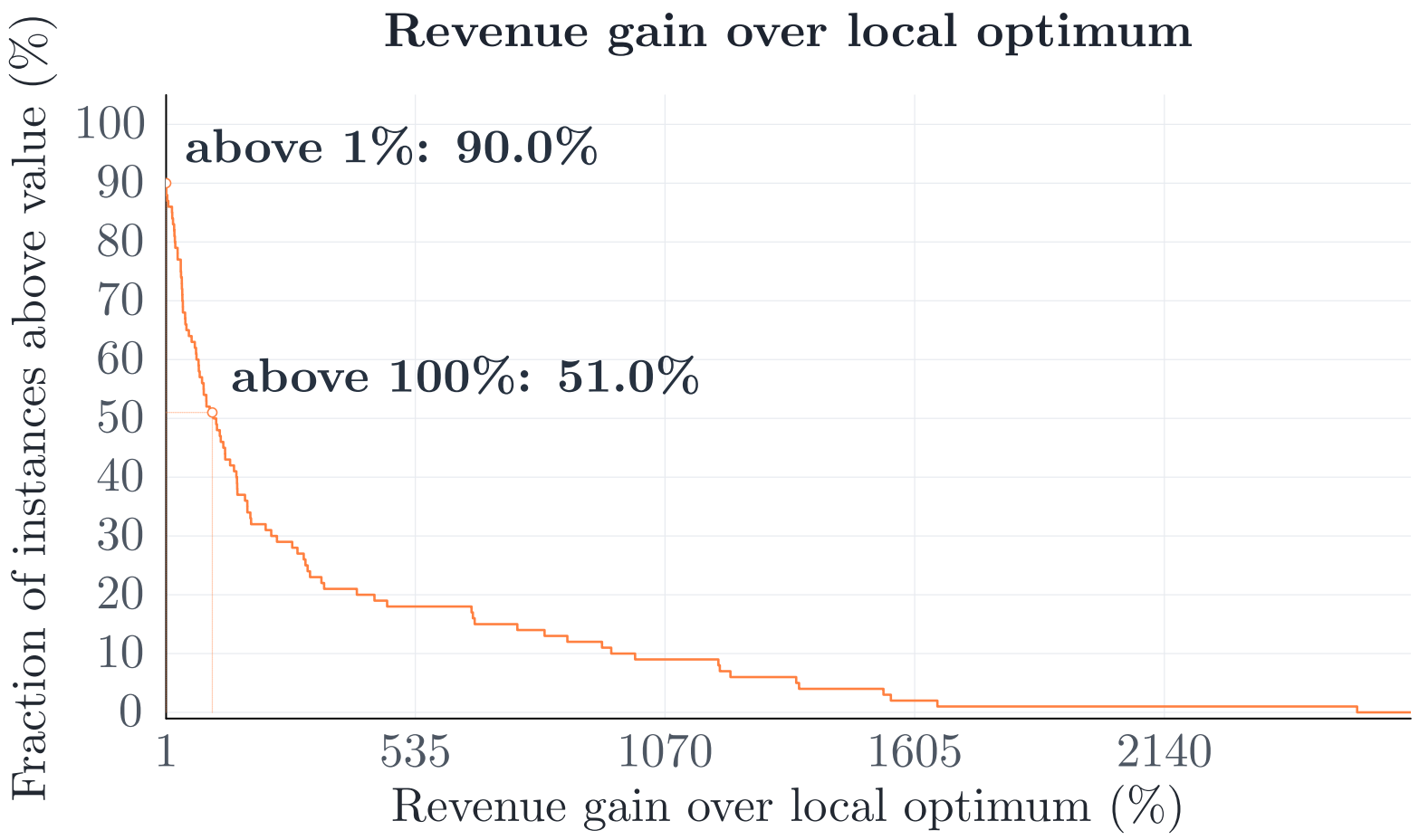}
\caption{Distribution of the revenue gain obtained by the global optimum over the local-search benchmark, defined as \(100(d-p)/p\), where \(d\) and \(p\) denote the globally optimal and local-search revenues, respectively.}
\label{fig:localvsglobal}
\end{figure}
Figure~\ref{fig:localvsglobal} presents the distribution of the revenue gain: for each value on the horizontal axis, the vertical axis reports the percentage of instances whose revenue gain exceeds that value. The distribution exhibits a pronounced right tail. Global optimization improves revenue by more than \(1\%\) in \(90.0\%\) of the instances and by more than \(100\%\) in \(51.0\%\) of the instances (doubles the revenue obtained by local search in over half of the instances). The largest improvement reaches \(2553\%\), demonstrating that local optimization frequently produces severely suboptimal pricing decisions, even for small instances.

\subsection{Synthetic Pricing Optimization Problem}
\label{syn}
We use synthetic instances for benchmarking. Our proposed approach explicitly accommodates heterogeneity in price sensitivity across customer segments while providing globally optimal solutions with certificates of optimality. Moreover, the results demonstrate that the method scales to large problem instances and substantially improves computational performance relative to the benchmark approaches.

For each segment $c$ and product $j$, we draw
$a_{cj} \sim \mathit{U}(-7,7)$, $b_{cj} \sim \mathit{U}(0.001,0.1)$, 
and normalize segment weights $\omega_c$ drawn from \(U[0,1]\). \(L\) breakpoints are selected within the interval \([\underline{x}_j,\overline{x}_j]\) given by Proposition~\ref{prop:price-bounds}. We choose $|M|\in\{15,20,25\}$, $L \in \{3,5,10\}$, and $|N|\in\{150,200,250\}$ and, for each combination, generate \(20\) instances. We compare the price decomposition (PD) method with solving~\eqref{price optimization} directly with \texttt{Gurobi} using unary binarization, itself a contribution of this paper, rather than the raw MINLP. The convergence criterion is a relative primal-dual gap below \(0.01\%\). The time limit is set to $3600$ seconds.

\begin{figure}
    \centering
    \includegraphics[width=\linewidth]{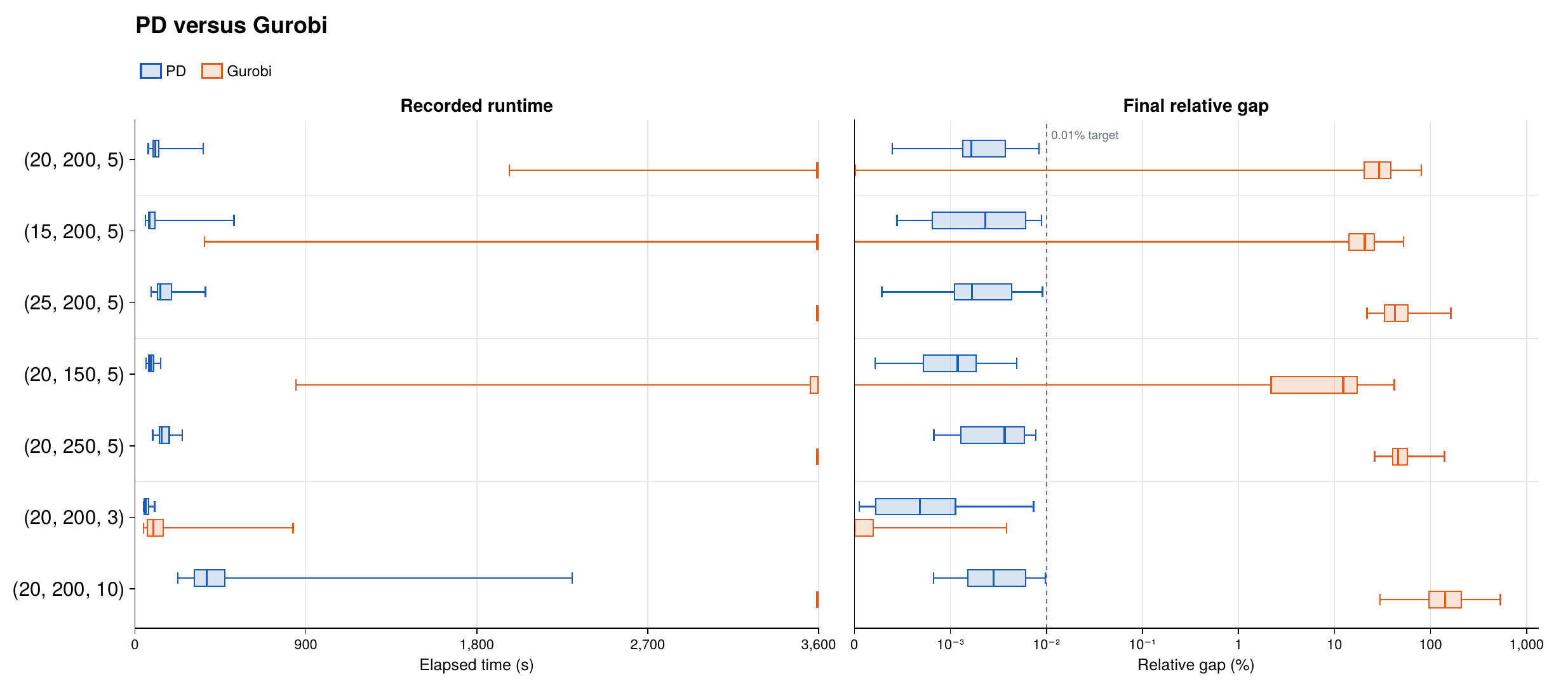}
    \caption{Runtime and relative gap of PD and Gurobi on synthetic instances, with 20 randomly generated instances per configuration. Starting from the baseline $(M,N,L)=(20,200,5)$, we vary one parameter at a time, fixing the other two.}
    \label{fig:PDGurobigap}
\end{figure}

\begin{figure}[htbp]
\centering
\includegraphics[width = \textwidth]{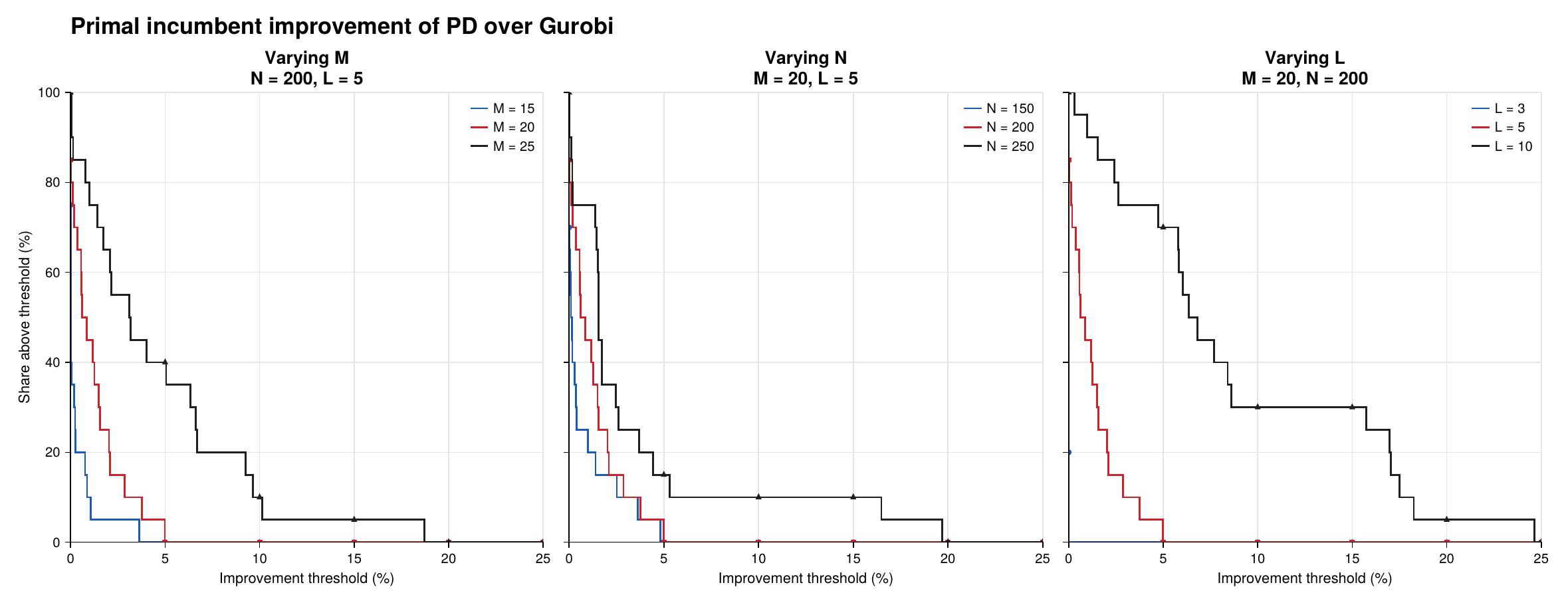}
\caption{Distribution of the primal incumbent improvement of PD over
\texttt{Gurobi} on synthetic instances. Each curve reports the percentage of instances whose incumbent improvement exceeds a given threshold. From left to right, the panels vary $M\in\{15,20,25\}$, $N\in\{150,200,250\}$, and $L\in\{3,5,10\}$, while holding the other two parameters at their baseline values $(M,N,L)=(20,200,5)$.}
\label{fig:synthetic-primal-improvement-tail}
\end{figure}

Figure~\ref{fig:PDGurobigap} shows that PD's advantage is substantial when direct optimization is difficult. At the baseline $(M,N,L)=(20,200,5)$, PD achieves a median final gap of $0.0016\%$ with a median runtime of $107.1$ seconds, whereas \texttt{Gurobi} has a median gap of $29.24\%$ after approximately $3600$ seconds. Across the five $L=5$ configurations, PD reaches the $0.01\%$ target on all $100$ instances, compared with only $9$ instances for \texttt{Gurobi}. At $L=10$, PD reaches the target on all $20$ instances, while \texttt{Gurobi} reaches it on none within the $3600$-second time limit. Both methods solve all $20$ of the easier $L=3$ instances. 

More importantly, Figure~\ref{fig:synthetic-primal-improvement-tail} shows that PD discovers better solutions. For half of the baseline instances, PD's solutions increase revenue relative to \texttt{Gurobi}'s best solution by more than $0.73\%$. As problem size increases, larger gains are more common: at $M=25$ (resp. $N=250$), $8$ (resp. $3$) of the $20$ instances gain more than $5\%$ in revenue. For $L=10$, the median gain is $6.58\%$, $14$ of $20$ instances gain more than $5\%$, and the largest gain is $24.63\%$. For the easy $L=3$ instances, the two methods yield identical revenues.

\subsection{Empirical Experiment over Real Data from Expedia}\label{exp4}
Now, we consider a real dataset.
We estimate a segmented MNL model using the Expedia Personalized Hotel Search data released for the ICDM 2013 competition~(\cite{expedia-personalized-sort}). 
This is a common approach for the estimation of FML in the literature \citep{ettl2020data}, where the authors cluster airline data into $7$ customer types before fitting MNL models for each segment. Given these segment assignments, we estimate each product's intrinsic preference and price sensitivity for each segment by maximum likelihood regularized with L2 penalty~(\cite{hoerl1970ridge}) based on the choice probabilities in~\ref{price optimization}.

We perform several steps to transform this data into a pricing problem. Each search session is treated as a choice observation over a set of available hotels and a no-purchase option. 
Since the data contains multiple destinations, we group geographically proximate destinations when they belong to the same local market. It is reasonable for a customer’s choice set to include hotels in both Minneapolis and St. Paul, the well-known Twin Cities located across the Mississippi River from each other. 
We construct a destination graph in which each node is a
destination and connect two destinations if they share at least one
hotel. For each selected center destination, we define its local market as the union of its connected component. This construction is meant to capture the local hotel market with sufficient observations in each market for estimation.

We cluster hotels using observable hotel characteristics, including star rating, review score, brand indicator, location score, and historical log price. These features are normalized prior to clustering, and each resulting cluster is treated as one product in the pricing model. Within each search session, a product type may contain multiple available
hotels. In that case, we represent the price of this product type as the minimum observed price among hotels of that type in that search session. If the booked hotel belongs to a retained product type, that product type is recorded as the chosen alternative. Sessions where no booking is recorded are assumed to have selected the no-purchase option. We remove extreme price observations by retaining only rows satisfying \(10 \leq \text{price} \leq 1000\). Product types that appear in fewer than ten search sessions and sessions where these products are selected are discarded.

Customers are classified into segments based on five characteristics: whether the visitor location country is $219$ (approximately 57\% of the customers), whether the length of their stay is one night, whether they search for one adult, whether children are present, and whether the stay includes a Saturday night. This gives at most $2^5=32$ customer segments. The segment weight $\omega_c$ is the empirical share of segment $c$. To obtain economically meaningful parameters, we estimate the model using the smooth reparametrizations
$a_{cj}=7\tanh(\theta_{cj})$,
and $b_{cj}=\operatorname{softplus}(\beta_{cj})+0.001$.
This ensures that intrinsic preferences are in $[-7,7]$ and price sensitivities are strictly positive and greater than $0.001$. The estimated parameters \((a_{cj},b_{cj},\omega_c)\) are then used directly as inputs to the pricing optimization problem.
We randomly divide the retained search sessions into training and validation sets containing $80\%$ and $20\%$ of the observations, respectively.
\begin{table}[htbp]
\centering
\begin{tabular}{ccccccc}
\toprule
Center
& Products
& Segments
& Choices
& Val. NLL
& Max weight\\
\midrule
8467  & 96 & 32 & 24686 & 2.2058  &  13.29\% \\
4562  & 99 & 32 & 24591 & 2.0199  &  10.24\% \\
23246 & 97 & 32 & 23660 & 2.2174 &  13.69\% \\
\bottomrule
\end{tabular}
\caption{Summary of empirical pricing instances calibrated from Expedia hotel-search data. Val. NLL is the average validation negative log-likelihood per choice observation. Max weight is $\max_c \omega_c$.}
\label{tab:expedia-estimation-summary}
\end{table}

Table~\ref{tab:expedia-estimation-summary} reports the three largest center-destination instances, corresponding to centers \(4562\), \(8467\), and \(23246\).
An overview of the estimated parameters is reported in Section~\ref{overviewofpara}.
The estimated models exhibit consumer heterogeneity, and the largest segment weight is below \(14\%\). For the pricing optimization problem, prices are restricted to a discrete set. For each product type \(j\), we first collect all distinct prices observed for that product type in the search sessions. We sort the distinct values in increasing order. If the number of distinct observed prices is at most \(L=5\), all of them are retained as feasible prices. Otherwise, we select \(L=5\) evenly spaced order statistics from the sorted list.

\begin{figure}[htbp]
\centering
\includegraphics[width = 1\textwidth]{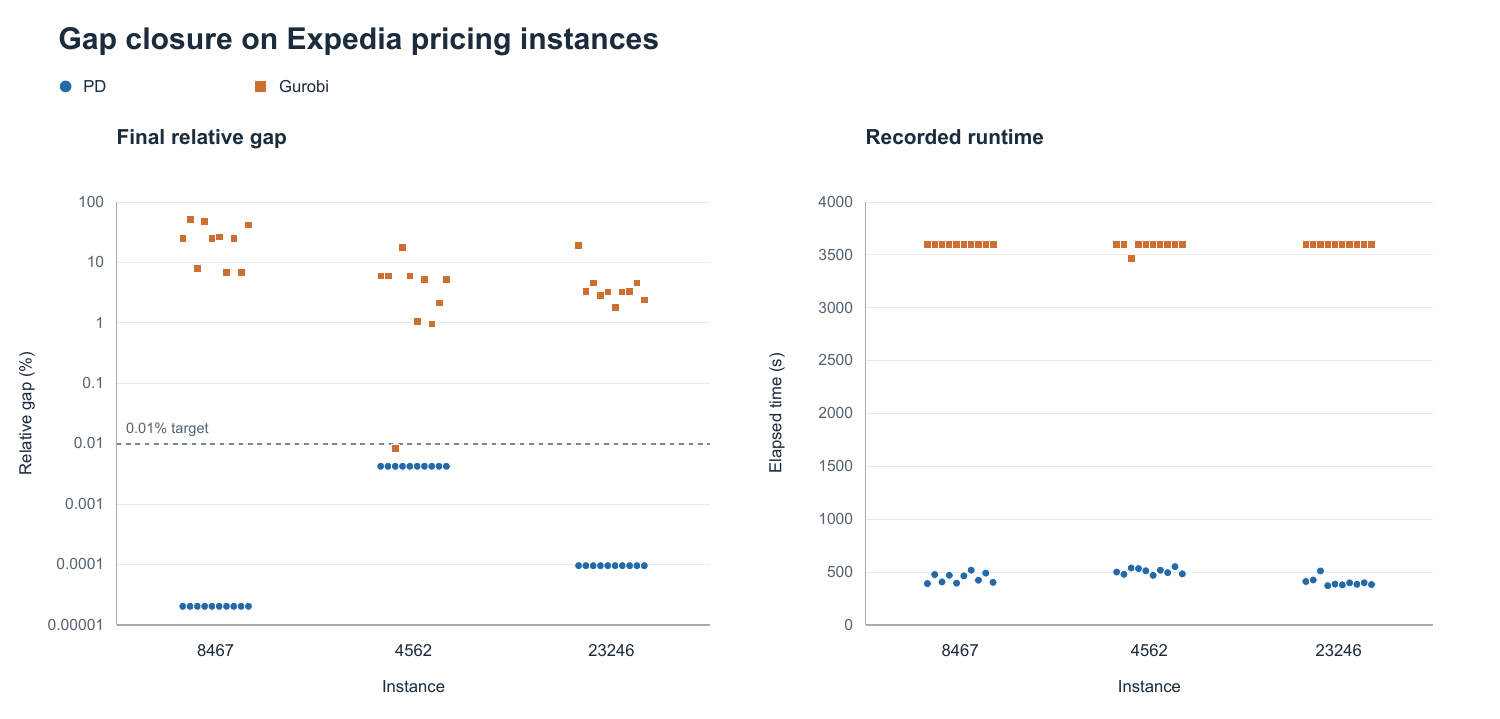}
\caption{Comparison of PD and Gurobi for ten repeated runs for each of three Expedia destination centers: 8467, 4562, and 23246. Both methods use identical model parameters and price grids. Each marker represents one run; horizontal offsets separate the repetitions for visibility. Left: relative gap. The dashed line marks the $0.01\%$ gap tolerance. Right: recorded elapsed time in seconds.}
\label{fig:expedia-gap-runtime}
\end{figure}

For each instance, we run PD and the Gurobi solver $10$ times. The time limit for each is set to \(3600\) seconds. We report the running time and the final relative optimality gap in Figure~\ref{fig:expedia-gap-runtime}. A problem is considered solved if the final relative gap is below \(0.01\%\). PD achieved this in all ten runs for each center, while Gurobi did not achieve this gap for centers 8467 and 23246, and achieved it once for center 4562.

\section{Conclusions}\label{sec:conclusions}

In this paper, we develop a new MIP formulation for pricing under the FML consumer choice model and a decomposition algorithm to solve this problem. This formulation uncovers a surprising convexification effect of using the unary binarization scheme along with homogenization of McCormick inequalities. To address continuous pricing, our methods can be naturally extended, by refining the discretization domain and ensuring that the upper bound is valid between the discrete points. We demonstrate that our techniques can solve practical sized instances where state-of-the-art approaches either ignore consumer price sensitivity or struggle to find optimal solutions. We also show that global optimization of this pricing problem leads to much better solutions than local optimization and rounding.

\bibliographystyle{plainnat} 

\begin{thebibliography}{61}
\providecommand{\natexlab}[1]{#1}
\providecommand{\url}[1]{\texttt{#1}}
\expandafter\ifx\csname urlstyle\endcsname\relax
  \providecommand{\doi}[1]{doi: #1}\else
  \providecommand{\doi}{doi: \begingroup \urlstyle{rm}\Url}\fi

\bibitem[Abdolhamidi and Lurkin(2025)]{abdolhamidi2025tactical}
Dorsa Abdolhamidi and Virginie Lurkin.
\newblock A tactical time slot management problem under mixed logit demand.
\newblock \emph{OR Spectrum}, pages 1--45, 2025.

\bibitem[Adam et~al.(2013)Adam, Hamner, Friedman, and
  SSA\_Expedia]{expedia-personalized-sort}
Adam, Ben Hamner, Dan Friedman, and SSA\_Expedia.
\newblock Personalize expedia hotel searches - icdm 2013.
\newblock \url{https://kaggle.com/competitions/expedia-personalized-sort},
  2013.
\newblock Kaggle.

\bibitem[Arslan et~al.(2022)Arslan, Easley, Wang, and
  Y{\i}lmaz]{arslan2022data}
Hayri~A Arslan, Robert~F Easley, Ruxian Wang, and {\"O}v{\"u}n{\c{c}}
  Y{\i}lmaz.
\newblock Data-driven sports ticket pricing for multiple sales channels with
  heterogeneous customers.
\newblock \emph{Manufacturing \& Service Operations Management}, 24\penalty0
  (2):\penalty0 1241--1260, 2022.

\bibitem[Benders(1962)]{benders1962partitioning}
Jacques~F Benders.
\newblock Partitioning procedures for solving mixed-variables programming
  problems.
\newblock \emph{Numerische mathematik}, 4\penalty0 (1):\penalty0 238--252,
  1962.

\bibitem[Bezanson et~al.(2017)Bezanson, Edelman, Karpinski, and
  Shah]{Julia-2017}
Jeff Bezanson, Alan Edelman, Stefan Karpinski, and Viral~B Shah.
\newblock \text{Julia}: A fresh approach to numerical computing.
\newblock \emph{SIAM {R}eview}, 59\penalty0 (1):\penalty0 65--98, 2017.
\newblock \doi{10.1137/141000671}.
\newblock URL \url{https://epubs.siam.org/doi/10.1137/141000671}.

\bibitem[Bitran and Caldentey(2003)]{bitran2003overview}
Gabriel Bitran and Ren{\'e} Caldentey.
\newblock An overview of pricing models for revenue management.
\newblock \emph{Manufacturing \& Service Operations Management}, 5\penalty0
  (3):\penalty0 203--229, 2003.

\bibitem[Block and Marschak(1959)]{block1959random}
Henry~David Block and Jacob Marschak.
\newblock Random orderings and stochastic theories of response.
\newblock 1959.

\bibitem[Bolusani et~al.(2024)Bolusani, Besan{\c{c}}on, Bestuzheva, Chmiela,
  Dion{\'\i}sio, Donkiewicz, van Doornmalen, Eifler, Ghannam, Gleixner,
  et~al.]{bolusani2024scip}
Suresh Bolusani, Mathieu Besan{\c{c}}on, Ksenia Bestuzheva, Antonia Chmiela,
  Jo{\~a}o Dion{\'\i}sio, Tim Donkiewicz, Jasper van Doornmalen, Leon Eifler,
  Mohammed Ghannam, Ambros Gleixner, et~al.
\newblock The {SCIP} optimization suite 9.0.
\newblock \emph{arXiv preprint arXiv:2402.17702}, 2024.

\bibitem[Boyd(1994)]{boyd1994fenchel}
E~Andrew Boyd.
\newblock Fenchel cutting planes for integer programs.
\newblock \emph{Operations Research}, 42\penalty0 (1):\penalty0 53--64, 1994.

\bibitem[Caro and Gallien(2012)]{caro2012clearance}
Felipe Caro and J{\'e}r{\'e}mie Gallien.
\newblock Clearance pricing optimization for a fast-fashion retailer.
\newblock \emph{Operations research}, 60\penalty0 (6):\penalty0 1404--1422,
  2012.

\bibitem[Charnes and Cooper(1962)]{charnes1962programming}
Abraham Charnes and William~W Cooper.
\newblock Programming with linear fractional functionals.
\newblock \emph{Naval Research logistics quarterly}, 9\penalty0 (3-4):\penalty0
  181--186, 1962.

\bibitem[Cho et~al.(2023)Cho, Ferguson, Pekg{\"u}n, and
  Vakhutinsky]{cho2023estimating}
Sanghoon Cho, Mark Ferguson, Pelin Pekg{\"u}n, and Andrew Vakhutinsky.
\newblock Estimating personalized demand with unobserved no-purchases using a
  mixture model: An application in the hotel industry.
\newblock \emph{Manufacturing \& Service Operations Management}, 25\penalty0
  (4):\penalty0 1245--1262, 2023.

\bibitem[Cohen et~al.(2017)Cohen, Leung, Panchamgam, Perakis, and
  Smith]{cohen2017impact}
Maxime~C Cohen, Ngai-Hang~Zachary Leung, Kiran Panchamgam, Georgia Perakis, and
  Anthony Smith.
\newblock The impact of linear optimization on promotion planning.
\newblock \emph{Operations Research}, 65\penalty0 (2):\penalty0 446--468, 2017.

\bibitem[Cohen et~al.(2021)Cohen, Kalas, and Perakis]{cohen2021promotion}
Maxime~C Cohen, Jeremy~J Kalas, and Georgia Perakis.
\newblock Promotion optimization for multiple items in supermarkets.
\newblock \emph{Management Science}, 67\penalty0 (4):\penalty0 2340--2364,
  2021.

\bibitem[Davis et~al.(2017)Davis, Topaloglu, and Williamson]{davis2017pricing}
James~M Davis, Huseyin Topaloglu, and David~P Williamson.
\newblock Pricing problems under the nested logit model with a quality
  consistency constraint.
\newblock \emph{INFORMS Journal on Computing}, 29\penalty0 (1):\penalty0
  54--76, 2017.

\bibitem[De~Loera et~al.(2010)De~Loera, Rambau, and
  Santos]{de2010triangulations}
J.~De~Loera, J.~Rambau, and F.~Santos.
\newblock \emph{Triangulations: Structures for Algorithms and Applications}.
\newblock Algorithms and Computation in Mathematics. Springer Berlin
  Heidelberg, 2010.

\bibitem[Ettl et~al.(2020)Ettl, Harsha, Papush, and Perakis]{ettl2020data}
Markus Ettl, Pavithra Harsha, Anna Papush, and Georgia Perakis.
\newblock A data-driven approach to personalized bundle pricing and
  recommendation.
\newblock \emph{Manufacturing \& Service Operations Management}, 22\penalty0
  (3):\penalty0 461--480, 2020.

\bibitem[Falmagne(1978)]{falmagne1978representation}
Jean-Claude Falmagne.
\newblock A representation theorem for finite random scale systems.
\newblock \emph{Journal of Mathematical Psychology}, 18\penalty0 (1):\penalty0
  52--72, 1978.

\bibitem[Ferreira et~al.(2016)Ferreira, Lee, and
  Simchi-Levi]{ferreira2016analytics}
Kris~Johnson Ferreira, Bin Hong~Alex Lee, and David Simchi-Levi.
\newblock Analytics for an online retailer: Demand forecasting and price
  optimization.
\newblock \emph{Manufacturing \& service operations management}, 18\penalty0
  (1):\penalty0 69--88, 2016.

\bibitem[Fischetti et~al.(2010)Fischetti, Salvagnin, and
  Zanette]{fischetti2010note}
Matteo Fischetti, Domenico Salvagnin, and Arrigo Zanette.
\newblock A note on the selection of benders’ cuts.
\newblock \emph{Mathematical Programming}, 124\penalty0 (1):\penalty0 175--182,
  2010.

\bibitem[Gallego and Wang(2014)]{gallego2014multiproduct}
Guillermo Gallego and Ruxian Wang.
\newblock Multiproduct price optimization and competition under the nested
  logit model with product-differentiated price sensitivities.
\newblock \emph{Operations Research}, 62\penalty0 (2):\penalty0 450--461, 2014.

\bibitem[Guan and Mi{\v{s}}i{\'c}(2026)]{guan2026randomized}
Xinyi Guan and Velibor~V Mi{\v{s}}i{\'c}.
\newblock Randomized robust price optimization.
\newblock \emph{Management Science}, 2026.
\newblock \doi{10.1287/mnsc.2023.01844}.

\bibitem[{Gurobi Optimization, LLC}(2022)]{gurobi}
{Gurobi Optimization, LLC}.
\newblock {Gurobi Optimizer Reference Manual}, 2022.
\newblock URL \url{https://www.gurobi.com}.

\bibitem[Hansen et~al.(1990)Hansen, Poggi~de Arag{\~a}o, and
  Ribeiro]{hansen1990boolean}
Pierre Hansen, Marcus~V Poggi~de Arag{\~a}o, and Celso~C Ribeiro.
\newblock Boolean query optimization and the 0-1 hyperbolic sum problem.
\newblock \emph{Annals of Mathematics and Artificial Intelligence}, 1:\penalty0
  97--109, 1990.

\bibitem[Hanson and Martin(1996)]{hanson1996optimizing}
Ward Hanson and Kipp Martin.
\newblock Optimizing multinomial logit profit functions.
\newblock \emph{Management Science}, 42\penalty0 (7):\penalty0 992--1003, 1996.

\bibitem[Harsha et~al.(2019)Harsha, Subramanian, and Ettl]{harsha2019practical}
Pavithra Harsha, Shivaram Subramanian, and Markus Ettl.
\newblock A practical price optimization approach for omnichannel retailing.
\newblock \emph{INFORMS Journal on Optimization}, 1\penalty0 (3):\penalty0
  241--264, 2019.

\bibitem[He and Tawarmalani(2022)]{he2022tractable}
Taotao He and Mohit Tawarmalani.
\newblock Tractable relaxations of composite functions.
\newblock \emph{Mathematics of Operations Research}, 47\penalty0 (2):\penalty0
  1110--1140, 2022.

\bibitem[He and
  Tawarmalani(2026)]{he2026discretenessconvexitypromotionplanning}
Taotao He and Mohit Tawarmalani.
\newblock Discreteness to convexity: Promotion planning via simplotope
  triangulation, 2026.
\newblock URL \url{https://arxiv.org/abs/2408.04562}.

\bibitem[He et~al.(2025)He, Liu, and Tawarmalani]{he2025convexification}
Taotao He, Siyue Liu, and Mohit Tawarmalani.
\newblock Convexification techniques for fractional programs.
\newblock \emph{Mathematical Programming}, 213\penalty0 (1):\penalty0 107--149,
  2025.

\bibitem[Hoerl and Kennard(1970)]{hoerl1970ridge}
Arthur~E Hoerl and Robert~W Kennard.
\newblock Ridge regression: Biased estimation for nonorthogonal problems.
\newblock \emph{Technometrics}, 12\penalty0 (1):\penalty0 55--67, 1970.

\bibitem[Hosseini and Turner(2025)]{hosseini2025deepest}
Mojtaba Hosseini and John Turner.
\newblock Deepest cuts for benders decomposition.
\newblock \emph{Operations Research}, 73\penalty0 (5):\penalty0 2591--2609,
  2025.

\bibitem[Hu et~al.(2022)Hu, Simchi-Levi, and Yan]{hu2022learning}
Yiqun Hu, David Simchi-Levi, and Zhenzhen Yan.
\newblock Learning mixed multinomial logits with provable guarantees.
\newblock \emph{Advances in Neural Information Processing Systems},
  35:\penalty0 9447--9459, 2022.

\bibitem[Kannan et~al.(2009)Kannan, Pope, and Jain]{kannan2009practice}
PK~Kannan, Barbara~Kline Pope, and Sanjay Jain.
\newblock Practice prize winner—pricing digital content product lines: A
  model and application for the national academies press.
\newblock \emph{Marketing Science}, 28\penalty0 (4):\penalty0 620--636, 2009.

\bibitem[Ke et~al.(2019)Ke, Zhang, and Zheng]{ke2019approximate}
Jiannan Ke, Dan Zhang, and Huan Zheng.
\newblock An approximate dynamic programming approach to dynamic pricing for
  network revenue management.
\newblock \emph{Production and Operations Management}, 28\penalty0
  (11):\penalty0 2719--2737, 2019.

\bibitem[Kim et~al.(2022)Kim, Richard, and Tawarmalani]{kim2022reciprocity}
Jongeun Kim, JP~Richard, and Mohit Tawarmalani.
\newblock A reciprocity between tree ensemble optimization and multilinear
  optimization.
\newblock \emph{Optimization Online, https://optimization-online.
  org/2022/03/8828}, 2022.

\bibitem[Lewis(2005)]{lewis2005research}
Michael Lewis.
\newblock Research note: A dynamic programming approach to customer
  relationship pricing.
\newblock \emph{Management science}, 51\penalty0 (6):\penalty0 986--994, 2005.

\bibitem[Li and Huh(2011)]{li2011pricing}
Hongmin Li and Woonghee~Tim Huh.
\newblock Pricing multiple products with the multinomial logit and nested logit
  models: Concavity and implications.
\newblock \emph{Manufacturing \& Service Operations Management}, 13\penalty0
  (4):\penalty0 549--563, 2011.

\bibitem[Li et~al.(2019)Li, Webster, Mason, and Kempf]{li2019product}
Hongmin Li, Scott Webster, Nicholas Mason, and Karl Kempf.
\newblock Product-line pricing under discrete mixed multinomial logit demand:
  winner—2017 m\&som practice-based research competition.
\newblock \emph{Manufacturing \& Service Operations Management}, 21\penalty0
  (1):\penalty0 14--28, 2019.

\bibitem[Linderoth and Wright(2003)]{linderoth2003decomposition}
Jeff Linderoth and Stephen Wright.
\newblock Decomposition algorithms for stochastic programming on a
  computational grid.
\newblock \emph{Computational Optimization and Applications}, 24\penalty0
  (2):\penalty0 207--250, 2003.

\bibitem[Lu et~al.(2020)Lu, Wang, Gu, and Yang]{lu2020modelling}
Jing Lu, Zheqing Wang, Yingjing Gu, and Wenhui Yang.
\newblock Modelling the air ticket purchase behavior incorporating latent class
  model.
\newblock \emph{Mathematical Problems in Engineering}, 2020\penalty0
  (1):\penalty0 2046106, 2020.

\bibitem[Lubin et~al.(2023)Lubin, Dowson, Garcia, Huchette, Legat, and
  Vielma]{lubin2023jump}
Miles Lubin, Oscar Dowson, Joaquim~Dias Garcia, Joey Huchette, Beno{\^\i}t
  Legat, and Juan~Pablo Vielma.
\newblock Jump 1.0: recent improvements to a modeling language for mathematical
  optimization.
\newblock \emph{Mathematical Programming Computation}, pages 1--9, 2023.

\bibitem[Marandi and Lurkin(2023)]{marandi2023exact}
Ahmadreza Marandi and Virginie Lurkin.
\newblock An exact algorithm for the static pricing problem under discrete
  mixed logit demand.
\newblock \emph{EURO Journal on Computational Optimization}, 11:\penalty0
  100073, 2023.

\bibitem[Markowitz and Manne(1957)]{markowitz1957solution}
Harry~M Markowitz and Alan~S Manne.
\newblock On the solution of discrete programming problems.
\newblock \emph{Econometrica: journal of the Econometric Society}, pages
  84--110, 1957.

\bibitem[McCormick(1976)]{mccormick1976computability}
Garth~P McCormick.
\newblock Computability of global solutions to factorable nonconvex programs:
  Part i — \text{Convex} underestimating problems.
\newblock \emph{Mathematical Programming}, 10\penalty0 (1):\penalty0 147--175,
  1976.

\bibitem[McFadden(1974)]{mcfadden1972conditional}
Daniel McFadden.
\newblock Conditional logit analysis of qualitative choice behavior.
\newblock In Paul Zarembka, editor, \emph{Fontiers in {Econometrics}}, pages
  105--142. Academic press, New York, 1974.

\bibitem[McFadden and Train(2000)]{mcfadden2000mixed}
Daniel McFadden and Kenneth Train.
\newblock Mixed {MNL} models for discrete response.
\newblock \emph{Journal of applied Econometrics}, 15\penalty0 (5):\penalty0
  447--470, 2000.

\bibitem[Mehmanchi et~al.(2019)Mehmanchi, G{\'o}mez, and
  Prokopyev]{mehmanchi2019fractional}
Erfan Mehmanchi, Andr{\'e}s G{\'o}mez, and Oleg~A Prokopyev.
\newblock Fractional 0--1 programs: links between mixed-integer linear and
  conic quadratic formulations: E. mehmanchi et al.
\newblock \emph{Journal of Global Optimization}, 75\penalty0 (2):\penalty0
  273--339, 2019.

\bibitem[Meyer et~al.(2024)Meyer, Kauermann, and Smith]{meyer2024interpretable}
Jan~Felix Meyer, Go{\"{}}~ran Kauermann, and Michael~Stanley Smith.
\newblock Interpretable modelling of retail demand and price elasticity for
  passenger flights using booking data.
\newblock \emph{Statistical Modelling}, 24\penalty0 (1):\penalty0 82--106,
  2024.

\bibitem[Namin et~al.(2020)Namin, Gauri, and Kwortnik]{namin2020improving}
Aidin Namin, Dinesh~K Gauri, and Robert~J Kwortnik.
\newblock Improving revenue performance with third-degree price discrimination
  in the cruise industry.
\newblock \emph{International Journal of Hospitality Management}, 89:\penalty0
  102597, 2020.

\bibitem[Schindler and Kibarian(1996)]{schindler1996increased}
Robert~M Schindler and Thomas~M Kibarian.
\newblock Increased consumer sales response though use of 99-ending prices.
\newblock \emph{Journal of Retailing}, 72\penalty0 (2):\penalty0 187--199,
  1996.

\bibitem[Schrijver(1983)]{schrijver1983short}
Alexander Schrijver.
\newblock Short proofs on the matching polyhedron.
\newblock \emph{Journal of Combinatorial Theory, Series B}, 34\penalty0
  (1):\penalty0 104--108, 1983.

\bibitem[{\c{S}}en et~al.(2018){\c{S}}en, Atamt{\"u}rk, and
  Kaminsky]{csen2018conic}
Alper {\c{S}}en, Alper Atamt{\"u}rk, and Philip Kaminsky.
\newblock A conic integer optimization approach to the constrained assortment
  problem under the mixed multinomial logit model.
\newblock \emph{Operations Research}, 66\penalty0 (4):\penalty0 994--1003,
  2018.

\bibitem[Smith et~al.(2007)Smith, Darrow, Elieson, Guenther, Rao, and
  Zouaoui]{smith2007travelocity}
Barry~C Smith, Ross Darrow, John Elieson, Dirk Guenther, B~Venkateshwara Rao,
  and Faker Zouaoui.
\newblock Travelocity becomes a travel retailer.
\newblock \emph{Interfaces}, 37\penalty0 (1):\penalty0 68--81, 2007.

\bibitem[Tawarmalani and Sahinidis(2002)]{tawarmalani2002convexification}
Mohit Tawarmalani and Nikolaos~V Sahinidis.
\newblock \emph{Convexification and global optimization in continuous and
  mixed-integer nonlinear programming: theory, algorithms, software, and
  applications}, volume~65.
\newblock Springer Science \& Business Media, 2002.

\bibitem[Tawarmalani and Sahinidis(2004)]{tawarmalani2004global}
Mohit Tawarmalani and Nikolaos~V Sahinidis.
\newblock Global optimization of mixed-integer nonlinear programs: A
  theoretical and computational study.
\newblock \emph{Mathematical programming}, 99\penalty0 (3):\penalty0 563--591,
  2004.

\bibitem[Tawarmalani et~al.(2002)Tawarmalani, Ahmed, and
  Sahinidis]{tawarmalani2002global}
Mohit Tawarmalani, Shabbir Ahmed, and Nikolaos~V Sahinidis.
\newblock Global optimization of 0-1 hyperbolic programs.
\newblock \emph{Journal of Global Optimization}, 24\penalty0 (4):\penalty0
  385--416, 2002.

\bibitem[Tawarmalani et~al.(2013)Tawarmalani, Richard, and
  Xiong]{tawarmalani2013explicit}
Mohit Tawarmalani, Jean-Philippe~P Richard, and Chuanhui Xiong.
\newblock Explicit convex and concave envelopes through polyhedral
  subdivisions.
\newblock \emph{Mathematical Programming}, 138\penalty0 (1-2):\penalty0
  531--577, 2013.

\bibitem[Thomas and Morwitz(2005)]{thomas2005penny}
Manoj Thomas and Vicki Morwitz.
\newblock Penny wise and pound foolish: the left-digit effect in price
  cognition.
\newblock \emph{Journal of consumer Research}, 32\penalty0 (1):\penalty0
  54--64, 2005.

\bibitem[Train(2009)]{train2009discrete}
Kenneth~E Train.
\newblock \emph{Discrete choice methods with simulation}.
\newblock Cambridge university press, 2009.

\bibitem[van~de Geer and den Boer(2022)]{van2022price}
Ruben van~de Geer and Arnoud~V den Boer.
\newblock Price optimization under the finite-mixture logit model.
\newblock \emph{Management Science}, 68\penalty0 (10):\penalty0 7480--7496,
  2022.

\bibitem[W{\"a}chter and Biegler(2006)]{wachter2006implementation}
Andreas W{\"a}chter and Lorenz~T Biegler.
\newblock On the implementation of an interior-point filter line-search
  algorithm for large-scale nonlinear programming.
\newblock \emph{Mathematical programming}, 106\penalty0 (1):\penalty0 25--57,
  2006.

\end{thebibliography}
 
\ECSwitch

\ECDisclaimer

\ECHead{Additional discussions, missing proofs, and formulations}

\section{Proof of Proposition~\ref{NP-Hard}}\label{prove NP-Hard}
This proof uses a construction similar to the proof of Theorem~1 in~\cite{hansen1990boolean}, reducing the partition problem to pricing under FML.
We give a polynomial-time reduction from \textsc{Partition}.
Given positive integers $c_1,\ldots,c_n$, the \textsc{Partition}
problem asks whether there exists a subset $J\subseteq[n]$ such that
\[
\sum_{j\in J}c_j=\sum_{j\in[n]\setminus J}c_j.
\]
Define
\[
S=\sum_{j=1}^n c_j,\qquad x_j=\frac{c_j}{S}.
\]
Thus, the instance is a yes-instance if and only if there exists
$J\subseteq[n]$ with $\sum_{j\in J}x_j=1/2$.

Consider two customer segments with weights $w_1=1/3$ and $w_2=2/3$.
For each product $j\in[n]$, set
\[
p_{j1}=3,\qquad p_{j2}=3+2x_j,
\]
and define the attraction values by
\[
v_{1j1}=1+(n+1)x_j,\qquad v_{1j2}=1,
\qquad
v_{2j1}=1+4(n+1)x_j,\qquad v_{2j2}=1.
\]
These attraction values are consistent with the assumed utility model.
Indeed, choosing
\[
\beta_{c j}
=\frac{\log v_{cj1}}{2x_j}>0,
\qquad
\alpha_{cj}=p_{j2}\beta_{cj},
\qquad c\in\{1,2\},
\]
gives
$v_{cjk}=\exp(\alpha_{cj}-\beta_{cj}p_{jk})$
for both $k\in\{1,2\}$.

Let $z_j\in\{0,1\}$ indicate whether the second price is selected,
and define
\[
t=\sum_{j:z_j=0}x_j.
\]
Since $\sum_{j=1}^n x_j=1$, direct substitution gives the segment
revenues
\[
R_1(z)=3-\frac{2}{(n+1)}+\frac{1}{(n+1)(1+t)},
\qquad
R_2(z)=3-\frac{1}{2(n+1)}-\frac{1}{2(n+1)(1+4t)}.
\]
Therefore, the expected revenue is
\[
\Phi(z)
=\frac13 R_1(z)+\frac23 R_2(z)
=3-\frac{1}{(n+1)}+\frac{1}{3(n+1)}
\left(\frac{1}{1+t}-\frac{1}{1+4t}\right).
\]
Set the decision threshold to
\[
K=3-\frac{8}{9(n+1)}.
\]
Then
\begin{align*}
\Phi(z)-K
&=\frac{1}{9(n+1)}
\left(\frac{3}{1+t}-\frac{3}{1+4t}-1\right)\\
&=-\frac{(2t-1)^2}{9(n+1)(1+t)(1+4t)}
\le 0.
\end{align*}
Because the denominator is positive, equality holds if and only if
$t=1/2$, which is equivalent to
\[
\sum_{j:z_j=0}c_j=\sum_{j:z_j=1}c_j.
\]
Hence,
\[
\exists z\in\{0,1\}^n:\ \Phi(z)\ge K
\quad\Longleftrightarrow\quad
\text{the given \textsc{Partition} instance is a yes-instance}.
\]
Thus, the decision version
of the pricing problem is NP-hard, and consequently its optimization
version is NP-hard.\Halmos

\section{Proof of Theorem~\ref{them:personalized}}
Since we can decide price for each segment independently,
it suffices to prove that the statement holds when there is only one segment $c\in M$.
Abbreviate 
\[D^c(\z) = 1 + \sum_{i\in N}(\delta^c_{i1}+\sum_{k=1}^{d_i-1}\delta^c_{i,k+1}z_{ik} )\]
With~\eqref{eq:unary}, the personalized pricing problem for segment $c$ can be written as:
\[\max_{(\rho^c,\y^c)\in  S^c} \rho^c\sum_{i\in N}\gamma^c_{i1} +\sum_{i\in N}\sum_{k=1}^{d_i-1} \gamma^c_{ik+1}y^c_{ik} \]
where 
\[ S^{c}:=
\left\{
\left(\frac{1}{D^c(\z)},\frac{\z}{D^c(\z)}\right)\;\middle|\;
\z\in\prod_{i\in N}\left(\Delta^{d_i-1}\cap\{0,1\}^{d_i-1}\right)
\right\}.\]
Since the objective is linear, maximizing it over \(S^c\) and over \(\conv(S^c)\) yields the same optimal value. Consider the projective transformation of $S^{c}$:
\[\Phi(S^c) = \left\{(D^c(\z),\z)\;\middle|\; \z\in\prod_{i\in N}\left(\Delta^{d_i-1}\cap\{0,1\}^{d_i-1}\right) \right\}\]
Since $D^c$ is affine and $\prod_{i\in N}\left(\Delta^{d_i-1}\right)$ is the convex hull of $\prod_{i\in N}\left(\Delta^{d_i-1}\cap\{0,1\}^{d_i-1}\right)$, we have 
\[\conv(\Phi(S^c)) = \left\{\bigl(D^c(\z),\z\bigr)\;\middle|\; \z\in\prod_{i\in N}\left(\Delta^{d_i-1}\right) \right\}\]
Theorem~1 of~\cite{he2025convexification} gives
\begin{align*}
\conv(S^{c})&=
\left\{(\rho^c,\y^c)
\;\middle|\;\rho^c\geq 0,\; (1,\y^c)\in
\rho^c\conv\bigl(\Phi(S^{c})\bigr)\right\}.
\end{align*}
Now we prove that when $\rho^c \geq 0$, $(1,\y^c)\in
\rho^c\conv\bigl(\Phi(S^{c})\bigr)$ if and only if $(\rho^c,\y^c)$ satisfies constraints~\eqref{eq:mip-1} and~\eqref{eq:mip-2}.
If $(\rho^c,\y^c)$ satisfies constraints~\eqref{eq:mip-1} and~\eqref{eq:mip-2}, $\rho^c = 0$ is infeasible, hence we can pick $\z = \y^c/\rho^c$. Then we have 
\[\z\in \prod_{i\in N}\left(\Delta^{d_i-1}\right),\quad \rho^cD^c(\z) =  \frac{\rho^c}{\rho^c} \left(\rho^c + \sum_{i\in N}\left(\delta^c_{i1}\rho^c+\sum_{k=1}^{d_i-1}\delta^c_{i,k+1}y^c_{ik} \right)\right)=1\]
Thus $(1,\y^c)\in \rho^c\conv\bigl(\Phi(S^{c})\bigr)$. Conversely, if $(1,\y^c)\in
\rho^c\conv\bigl(\Phi(S^{c})\bigr)$, there exists $\z\in \prod_{i\in N}\left(\Delta^{d_i-1}\right)$ such that 
\[
\begin{aligned}
    &\y^c = \rho^c \z\\
    &1 = \rho^c D^c(\z) = \rho^c + \sum_{i\in N}(\delta^c_{i1}\rho^c+\sum_{k=1}^{d_i-1}\delta^c_{i,k+1}\rho^cz_{ik} )
\end{aligned}
\]
Replacing $\rho^cz_{ik}$ by $y^c_{ik}$ in the second equality above yields 
\[1 = \rho^c + \sum_{i\in N}(\delta^c_{i1}\rho^c+\sum_{k=1}^{d_i-1}\delta^c_{i,k+1}y^c_{ik} )\]
Thus $(\rho^c,\y^c)$ satisfies~\eqref{eq:mip-1} and~\eqref{eq:mip-2}.
Now the personalized pricing problem for segment $c$ has the same optimal value as the following LP:
\[\max_{(\rho^c,\y^c)} \rho^c\sum_{i\in N}\gamma^c_{i1} +\sum_{i\in N}\sum_{k=1}^{d_i-1} \gamma^c_{ik+1}y^c_{ik},\quad s.t. ~\eqref{eq:mip-1} \text{ and }~\eqref{eq:mip-2} \]

Let $(\rho^c,\mathbf y^c)$ be an optimal extreme point solution of the LP. Since every extreme point of $\conv(S^c)$ belongs to $S^c$, the vector $\y^c/\rho^c\in \prod_{i\in N}\left(\Delta^{d_i-1}\cap\{0,1\}^{d_i-1}\right)$. Then, an optimal personalized price vector is recovered by
\[
x_i^c=p_{i1}+\sum_{k=1}^{d_i-1}(p_{i,k+1}-p_{ik})\frac{y_{ik}^c}{\rho^c},
\qquad i\in N.
\]\Halmos

The proof of Theorem~\ref{them:personalized} also follows from pseudolinearity of linear fractional functions and Charnes-Cooper transformation \citep{charnes1962programming}. However, our proof uses techniques that will be relevant for the coordinated price setting as well.

\section{Proof of Theorem~\ref{them:exactnessMIP}}
Let \(r^*_{\rm PC}\) and \(r^*\) denote the optimal values of~\eqref{eq:MIP} and~\eqref{price optimization}. To show \(r^*_{\rm PC}\geq r^*\), take any feasible price vector \(\x\), let \(\z\) be its unary representation, and set
\[
\rho^c:=\frac{1}{D^c(\z)},\qquad \y^c:=\rho^c \z,\qquad c\in M.
\]
Then~\eqref{eq:mip-1} and~\eqref{eq:mip-2} hold directly. Since \(\z\) is a vertex of the unary simplex, \(z_{ik}z_{ir}=z_{i,\max\{k,r\}}\) and
\[
\min\{z_{ik},z_{jr}\}=\max\{0,z_{ik}+z_{jr}-1\}=z_{ik}z_{jr}.
\]
Thus, both linking constraints reduce to \(z_{ik}=\rho^cD^c(\z)z_{ik}\) and hold at equality. The lifted objective equals the original fractional objective by construction. Hence every feasible solution of~\eqref{price optimization} has a feasible lift with the same value, proving \(r^*_{\rm PC}\geq r^*\).

Now we show the converse. Let $\underline D^c = (1 + \sum_{\ell \in N} \delta^c_{\ell 1})$, and $\overline D^c = (1+\sum_{\ell \in N}\sum_{t\in [d_\ell]} \delta^c_{\ell t+1})$, where $\underline D^c$ and $\overline D^c$ are the lower and upper bounds of $D^c(\z)$ over $\prod_{\ell \in N}\Delta^{d_\ell - 1}$ respectively.
We will show that the following relations are implied by \eqref{eq:MIP}:
\begin{subequations}
\begin{align}
&\overline D^c y^c_{ik} \ge z_{ik}\label{eq:McY1}\\
&\underline D^c y^c_{ik} - \underline D^c\rho^c + 1 \ge z_{ik}\label{eq:McY2}\\ 
&\underline D^c y^c_{ik} \le z_{ik}\label{eq:McY3}\\
&\overline D^c y^c_{ik} - \overline D^c\rho^c + 1 \le z_{ik}\label{eq:McY4}.
\end{align}
\end{subequations}
For \eqref{eq:McY1}, we have:
\begin{align*}
    z_{ik } &\le \left(1+\sum_{\ell \in N}\delta^c_{\ell 1}\right) y^c_{ik} + \sum_{t\in [d_i -1]}\delta^c_{i t+1}y^c_{i \max\{t,k\}} + \sum_{\ell \in N\backslash\{i\}}\sum_{t\in [d_\ell -1]}\delta^c_{\ell t+1}\max\{0, y^c_{ik} + y^c_{\ell t} - \rho^c\}\\
    & \le \left(1+\sum_{\ell \in N}\delta^c_{\ell 1}\right) y^c_{ik} = \overline D^cy^c_{ik}.
\end{align*}
where the first inequality follows from \eqref{eq:coordinating-1} and the second inequality follows since $\max\{0,y^c_{ik}+y^c_{\ell t} - \rho^c\}, y^c_{i \max\{t,k\}} \ge 0$ and $\delta^c_{\ell t+1}\le 0$. To see \eqref{eq:McY3}, consider the following inequalities:
\begin{align*}
    z_{ik } &\ge \left(1+\sum_{\ell \in N}\delta^c_{\ell 1}\right) y^c_{ik} + \sum_{t\in [d_i -1]}\delta^c_{i t+1}y^c_{i \max\{t,k\}} + \sum_{\ell \in N\backslash\{i\}}\sum_{t\in [d_\ell -1]}\delta^c_{\ell t+1}\min\{y^c_{\ell t}, y^c_{ik}\}\\
    & \ge \left(1+\sum_{\ell \in N}\delta^c_{\ell 1}\right) y^c_{ik} + \sum_{t\in [d_i -1]}\delta^c_{i t+1}y^c_{ik} + \sum_{\ell \in N\backslash\{i\}}\sum_{t\in [d_\ell -1]}\delta^c_{\ell t+1}y^c_{ik} = \underline D^c y^c_{ik},
\end{align*}
where the first inequality follows from \eqref{eq:coordinating-2} and the second inequality follows since, by \eqref{eq:mip-2}, $y^c_{ik}\ge y^c_{i\max\{k,t\}}$, $y^c_{ik}\ge \min\{y^c_{\ell t}, y^c_{ik}\}$, and $\delta^c_{\ell t+1}\le 0$. Now, we show \eqref{eq:McY2} using the following inequalities:
\begin{align*}
    z_{ik } &\le \left(1+\sum_{\ell \in N}\delta^c_{\ell 1}\right) y^c_{ik} + \sum_{t\in [d_i -1]}\delta^c_{i t+1}y^c_{i \max\{t,k\}} + \sum_{\ell \in N\backslash\{i\}}\sum_{t\in [d_\ell -1]}\delta^c_{\ell t+1}\max\{0, y^c_{ik} + y^c_{\ell t} - \rho^c\}\\
    &\le \left(1+\sum_{\ell \in N}\delta^c_{\ell 1}\right) y^c_{ik} + \sum_{t\in [d_i -1]}\delta^c_{i t+1}y^c_{i \max\{t,k\}} + \sum_{\ell \in N\backslash\{i\}}\sum_{t\in [d_\ell -1]}\delta^c_{\ell t+1}(y^c_{ik} + y^c_{\ell t} - \rho^c)\\
    & \le \left(1+\sum_{\ell \in N}\delta^c_{\ell 1}\right) (y^c_{ik}-\rho^c) + \sum_{t\in [d_i -1]}\delta^c_{i t+1}(y^c_{i \max\{t,k\}}-y^c_{it}) + \sum_{\ell \in N\backslash\{i\}}\sum_{t\in [d_\ell -1]}\delta^c_{\ell t+1}(y^c_{ik} - \rho^c) + 1\\
    & \le \left(1+\sum_{\ell \in N}\delta^c_{\ell 1}\right) (y^c_{ik}-\rho^c) + \sum_{t\in [d_i -1]}\delta^c_{i t+1}(y^c_{ik} - \rho^c) + \sum_{\ell \in N\backslash\{i\}}\sum_{t\in [d_\ell -1]}\delta^c_{\ell t+1}(y^c_{ik} - \rho^c) + 1\\
    &=\underline D^c(y^c_{ik} - \rho^c) + 1.
\end{align*}
Here, the first inequality follows from \eqref{eq:coordinating-1} and the second inequality follows since $\delta^c_{\ell t+1}\le 0$.
The third inequality follows from \eqref{eq:mip-1}. 
By \eqref{eq:mip-2}, it follows that 
$(y^c_{i \max\{t,k\}} - y^c_{ik}) - (y^c_{it} - \rho^c) = -y^c_{i \min\{t,k\}} + \rho^c \ge 0$.
Then, the fourth inequality follows since $\delta^c_{\ell t+1}\le 0$. Now, we show \eqref{eq:McY4} as follows:
\begin{align*}
    z_{ik } &\ge \left(1+\sum_{\ell \in N}\delta^c_{\ell 1}\right) y^c_{ik} + \sum_{t\in [d_i -1]}\delta^c_{i t+1}y^c_{i \max\{t,k\}} + \sum_{\ell \in N\backslash\{i\}}\sum_{t\in [d_\ell -1]}\delta^c_{\ell t+1}\min\{y^c_{\ell t}, y^c_{ik}\}\\
    &= \left(1+\sum_{\ell \in N}\delta^c_{\ell 1}\right) (y^c_{ik}-\rho^c) + 1 + {} \\
    &\qquad\qquad \sum_{t\in [d_i -1]}\delta^c_{i t+1}(y^c_{i \max\{t,k\}} - y^c_{it}) + \sum_{\ell \in N\backslash\{i\}}\sum_{t\in [d_\ell -1]}\delta^c_{\ell t+1}(\min\{y^c_{\ell t}, y^c_{ik}\} - y^c_{\ell t})\\
    &\ge \left(1+\sum_{\ell \in N}\delta^c_{\ell 1}\right) (y^c_{ik}-\rho^c) + 1,\\
    & = \overline D^c (y^c_{ik} - \rho^c) + 1
\end{align*}
where the first inequality follows from \eqref{eq:coordinating-2}, the first equality from \eqref{eq:mip-1}, and the inequality because $\delta^c_{\ell t+1}\le 0$, $y^c_{i\max\{t,k\}}-y^c_{it}\le 0$, and $\min\{y^c_{\ell t},y^c_{ik}\}-y^c_{\ell t} \le 0$.

Now, it follows from \eqref{eq:McY1} and \eqref{eq:McY3} that if $z_{ik}=0$ then $y^c_{ik}=0$. Similarly, it follows from \eqref{eq:McY2} and \eqref{eq:McY4} that if $z_{ik}=1$, $y^c_{ik}=\rho^c$. Then, it follows from \eqref{eq:mip-1} that $(\rho^c,\y^c)$ satisfy \eqref{eq:CC}. Therefore, $x_i = p_{i1} + \sum_{t\in [d_i-1]} p_{it+1}z_{it}$ is feasible to \eqref{price optimization} and its objective value matches that of $(\z,\rho^c,\y^c)$ in \eqref{eq:MIP}. It follows that $r^*_{\text{PC}} \le r^*$. Combining, $r^*_{\text{PC}} = r^*$.

\Halmos
\section{Proof of Theorem~\ref{them:indhull}}\label{app:indhull}
\def \B {\mathcal{B}}

For a set $S$ and a scalar $\rho \geq 0$, we denote its scaling by $\rho$ as $\rho \cdot S$ as $\{ \rho \x \mid \x \in S,\ \rho \geq 0 \}$. If $S$ is a polyhedron given by a system of linear inequalities $\BFA\x \leq \BFb$, its scaled set $\rho \cdot S$ admits an polyhedral representation $\{\x'\mid \BFA \x' \leq \BFb \rho, \rho\geq 0 \}$. 

\begin{lemma}\label{lemma:projective}
For each segment $c$ and item $i$, let $\B^c_i$ denote the convex hull of the following bilinear set:
\begin{equation}\label{eq:bilinear}
\Bigl\{\bigl(1, \z, \phi_i(\z_i) \cdot D^c(\z) \bigr)   \Bigm| \z_j \in \Delta^{d_j-1} \cap\{0,1\}^{d_j-1} \for j \in N \Bigr\}.
\end{equation}
Then,~\eqref{eq:hull} can be described as follows:
\begin{equation}\label{eq:projective}
    \Bigl\{(\mu^c_i,\rho^c, \y^c) \Bigm|  (\rho^c, \y^c,\mu^c_i) \in \rho^c \cdot \B^c_i,\ (\rho^c,\y^c) \text{ satisfies }\eqref{eq:mip-1}  \Bigr\}.
\end{equation}
\end{lemma} 
\noindent{\bf Proof of Lemma~\ref{lemma:projective}.} 
We use Theorem 2 of~\cite{he2025convexification} to prove this result. Specifically, we treat the set in~\eqref{eq:bilinear} as the set of base functions (denoted as $\mathcal{S}$ in~\cite{he2025convexification}). The denominator $D^c(\z)$ is strictly positive on $\Delta = \prod_{i=1}^n \Delta^{d_i-1}$ and is a linear combination of the first two arguments $(1,\z)$ of the set in~\eqref{eq:bilinear}. Moreover,  both $\B^c_i$ and~\eqref{eq:hull} are polytopes. Therefore, all the conditions of Theorem 2 in~\cite{he2025convexification} are satisfied, and~\eqref{eq:projective} describes~\eqref{eq:hull}. \Halmos

\begin{lemma}\label{lemma:common_simplex}
    For each segment $c$ and a pair of items $i$ and $j$, let $\B^c_{ij}$ denote the convex hull of 
    \begin{equation}\label{eq:common_simplex}
        \Bigl\{(\z_i, \z_j, \phi_i(\z_i) \cdot \ell^c_j(\z_j) )   \Bigm| \z_i \in \Delta_i \cap\{0,1\}^{d_i-1}, \z_j \in \Delta_j \cap\{0,1\}^{d_j-1} \Bigr\}. 
    \end{equation}
Suppose that  $\phi_i(\z_i) =  \sum_{k=1}^{d_i-1} \alpha_{ik}z_{ik}$.   Then, the convex hull of~\eqref{eq:bilinear}, denoted as $\B^c_{i}$, is given as follows:
    \[
    \Bigl\{(1,\z, \mu^c_i)  \Bigm| \mu^c_i = \phi_i(\z_i)(1+\delta^c_{i1}) + \sum_{k=1}^{d_i-1} \alpha_{ik} \sum_{r=1}^{d_i-1} \delta^c_{ir+1}z_{i, \max\{k,r\}} +  \sum_{j  \neq i } \mu^c_{ij},\ (\z_i, \z_j, \mu^c_{ij}) \in \B^c_{ij} \for j  \neq i \Bigr\}. 
    \]
\end{lemma}
\noindent{\bf Proof of Lemma~\ref{lemma:common_simplex}.} 
By using the relation $z_{ik}z_{ir} = z_{i,\max\{k,r\}}$, we obtain
\[
\phi_i(\z_i) \cdot \ell^c_i(\z_i) =\phi_i(\z_i)\delta^c_{i1} + \sum_{k=1}^{d_i-1} \alpha_{ik} \sum_{r=1}^{d_i-1} \delta^c_{ir+1}z_{i, \max\{k,r\}}.  
\]
For each $c$ and a pair of items $i,j$, let $\mathcal{S}^c_{ij}$ denote the set in~\eqref{eq:common_simplex}. Thus, the set of our interest~\eqref{eq:bilinear} can be obtained as an affine transformation $T$ of the following bilinear set
\begin{equation}\label{eq:multi-bilinear}
    \Bigl\{\bigl(\z, (\mu^c_{ij})_{j\neq i} \bigr) \Bigm| (\z_i,\z_j, \mu^c_{ij})  \in \mathcal{S}^c_{ij}  \for j \in N \setminus \{i\}  \Bigr\}.
\end{equation}
Since convexification commutes with affine transformation, it suffices to derive the convex hull of~\eqref{eq:multi-bilinear}.  To decompose this convex hull of~\eqref{eq:multi-bilinear} into the intersection of the convex hulls of the sets $\mathcal{S}^c_{ij}$, we invoke an established decomposition result from the literature~\citep{schrijver1983short,he2026discretenessconvexitypromotionplanning}. In particular, the projection of each $\mathcal{S}^c_{ij}$ onto the space of  the common unary variable $\z_i$ is $\Delta^{d_i-1} \cap \{0,1\}^{d_i}$ $j\neq i$, whose elements are affinely independent. Thus, by recursively using Lemma 7 of~\cite{he2026discretenessconvexitypromotionplanning}, we can establish that the convex hull of~\eqref{eq:multi-bilinear} is given by
\[
    \Bigl\{\bigl(\z, (\mu^c_{ij})_{j\neq i} \bigr) \Bigm| (\z_i,\z_j, \mu^c_{ij})  \in \B^c_{ij}  \for j \in N \setminus \{i\}  \Bigr\}.
\]
Combining this characterization with the affine transformation $T$ yields the claimed representation of $\B_i^c$, completing the proof.\Halmos
\begin{lemma}\label{lemma:common_triangulation}
    Let $f(\boldsymbol{u}) = \sum_{r =1}^m \alpha_ru_r$ where $\alpha_r \geq 0 $ and $\boldsymbol{u} \in \Delta^{m}$, and let $g(\boldsymbol{v}) = \sum_{t=1}^n \beta_t v_t$, where $\beta_t \leq 0$ and $\boldsymbol{v} \in \Delta^n$. Then, the convex hull of 
    \begin{equation}\label{eq:bilinear-graph-hull}
        \bigl\{( \boldsymbol{u}, \boldsymbol{v},\mu) \bigm| \mu = f(\boldsymbol{u})g(\boldsymbol{v}),\ \boldsymbol{u} \in \Delta^{m} \cap \{0,1\}^{m} \text{ and } \boldsymbol{v} \in \Delta^n \cap\{0,1\}^n \bigr\}
    \end{equation}
     is given by $\boldsymbol{u} \in \Delta^m$, $\boldsymbol{v} \in \Delta^n$ and 
    \[
     \sum_{r=1}^{m}\sum_{t=1}^{n}\alpha_{r}\beta_{t} \min \{u_{r},v_{t}\}  \leq   \mu \leq  \sum_{r=1}^{m}\sum_{t=1}^{n}\alpha_r \beta_t  \max\{0,u_{r} + v_{t} -1\}.
    \]
\end{lemma}
\noindent{\bf Proof of Lemma~\ref{lemma:common_triangulation}.} Let $g: \bigl(\Delta^m \cap\{0,1\}^m\bigr) \times \bigl(\Delta^n \cap \{0,1\}^n \bigr) \to \R$ defined as:
\[
g(\boldsymbol{u},\boldsymbol{v}) := f(\boldsymbol{u})g(\boldsymbol{v}) = \sum_{r=1}^m\sum_{t=1}^n\alpha_r\beta_t u_rv_t  
\]
To derive the convex hull of its graph, that is the set in~\eqref{eq:bilinear-graph-hull}, it suffices to derive its convex envelope and concave envelope. The convex (resp. concave) envelope of $g$, denoted as $\conv(g)$ (resp. $\conc(g)$), is the largest (resp. smallest) convex under-estimator (resp. concave over-estimator) of $g$ over $ \Delta:= \Delta^m \times \Delta^n$ (see Definition 2.1 of~\cite{tawarmalani2013explicit}).

First, we derive the convex envelope of $g$ by individually constructing that of $g_{rt}(\boldsymbol{u},\boldsymbol{v})= \alpha_r\beta_t u_rv_t$. For each $r$ and $t$, we have $g_{rt}(\boldsymbol{u},\boldsymbol{v}) \geq \alpha_r\beta_t\min \{u_r,v_t\} =: L_{rt}(\boldsymbol{u},\boldsymbol{v})$ for binary $\boldsymbol{u}$ and $\boldsymbol{v}$. Since $L_{rt}(\cdot)$ is convex, we conclude that $\conv(g_{rt})(\boldsymbol{u},\boldsymbol{v}) \geq L_{rt}(\boldsymbol{u},\boldsymbol{v})$ for $(\boldsymbol{u},\boldsymbol{v}) \in \Delta^m \times \Delta^n$. To establish the reverse direction, we introduce a subdivision of $\Delta^m \times \Delta^n$. For a permutation $\sigma$ of the index set $\{1,2 \ldots, m+n\}$, define 
\[
\begin{aligned}
&K_\sigma:=\Bigl\{\boldsymbol{w} := (\boldsymbol{u},\boldsymbol{v}) \in \Delta^m \times \Delta^n \Bigm| w_{\sigma(1)} \geq \cdots \geq w_{\sigma(m+n)}  \Bigr\},\\ 
&\mathcal{K} := \bigl\{K_\sigma \bigm| \sigma^{-1}(1) \geq \cdots \geq \sigma^{-1}(m) \text{ and }\sigma^{-1}(m+1) \geq \cdots \geq \sigma^{-1}(m+n) \bigr\}.    
\end{aligned}
\]
It follows from Section 2.2 in~\cite{he2022tractable} or Theorem 6.2.13 in~\cite{de2010triangulations}  that
 $\mathcal{K}$ is a triangulation of $\Delta^m \times \Delta^n$, \textit{e.g.}, $\Delta^m \times \Delta^n = \cup_{K \in \mathcal{K}} K$ and for two $K',K'' \in \mathcal{K}$, $K'\cap K''$ is a face of both $K'$ and $K''$. Now, for a given $\bar{\boldsymbol{w}}:=(\bar{\boldsymbol{u}},\bar{\boldsymbol{v}}) \in \Delta$, there exists $\bar{K} \in \mathcal{K}$ such that  $(\bar{\boldsymbol{u}},\bar{\boldsymbol{v}}) \in \bar{K}$. More specifically, there exists a vector of convex multipliers  $\lambda_{\boldsymbol{w}}$ such that $\bar{\boldsymbol{w}} = \sum_{\boldsymbol{w} \in \vertex (\bar{K})} \lambda_{\boldsymbol{w}}  \boldsymbol{w}$. Moreover, for each $r$ and $t$, $L_{rt}(\cdot)$ is linear over $\bar{K}$. Therefore, we conclude that
 \[
 \begin{aligned}
     \conv(g)(\bar{\boldsymbol{w}}) &\geq \sum_{r,t} \conv(g_{rt})(\bar{\boldsymbol{w}}) \geq \sum_{r,t}L_{rt}(\bar{\boldsymbol{w}}) \\
     &= \sum_{\boldsymbol{w} \in \vertex(\bar{K})} \lambda_{\boldsymbol{w}}\sum_{r,t} L_{rt}(\boldsymbol{w}) = \sum_{\boldsymbol{w} \in \vertex(\bar{K})} \lambda_{\boldsymbol{w}} g(\boldsymbol{w}) \geq \conv(g)(\bar{\boldsymbol{w}}),
 \end{aligned} 
 \]
where the first inequality holds since the function on the right-hand-side is a convex under-estimator of $g$, the second inequality has been established above, the first equality holds since $L_{rt}$  is linear over $\bar{K}$ for each $r$ and $t$, the second equality holds since $L_{rt}(\boldsymbol{u},\boldsymbol{v}) = g_{rt}(\boldsymbol{u},\boldsymbol{v})$ for binary $(\boldsymbol{u},\boldsymbol{v})$ and $\vertex(\bar{K})$ are binary, and the last inequality holds by the convexity of $\conv(g)$. Therefore, inequalities hold throughout. In particular, $\conv(g)(\boldsymbol{u},\boldsymbol{v}) = \sum_{r=1}^m\sum_{t=1}^nL_{rt}(\boldsymbol{u},\boldsymbol{v})$, and $\conv(g_{rt})(\boldsymbol{u},\boldsymbol{v}) = L_{rt}(\boldsymbol{u},\boldsymbol{v})$ for all $(\boldsymbol{u},\boldsymbol{v}) \in \Delta$.

Next, we derive the concave envelope of $g$. Consider a switched function $g'$ of $g$, defined as follows: $g'(\boldsymbol{u}, \boldsymbol{v}')= g(\boldsymbol{u}, \mathbf{1} - \boldsymbol{v'}) = \sum_{r=1}^m \sum_{t=1}^n\alpha_r\beta_tu_r(1-v'_t) = -\sum_{r=1}^m\sum_{t=1}^n\alpha_r\beta_tu_rv'_t + \sum_{r=1}^m\sum_{t=1}^n\alpha_r\beta_tu_r$. Then, 
\[
\begin{aligned}
\conc(g')(\boldsymbol{u},\boldsymbol{v}') &= -\sum_{r=1}^m\sum_{t=1}^n\alpha_r\beta_t \Bigl( \min \{u_r,v'_t\} + u_r\Bigr) \\
& = \sum_{r=1}^m\sum_{t=1}^n\alpha_r\beta_t \Bigl( \max \{0,u_r- v'_t\} \Bigr).
\end{aligned}
\]
The proof is complete by using $\conc(g)(\boldsymbol{u},\boldsymbol{v}) = \conc(g')(\boldsymbol{u}, \mathbf{1} -  \boldsymbol{v})$. \Halmos

\noindent{\bf Proof of Theorem~\ref{them:indhull}.} 
Since $\phi_i(\z_i)$ is a non-decreasing linear function,  it can be represented as   $\phi_i(\z_i) =  \sum_{k=1}^{d_i-1} \alpha_{ik}z_{ik}$ where $\alpha_{ik} \geq 0$.  Recall that $\ell^c_j(\z_j):= \delta^c_{j1} +  \sum_{k =1 }^{d_j-1}\delta^c_{jk+1}z_{jk}$. Since the incremental attraction values $\delta^c_{jk}$ are non-positive, it follows from Lemma~\ref{lemma:common_triangulation} that, for $j \neq i$, $\B^c_{ij}$ can be described by the following inequalities:
\[
\begin{aligned}
\mu^c_{ij} &\geq \delta^c_{j1} \phi_i(\z_i) +  \sum_{k=1}^{d_i-1} \alpha_{ik} \sum_{r =1}^{d_j-1} \delta^c_{jr+1}\min\{ z_{ik},z_{jr}\}   \\
\mu^c_{ij} &\leq  \delta^c_{j1} \phi_i(\z_i) +  \sum_{k=1}^{d_i-1} \alpha_{ik} \sum_{r=1}^{d_j-1}\delta^c_{jr+1} \max\{0, z_{ik} + z_{jr} - 1\} .    
\end{aligned}
\]
Therefore, by Lemma~\ref{lemma:common_simplex}, these inequalities allows us to describe $\B^c_i$. In particular, we obtain the following expression of $\rho^c \cdot \B^c_i$:
\[
\begin{aligned}
&(\rho^c, \y^c) \text{ satisfies } \eqref{eq:mip-2}  \\
&\mu^c_i \leq   \sum_{k=1}^{d_i-1} \alpha_{ik} \biggl( \Bigl(1+ \sum_{j \in N} \delta^c_{j1}\Bigr)y^c_{ik} +  \sum_{r=1}^{d_i-1} \delta^c_{ir+1}y^c_{i, \max\{k,r\}} + \sum_{j \neq i} \sum_{r=1}^{d_j-1}\delta^c_{jr+1} \max\{0, y^c_{ik} + y^c_{jr} - \rho^c\}\biggr) \\
&\mu^c_i \geq  \sum_{k=1}^{d_i-1} \alpha_{ik} \biggl( \Bigl(1+ \sum_{j \in N} \delta^c_{j1}\Bigr)y^c_{ik} +   \sum_{r=1}^{d_i-1} \delta^c_{ir+1}y^c_{i, \max\{k,r\}} + \sum_{j \neq i} \sum_{r=1}^{d_j-1}\delta^c_{jr+1} \min\{y^c_{ik},y^c_{jr}\} \biggr), 
\end{aligned}
\]
By Lemma~\ref{lemma:projective}, these inequalities, joint with~\eqref{eq:mip-1}, provide an explicit expression of~\eqref{eq:hull}. This expression is equivalent to~\eqref{eq:indhull_rep} since the coefficients $\alpha_{ik}$ are nonnegative. \Halmos

\section{Reversed inequalities in Example~\ref{ex:onehot-vs-unary}}\label{ex:reverse}

\paragraph{Reversed inequalities on $\mu$} If we consider the set where 
\begin{align*}
    \mu_1 &\ge 2\lambda_{13}\lambda_{21}+\lambda_{13}\lambda_{22}
        = \lambda_{13}\bigl(2\lambda_{21}+\lambda_{22}\bigr),\\[2pt]
    \mu_2 &\ge 2\lambda_{12}\lambda_{21}+\lambda_{12}\lambda_{22}
        +2\lambda_{13}\lambda_{21}+\lambda_{13}\lambda_{22}
        = (\lambda_{12}+\lambda_{13})\bigl(2\lambda_{21}+\lambda_{22}\bigr),
\end{align*}
then the RLT set analogous to Set~1 with reverse linking inequalities satisfies 
\begin{equation*}
2\lambda_{12} - \lambda_{22} + 2\lambda_{13} - 2\lambda_{23} - \mu_2 \le 0.
\end{equation*}
However, the one-hot encoding base Set 2 type relaxation does not identify infeasibility of $\lprod_{kt} = 0$ for all $k,t$, $\blambda_1=(\tfrac12,0,\tfrac12)$, $\blambda_2=(\tfrac12,\tfrac12,0)$, $\mu_1=3$, and $\mu_2=0$. The unary encoding based relaxation proves this inequality using the following reduction:
\begin{equation*}
\begin{split}
&(\zprod_{11}\le z_{21}) + (\zprod_{12}\le z_{22}) + (2z_{11}-\zprod_{11}-\zprod_{12}\le \mu_2)\\
&\qquad+ (z_{22}=\lambda_{23}) + (z_{21} = \lambda_{22} + \lambda_{23}) + 2(z_{11} = \lambda_{12} + \lambda_{13}) 
\end{split}
\end{equation*}
\section{Proof of Proposition~\ref{prop:onehot-convexhull}}
We will show that \eqref{eq:onehot-convexhull} is the convex hull of points defined as:
\begin{equation*}
\begin{split}
    X = &\Bigl\{\bigl(x_i, \blambda_i, \rho^c, \s^c, \boldsymbol{\Gamma}^c\bigr) \Bigm| \boldsymbol{\Gamma}^c_j = \rho^c\blambda_i\blambda_j^\intercal, \s^c_j = \rho^c\blambda_j, \rho^c = \frac{1}{1+\sum_{j\in N}\exp(a_{cj}-b_{cj}x_j)},\\
    &\qquad (x_j, \blambda_j) \in \eqref{eq:full} \for j\in N\Bigr\}.
\end{split}
\end{equation*}
First, observe that $\rho^c > 0$. If $\rho^c=0$, it follows that $\s^c_i=0$ which implies $\Gamma^c_{jkt} = 0$ for all $j,k,t$, and $s^c_{jt} = 0$ for all $j$. Then, $\rho^c + \sum_{j}\sum_{t=1}^{d_j} v^c_{jt}s^c_{jt} = 0 \ne 1$.

For convexification purposes, we will first ignore \eqref{eq:onehot-definex} and \eqref{eq:onehot-defn} since convexification commutes with affine transformation. Now, observe that the constraints \eqref{eq:onehot-beginflow}--\eqref{eq:onehot-nonnegative} define a convex cone and \eqref{eq:onehot-normalize} can be treated as a normalization hyperplane. Therefore, it suffices to study the extreme points of \eqref{eq:onehot-beginflow}-\eqref{eq:onehot-nonnegative} with $\rho^c = 1$ and scale to satisfy \eqref{eq:onehot-normalize}. By Simplex Lemma~\citep{schrijver1983short,kim2022reciprocity}, it follows easily that the extreme rays of the cone are such that, for each $j\in N$, there is a $t_j$ such that $s^c_{jt_j} = \rho^c$ and $s^c_{jt} = 0$ for $t\ne t_j$. For $j\in N\backslash\{i\}$, $\Gamma^c_{jt_it_j} = \rho^c$ with $\Gamma^c_{jkt} = 0$ unless $k=t_i$ and $t=t_j$. Then, by \eqref{eq:onehot-normalize}, it follows that an extreme point satisfies $1 = \rho^c\left(1 + \sum_j v^c_{jt_j}\right)$, hence $\rho^c = (1+\sum_j v^c_{jt_j})^{-1}$. Now, consider two cases: $t_i=k$ or $t_i\ne k$. In the first case, $\lambda_{ik} = \rho^c\left(1 + \sum_{j\in N}v^c_{jt_j}\right) = 1$. In the second case, $s^c_{ik} = \Gamma^c_{jkt} = 0$ for all $j\ne i$, $t$, hence $\lambda_{ik}=0$. Finally, \eqref{eq:onehot-definex} just defines $x_i$ using $\blambda_i$. For $j\in N\backslash\{i\}$, define $\blambda_{j} = e_{t_j}$ and $x_j = p_{jt_j}$. Then, it is easy to verify that these lifted extreme points are feasible to the set defining $X$. It follows that \eqref{eq:onehot-convexhull} is included in $\conv(X)$. The opposite inclusion follows since the validity of the constraints can be verified directly from the definitions of $\Gamma^c_{jkt}$ and $s^c_{jt}$. The proposition then follows since it considers projections of $X$, and convexification commutes with projection.
\Halmos
\section{Proof of Proposition~\ref{prop:PD-dominates-FP}}
For each $(\z,\boldsymbol{\rho},\y)$ feasible to the LP relaxation of \eqref{eq:MIP}, we will construct $(\z,\bold d,\boldsymbol{\mu})$ feasible to the LP relaxation of \eqref{fp} with the same objective function value. Let
\begin{equation}\label{eq:factprog}
\begin{aligned}
&N^c(\z) := \sum_{i\in N}\left(\gamma^c_{i1}+\sum_{k=2}^{d_i}\gamma^c_{ik}z_{i,k-1}\right), &D^c(\z) := 1+\sum_{i\in N}\left(\delta^c_{i1}+\sum_{k=2}^{d_i}\delta^c_{ik}z_{i,k-1}\right),\\
&\mathcal M^c(\y^c) := \sum_{i\in N} \biggl(\gamma^c_{i1} \cdot \rho^c +\sum_{k=2}^{d_i}\gamma^c_{ik} \cdot y^c_{ik-1}\biggr).
\end{aligned}
\end{equation}
Given $(\z,\boldsymbol\rho,\y)$, the point $(\z,\bold d,\boldsymbol \mu)$ uses the same $\z$. We set the remaining coordinates $\bold d=(D^c(\z))_{c\in M}$ and $\boldsymbol\mu = (\mathcal M^c(\y^c))_{c\in M}$. The objective of \eqref{eq:MIP}, $\sum_{c\in M}\omega^c M^c(\y^c)$, matches the objective value of \eqref{fp} at the constructed solution. 

By construction, $d^c=D^c(\z)$, and since $\z\in\Delta$, we have $\underline d^c\le d^c\le\overline d^c$. To verify the bounds on $\mu^c$, constraints~\eqref{eq:mip-1} and~\eqref{eq:mip-2} imply $\rho^c>0$.
Thus, we can define $\mathbf u^c:=\y^c/\rho^c$. 
By~\eqref{eq:mip-2}, $\mathbf u^c\in\Delta$, and substituting $\y^c$ in~\eqref{eq:mip-1} by $\rho^c\mathbf u^c$ 
gives $\rho^c D^c(\mathbf u^c)=1$. Hence,
\[
\mu^c=\mathcal M^c(\y^c)
=\rho^c N^c(\mathbf u^c)
=\frac{N^c(\mathbf u^c)}{D^c(\mathbf u^c)}.
\]
Since $\mathbf u^c$ lies in the unit box, $\underline\mu^c\le\mu^c\le\overline\mu^c$.

It remains to verify that the constructed solution satisfies \eqref{FPMCC}. Using \eqref{eq:factprog}, we have
\begin{align}
&N^c(\z)-\underline d^c\mu^c
=(1-\underline d^c\rho^c)\sum_{j\in N}\gamma^c_{j1}
+\sum_{j\in N}\sum_{k=2}^{d_j}\gamma^c_{jk}
\bigl(z_{j,k-1}-\underline d^c y^c_{j,k-1}\bigr)\label{eq:nc-dmu}\\
&D^c(\z)-\underline d^c=\left(1+\sum_{j\in N}\delta^c_{j1}\right)(1-\underline d^c\rho^c)
+\sum_{j\in N}\sum_{k=2}^{d_j}\delta^c_{jk}(z_{j,k-1}-\underline d^cy^c_{j,k-1}),\label{eq:dc-lower}
\end{align}
where we have used \eqref{eq:mip-1} to write the second relation. By Theorem~\ref{them:indhull}, the LP relaxation of \eqref{eq:MIP} includes the convex hull of $(z_{jk},\rho^c,y^c_{jk})$. Therefore, \eqref{eq:MIP} satisfies:
\begin{equation}\label{eq:thmimpliedMC}
\underline d_cy^c_{jk}\le z_{jk}\le\overline d_cy^c_{jk},\qquad
\underline d_c(\rho^c-y^c_{jk})\le1-z_{jk}\le\overline d_c(\rho^c-y^c_{jk}).
\end{equation}
This implies 
\begin{equation}\label{eq:squish}
0\le z_{jk}-\underline{d}^cy^c_{jk} \le 1-\underline{d}^c\rho^c
\end{equation}
for $k\in [d_j-1]$. Since $\delta^c_{jk}\le 0$, it follows from \eqref{eq:squish} that
\begin{equation*}
    D^c(\z)-\underline{d}^c \ge \left(1+\sum_{j\in N}\left(\delta^c_{j1} +\sum_{k=1}^{d_j}\delta^c_{jk}\right)\right)(1-\underline{d}^c\rho^c) = \exp(a_{cj}-b_{cj}p_{jd_j})(1-\underline{d}^c\rho^c)\ge 0.
\end{equation*}
If $D^c(\z)-\underline{d}^c = 0$, the equality holds throughout. Since $\exp(a_{cj}-b_{cj}p_{jd_j}) > 0$, it follows that $(1-\underline{d}^c\rho^c) = 0$, and by \eqref{eq:squish}, $z_{jk}-\underline{d}^cy^c_{jk}=0$. This implies that $N^c(\z)-\underline{d}^c\mu^c = 0$. It follows that 
\begin{equation*}
\underline{\mu}^c(D^c(\z)-\underline{d}^c) = N^c(\z) - \underline{\mu}^cd^c = \overline{\mu}^c(D^c(\z)-\underline{d}^c).
\end{equation*}
On the other hand, if $D^c(\z)- \underline{d}^c > 0$, we have 
\begin{equation*}
\underline\mu^c\le \frac{N^c(\z) - \underline d^c\mu^c}{D^c(\z) - \underline{d}^c} = 
\frac{\sum_{j\in N}\gamma^c_{j1}
+\sum_{j\in N}\sum_{k=2}^{d_j}\gamma^c_{jk}p^+_{j,k-1}}{\sum_{j\in N}\delta^c_{j1}
+\sum_{j\in N}\sum_{k=2}^{d_j}\delta^c_{jk}p^+_{j,k-1}} = \frac{N^c(\bold p^+)}{D^c(\bold p^+)}\le \overline{\mu}^c,
\end{equation*}
where $p^+_{jk} := \frac{z_{jk}-\underline{d}^cy^c_{jk}}{1-\underline{d}^c\rho^c}$. The inequalities follow since $0\le p^+_{jk}\le 1$ by \eqref{eq:squish} and $\underline \mu^c$ and $\overline\mu^c$ are bounds of $\frac{N^c(\z)}{D^c(\z)}$ over $[0,1]^{\sum_{j\in N}d_{j-1}}$. This shows the validity of two out of four inequalities in \eqref{FPMCC}. 

Now, we consider the two remaining inequalities.
\begin{align}
&\overline d^c\mu^c - N^c(\z)
=(\overline d^c\rho^c - 1)\sum_{j \in N}\gamma^c_{j1}
+\sum_{j\in N}\sum_{k=2}^{d_j}\gamma^c_{jk}
\bigl(\overline d^c y^c_{j,k-1} - z_{j,k-1}\bigr)\label{eq:dmu-nc}\\
&\overline d^c - D^c(\z)=\left(1+\sum_{j\in N}\delta^c_{j1}\right)(\overline d^c\rho^c-1)
+\sum_{j\in N}\sum_{k=2}^{d_j}\delta^c_{jk}(\overline d^cy^c_{j,k-1}-z_{j,k-1}),\label{eq:dc-upper}
\end{align}
By \eqref{eq:thmimpliedMC}, we have
\begin{equation}\label{eq:squish2}
0\le \overline{d}^cy^c_{jk} - z_{jk}\le \overline{d}^c\rho^c - 1
\end{equation}
Since $\delta^c_{jk}\le 0$, it follows from \eqref{eq:squish2} that:
\begin{equation*}
    \overline{d}^c - D^c(\z) \ge \left(1+\sum_{j\in N}\left(\delta^c_{j1} +\sum_{k=1}^{d_j}\delta^c_{jk}\right)\right)(\overline{d}^c\rho^c-1) = \exp(a_{cj}-b_{cj}p_{jd_j})(\overline{d}^c\rho^c-1)\ge 0.
\end{equation*}
Then, using an argument similar to above, it follows that if $\overline{d}^c-D^c(\z) = 0$, we have:
\begin{equation*}
\underline{\mu}^c(\overline{d}^c - D^c(\z)) = \overline{\mu}^cd^c - N^c(\z) = \overline{\mu}^c(\overline{d}^c - D^c(\z)).
\end{equation*}
If $\overline{d}^c-D^c(\z) > 0$, we have:
\begin{equation*}
\underline\mu^c\le \frac{\overline d^c\mu^c - N^c(\z)}{\overline{d}^c - D^c(\z)} = 
\frac{\sum_{j\in N}\gamma^c_{j1}
+\sum_{j\in N}\sum_{k=2}^{d_j}\gamma^c_{jk}p^-_{j,k-1}}{\sum_{j\in N}\delta^c_{j1}
+\sum_{j\in N}\sum_{k=2}^{d_j}\delta^c_{jk}p^-_{j,k-1}} = \frac{N^c(\bold p^-)}{D^c(\bold p^-)}\le \overline{\mu}^c,
\end{equation*}
where $p^-_{jk} := \frac{\overline{d}^cy^c_{jk}-z_{jk}}{\overline{d}^c\rho^c-1}$. The inequalities follow since $0\le p^-_{jk}\le 1$ by \eqref{eq:squish2} and $\underline \mu^c$ and $\overline\mu^c$ are bounds of $\frac{N^c(\z)}{D^c(\z)}$ over $[0,1]^{\sum_{j\in N}d_{j-1}}$. This shows the validity of the remaining inequalities in \eqref{FPMCC}. 
\Halmos
\section{Computational Benefit of Unary formulation over One-Hot Approach}\label{sec:example_projection}
\begin{example}\label{separation efficiency comparison}
Consider $n$ items, each with $d$ admissible price levels, and let $\blambda_i=(\lambda_{i1},\ldots,\lambda_{id})\in\operatorname{vert}(\Sigma^d)$ denote the one-hot encoding for item $i\in[n]$. For each instance, we randomly generate coefficients satisfying $c_{i1}\geq c_{i2}\geq\cdots\geq c_{id}\geq0$ for $i\in\{2,\ldots,n\}$. We study separation over the convex hull of the graph
\[
\mu=\left(\sum_{j=\lceil d/2\rceil}^{d}\lambda_{1j}\right)\left(\sum_{i=2}^{n}\sum_{j=1}^{d}c_{ij}\lambda_{ij}\right).
\]

For each instance, we randomly generate candidate points with $\blambda_i\in\Sigma^d$ for $i\in[n]$ and $\mu\in[0,\sum_{i=2}^{n}c_{i1}]$. We compare LP-based separation using the one-hot formulation with RLT balance constraints against separation in unary coordinates. Specifically, we transform each candidate point using
\[
z_{i,k-1}=\sum_{j=k}^{d}\lambda_{ij},
\qquad i\in[n],\quad k\in\{2,\ldots,d\},
\]
and evaluate the lower and upper envelopes obtained from the termwise McCormick inequalities, without solving an LP. Both procedures therefore test the same candidate points in their respective coordinate systems.
Using $(n,d)=(10,10)$ as the baseline, we examine how separation time changes as $n$ increases with $d=10$ fixed, and as $d$ increases with $n=10$ fixed in Figure~\ref{fig:binarization_comp_speedup}.
\begin{figure}
    \centering
    \includegraphics[width=1\linewidth]{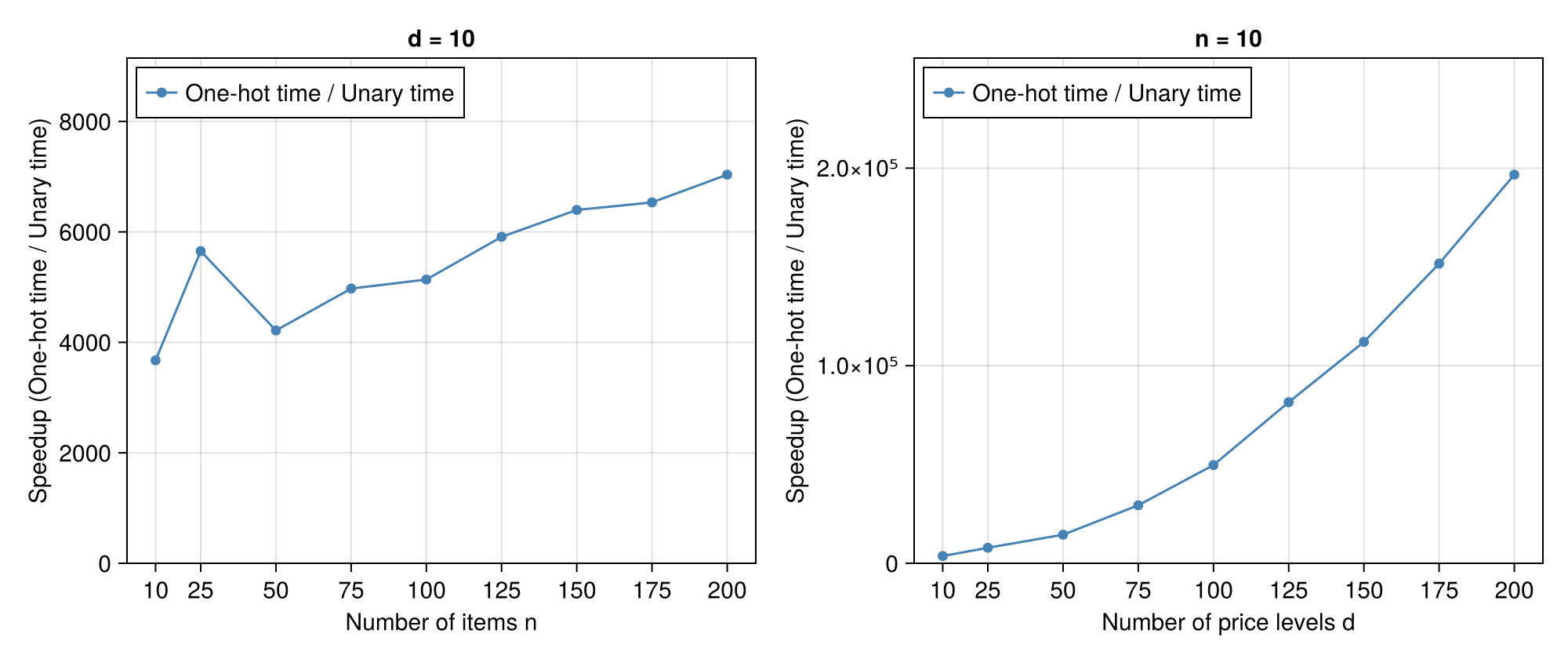}
    \caption{Speedup of unary separation over LP-based one-hot separation as the number of items $n$ increases with $d=10$ (left) and the number of price levels $d$ increases with $n=10$ (right). Speedup is measured as the ratio of one-hot to unary separation time.}
    \label{fig:binarization_comp_speedup}
\end{figure}

\end{example}

\section{Proof of Theorem~\ref{hierarchy}}

First, we present a lemma to show that flow polytopes are indeed relaxations of constraint~\eqref{eq:highorderScaleF}, and yield a convex hull description at the $n^{\text{th}}$ level. 

\begin{lemma}\label{lemma:multilinearhull}
 For $\tau \geq 2$, we have $W^c_\emptyset \cdot \mathcal{M}^\tau \subseteq L \bigl(P^\tau(W^c_{\emptyset})\bigr)$. Moreover, $\conv(W^c_\emptyset \cdot \mathcal{M}^n) = L\bigl(P^n(W^c_\emptyset)\bigr)$. 
\end{lemma}
\textbf{Proof of Lemma~\ref{lemma:multilinearhull}}
Let $\tau \geq 2$. To show $W^c_\emptyset \cdot \mathcal M^{ \tau} \subseteq L\bigl(P(W^c_\emptyset)\bigr)$, we consider  $\boldsymbol{W}^c\in W^c_\emptyset\mathcal M^{\tau}$. By definition, there is $\z_i \in \Delta_i \cap \{0,1\}^{d_i-1}$ for $i \in N$ such that $W^c_\alpha=W^c_\emptyset\prod_{(i,r)\in\alpha}z_{ir}$ for $\alpha \in \mathcal{A}$ such that $\vert \alpha \vert \leq \tau$.  Suppose for each $i$, there exists $\ell_i\in [d_i]$ such that $z_{ir} = 1$ for $r<\ell_i$ and $z_{ir}=0$ otherwise. For $i \in N$, define the corresponding one-hot variable by $\lambda_{ik} = 1$ if $k = \ell_i$ and $\lambda_{ik} = 0$ otherwise. For each $\bar{\alpha} \in \bar{\mathcal{A}}$, let 
\[
\Gamma^c_{\bar{\alpha}} = W^c_{\emptyset} \prod_{(i,k) \in \bar{\alpha}} \lambda_{ik},
\]
where the product over the empty index set is $1$. Thus, $\Gamma^c_{\emptyset} = W^c_{\emptyset}$ and $\boldsymbol{\Gamma}^c \geq 0$. To verify the flow conservation equations, consider any
$\bar\alpha\in\bar{\mathcal A}$ and $i\in N$ such that
$|\bar\alpha|<\tau$ and
$\bar\alpha\cap(\{i\}\times[d_i])=\emptyset$. We have
\[
\sum_{k=1}^{d_i} \Gamma^c_{\bar{\alpha} \cup\{(i,k)\} } = W^c_{\emptyset} \prod_{(j,r) \in \bar{\alpha}} \lambda_{jr} \sum_{k=1}^{d_i} \lambda_{ik}  =  W^c_{\emptyset} \prod_{(j,r) \in \bar{\alpha}} \lambda_{jr} = \Gamma^c_{\bar{\alpha}},
\]
where the second equality follows from
$\sum_{k=1}^{d_i}\lambda_{ik}=1$.
Therefore, $\boldsymbol\Gamma^c\in P^\tau(W^c_\emptyset)$. To complete the proof, we show $\boldsymbol{W}^c = L(\boldsymbol{\Gamma}^c)$. For $\alpha = \{(i,r_i) : i \in I\} \in \mathcal{A}$ with $\vert \alpha \vert \leq \tau$, we observe that 
\[
L(\boldsymbol{\Gamma}^c)_\alpha = \sum_{\bar{\alpha} \in \{(i,\ell_i) \mid i \in I,r_i< \ell_i \leq d_i\}} \Gamma^c_{\bar{\alpha}} = W^c_{\emptyset} \prod_{i \in I}\biggl(\sum_{\ell_i = r_i+1}^{d_i} \lambda_{i\ell_i}\biggr) = W^c_{\emptyset} \prod_{i \in I} z_{ir_i} = W^c_{\alpha}, 
\]
where the second equality holds by the construction of $\boldsymbol{\Gamma}^c$, and the third equality follows from the definition of the one-hot variables and unary variables. 

We next prove the convex hull characterization under $\tau = n$. Since $P^n(W^c_\emptyset)$ is convex and $L$ is linear, the first statement implies that $\conv(W^c_{\emptyset} \cdot \mathcal{M}^n) \subseteq L\bigl(P^n(W^c_{\emptyset})\bigr)$. To prove the reverse inclusion, we consider  $\boldsymbol \Gamma^c\in P^n(W^c_\emptyset)$. For each complete  path $ \boldsymbol{\ell} \in \prod_{i\in N}[d_i]$, define its weight as 
\[
\theta^{\boldsymbol{\ell}} = \Gamma^c_{\{(i,\ell_i) : i \in N\} }.
\]
Clearly, $\theta^{\boldsymbol{\ell}} \geq 0$. Moreover, the flow conservation equations imply that for $\bar{\alpha} \in \bar{\mathcal{A}}$,
\begin{equation}\label{eq:flow-distribution}
    \Gamma^c_{\bar{\alpha}} = \sum_{\substack{\boldsymbol{\ell} \in\prod_{i\in N}[d_i]\\
\ell_i=k_i\ \forall (i,k_i)\in\bar\alpha}} \theta^{\boldsymbol{\ell}}. 
\end{equation}
Indeed, this equality holds trivially when $\vert \bar{\alpha} \vert  =n$. For any partial path, $\bar{\alpha}$ with $\vert \bar{\alpha} \vert  <n$, it follows by repeatedly using the flow conservation equations to the products whose layers are not yet represented in $\bar{\alpha}$. In particular, taking $\bar{\alpha} = \emptyset$ yields   $ \Gamma^c_{\emptyset} = \sum_{\boldsymbol{\ell} \in \prod_{i \in N}[d_i]} \theta^{\boldsymbol{\ell}} $. For each complete  path $\boldsymbol{\ell} \in \prod_{i\in N}[d_i]$, we also introduce $\z^{\boldsymbol{\ell}}$ so that for $i \in N$
\[
z^\ell_{ik} = 1 \quad \for k < \ell_i \quad  \text{ and } \quad z^\ell_{ik} = 0 \quad \for  k \geq \ell_i. 
\]
Then, $\mathbf \z^{\boldsymbol{\ell}}_i \in \Delta_i \cap \{0,1\}^{d_i-1}$ for every $i\in N$. For $\alpha = \{(i,r_i) \mid i \in I \} \in \mathcal{A}$, we have
\[
L(\boldsymbol{\Gamma}^c)_{\alpha} = \sum_{\substack{k_i>r_i\\i\in I}}
\Gamma^c_{\{(i,k_i):i\in I\}} = \sum_{\boldsymbol{\ell} \in \prod_{i\in N}[d_i]}
\theta^{\boldsymbol{\ell}}
\prod_{i\in I}\mathbbm{1}\{\ell_i>r_i\} = \sum_{\boldsymbol{\ell} \in \prod_{i\in N}[d_i]}
\theta^{\boldsymbol{\ell}}\prod_{(i,k) \in \alpha}
z^{\boldsymbol{\ell}}_{ik},
\]
where the second equality holds from~\eqref{eq:flow-distribution}, and the first and last equality follows by definition. Therefore, we can conclude that 
\[
L(\boldsymbol{\Gamma}^c) = \sum_{\boldsymbol{\ell} \in \prod_{i \in N}[d_i]} \theta^{\boldsymbol{\ell}} \biggl( \prod_{(i,r) \in \alpha} z^{\boldsymbol{\ell}}_{ir} \biggr)_{\alpha \in \mathcal{A}} \in \conv(W^c_\emptyset \cdot \mathcal{M}^n),
\]
where the membership follows from that $\theta^{\boldsymbol{\ell}} \geq 0$ and $\sum_{\boldsymbol{\ell}} \theta^{\boldsymbol{\ell}} = W^c_{\emptyset}$, and $ \bigl( \prod_{(i,r) \in \alpha} z^{\boldsymbol{\ell}}_{ir} \bigr)_{\alpha \in \mathcal{A}} \in \mathcal{M}^n$.  \Halmos

\begin{lemma}\label{lemma:independent}
    The set 
    \[
    \mathcal{M}^{n} := \biggl\{ \boldsymbol{Z} \biggm| \z_i \in \Delta^{d_i-1} \cap \{0,1\}^{d_i-1} \for i \in N,\  Z_\alpha = \prod_{(i,k) \in \alpha}z_{ik} \for \alpha \in  \mathcal{A}  \biggr\}.
\]
contains $\prod_{i\in N}d_i$ affinely 
points. 
\end{lemma}
\noindent \textbf{Proof of Lemma~\ref{lemma:independent}} First, we introduce a matrix representation of points in $\mathcal{M}^n$.  For each item $i \in N$ and each price level $\ell_i \in [d_i]$, define a binary variable $\z^{\ell_i}_i$ by 
\[
\z^{\ell_i}_{ik} = 1 \quad \for k < \ell_i \qquad \text{ and } \qquad \z^{\ell_i}_{ik} = 0 \quad \for k \geq \ell_i. 
\]
For $\boldsymbol{\ell} \in \mathcal{L}:=  \prod_{i \in N} [d_i]$, define $\boldsymbol{Z}^{\boldsymbol{\ell}}$ as follows
\[
Z^{\boldsymbol{\ell}}_{\alpha} = \prod_{(i,k) \in \alpha} z^{\ell_i}_{ik} \quad \for \alpha \in \mathcal{A}. 
\]
Let $\boldsymbol{Q}$ the matrix representation of $\boldsymbol{Z}^{\boldsymbol{\ell}}$ with $\boldsymbol{\ell} \in \mathcal{L}$. More specifically, $\boldsymbol{Q}$ is the $\prod_{i\in N}d_i$ by $\prod_{i\in N}d_i$ square matrix whose columns are indexed by $\boldsymbol{\ell} \in \mathcal{L}$, whose rows are indexed by $\alpha \in \mathcal{A}$, and whose $(\alpha, \boldsymbol{\ell})$ entry is $Z^{\boldsymbol{\ell}}_\alpha$. We will prove that $ Q$ is non-singular. This implies that  $\{\boldsymbol{Z}^{\boldsymbol{\ell}} \mid \boldsymbol{\ell} \in  \mathcal{L} \}$ are affinely independent since the row corresponding with $\alpha = \emptyset$ of $\boldsymbol{Q}$ is the all-ones vector. 

Now, we prove by induction on the number of items $n$ that  $\boldsymbol{Q}$ is non-singular. For $n=1$, the matrix is the upper triangular:
\[
\boldsymbol{Q}=
\begin{pmatrix}
1&1&\cdots&1\\
0&1&\cdots&1\\
\vdots&\vdots&\ddots&\vdots\\
0&0&\cdots&1
\end{pmatrix},
\]
which is non-singular. Suppose the claim holds for $n-1$ items. Now, we order the rows $\alpha  \in \mathcal{A}$ associated with product $n$ by $k=0,1,\ldots,d_n-1$, where \(k=0\) denotes that no unary coordinate of product $n$ is selected, obtaining
\[
\boldsymbol{Q}=
\begin{pmatrix}
\boldsymbol{Q}^{(-n)}&\boldsymbol{Q}^{(-n)}&\cdots& \boldsymbol{Q}^{(-n)}\\
0&\boldsymbol{Q}^{(-n)}&\cdots& \boldsymbol{Q}^{(-n)}\\
\vdots&\vdots&\ddots&\vdots\\
0&0&\cdots&\boldsymbol{Q}^{(-n)}
\end{pmatrix},
\]
where $\boldsymbol{Q}^{(-n)}$ is the matrix associated with the first \(n-1\) products, and block in row $k$ and column $\ell_n$ equals \(\boldsymbol Q^{(-n)}\) if \(k<\ell_n\), and zero otherwise. Hence, $\det(\boldsymbol{Q}) = \det(\boldsymbol{Q}^{(-n)})^{d_n} \neq 0$, showing the non-singularity of $\boldsymbol{Q}$.  \Halmos

\noindent \textbf{Proof of Theorem~\ref{hierarchy}}  For each segment $c \in M$  and degree $\tau \geq 2$, let 
\[
\mathcal{S}^{\tau,c}: = \Bigl\{(\boldsymbol{Z}, \boldsymbol{W}^c) \Bigm| \z_i \in \Delta^{d_i-1} \cap \{0,1\}^{d_i-1} \for i \in N,\ (Z_{\alpha}, W^c_{\alpha}) \text{ satisfies }~\eqref{eq:high-order} \for \alpha \in \mathcal{A} \text{ with } \vert \alpha \vert \leq \tau \Bigr\},
\]
and define 
\[
\mathcal{S}^\tau := \Bigl\{(\boldsymbol{Z}, \boldsymbol{W}^1, \ldots, \boldsymbol{W}^m ) \Bigm| (\boldsymbol{Z}, \boldsymbol{W}^c) \in \mathcal{S}^{\tau, c} \for c \in M \Bigr\}. 
\]
By construction, the pricing problem~\eqref{price optimization}
is equivalent to maximizing the objective function in
Theorem~\ref{hierarchy}  over $\mathcal{S}^\tau$. 

We first establish that $R^*_\tau \geq R^*$. For each segment $c \in M$, by construction, constraints in~\eqref{eq:highorder} yield a relaxation of $\mathcal{S}^{\tau, c}$. Moreover,  by Lemma~\ref{lemma:multilinearhull}, $L(P^{\tau+1}(W^c_{\emptyset}))$ is a relaxation of constraint~\eqref{eq:highorderScaleF}. Therefore, the feasible region of the formulation defining $R_\tau^*$ is a relaxation of $\mathcal{S}^{\tau}$. Hence,  $R^*_\tau \geq R^*$. 

Second,  $R^*_{\tau} \geq R^*_{\tau + 1}$ holds because the formulation defining $R^*_{\tau+1}$ is obtained by imposing additional constraints on the formulation defining $R^*_{\tau}$. 

Last, to show $R^*_{n} = R^*$, it suffices to show that the feasible region of the formulation defining $R^*_n$  describes the convex hull of $\mathcal{S}^n$.  Observe that for each $c \in M$, the projection of $\mathcal{S}^{n,c}$ onto the $\boldsymbol{Z}$ variable space is $\mathcal{M}^n$, which, by Lemma~\ref{lemma:independent}, is a set of affinely independent points. Thus, it follows from Lemma 7 of~\cite{he2026discretenessconvexitypromotionplanning} that 
\[
\conv(\mathcal{S}^n) := \Bigl\{(\boldsymbol{Z}, \boldsymbol{W}^1, \ldots, \boldsymbol{W}^m ) \Bigm| (\boldsymbol{Z}, \boldsymbol{W}^c) \in \conv(\mathcal{S}^{n, c}) \for c \in M \Bigr\}. 
\]
By using Theorem 2 of~\cite{he2025convexification}, we obtain the following description of $\conv(\mathcal{S}^{n,c})$,
\begin{equation}\label{eq:multi}
\begin{aligned}
&( \boldsymbol{Z}, \boldsymbol{W}^c) \in W^c_{\emptyset} \cdot \conv \Bigl\{\bigl( \boldsymbol{Z} \cdot D^c(\z), \boldsymbol{Z} \bigr)   \Bigm|  \boldsymbol{Z} \in \mathcal{M}^n \Bigr\}   \\
&W^c_{\emptyset} \geq 0 \text{ and } W^c_{\emptyset} + \sum_{i \in N} \Bigl( \delta^c_{i1} \cdot W^c_{\emptyset} + \sum_{k = 2}^{d_i}\delta^c_{ik} \cdot W^c_{i,k-1}  \Bigr) = 1.
\end{aligned}
\end{equation}
By using $z_{ik}z_{ir} = z_{i, \max\{k,r\}}$, $\boldsymbol{Z}\cdot D^c(\z)$ can be expressed as a linear transformation of $\boldsymbol{Z}$. Since convexification commutes with linear transformation, the conic constraint in~\eqref{eq:multi} can be described as a linear transformation of $\conv(\mathcal{M}^n)$. This yields precisely the constraints in~\eqref{eq:highorderScaleF} and~\eqref{eq:highorderNormalize}. Finally, by Lemma~\ref{lemma:multilinearhull}, the flow polytope provides an exact description of the scaled multilinear constraints at degree $n$. Hence, the formulation defining $R_n^*$ describes $\conv(\mathcal{S}^n)$, showing $R^*_n = R^*$. \Halmos

\section{Proof of Proposition~\ref{cpd better than mcc}}
It suffices to prove that each inequality in~\eqref{mcc} is implied by an inequality in~\eqref{eq:coordinating-1-1}-\eqref{eq:coordinating-2-2}.
Since $\delta^c_{jr}\leq 0$ for every $c\in M$,
$j\in N$, $r\in\{2,\cdots,d_j\}$, and $1+\sum_{j\in N}\sum_{r\in [d_j]}\delta^c_{jr} > 0$ for every $c\in M$, we have 
\[\underline \rho^c = \frac{1}{1+\sum_{j\in N}\delta^c_{j1}},\qquad \overline \rho^c = \frac{1}{1+\sum_{j\in N}\sum_{r\in [d_j]}\delta^c_{jr}}\]
With the explicit bounds, we check the inequalities in~\eqref{mcc} one by one.
$\underline\rho^c z_{ik}\leq y^c_{ik}$ can be written as $z_{ik}\leq (1+\sum_{j\in N}\delta^c_{j1}) y^c_{ik}$. 
Since $\sum_{r=1}^{d_i-1}\delta^c_{ir+1}y^c_{i,\max\{k,r\}}\leq 0$, $\underline\rho^c z_{ik}\leq y^c_{ik}$ is implied by constraint~\eqref{eq:coordinating-1-2}.
$y^c_{ik}\leq\overline\rho^c z_{ik}$ can be written as $z_{ik}\geq (1+\sum_{j\in N}\sum_{r\in [d_j]}\delta^c_{jr})y^c_{ik}$. Since $\sum_{r=1}^{d_i-1}\delta^c_{ir+1}y^c_{i,\max\{k,r\} } \geq \sum_{r=1}^{d_i-1}\delta^c_{ir+1}y^c_{ik}$, $y^c_{ik}\leq\overline\rho^c z_{ik}$ is implied by~\eqref{eq:coordinating-2-1}.
$\rho^c-\overline\rho^c(1-z_{ik})\leq y^c_{ik}$ can be written as 
\[1-z_{ik} \geq (1+\sum_{j\in N}\sum_{r\in [d_j]}\delta^c_{jr})(\rho^c-y^c_{ik})\]
Using the constraint~\eqref{eq:mip-1}, it can be further rewritten as follows:
\[\rho^c + \sum_{j \in N} \Bigl( \delta^c_{j1} \cdot \rho^c+ \sum_{r = 1}^{d_j-1}\delta^c_{jr+1} \cdot y^c_{jr} \Bigr) - z_{ik}\geq (1+\sum_{j\in N}\sum_{r\in [d_j]}\delta^c_{jr})(\rho^c-y^c_{ik})\]
We can then organize it into
\[z_{ik}\leq (1+\sum_{j\in N}\delta^c_{j1})y^c_{ik}+\sum_{j\in N}\sum_{r=1}^{d_j-1} \delta^c_{jr+1}(y^c_{ik}+y^c_{jr} - \rho^c)\]
Since $y^c_{i,\max\{k,r\}} = \min\{y^c_{ik},y^c_{ir}\} = y^c_{ik}+y^c_{ir}-\max\{y^c_{ik},y^c_{ir}\}\geq y^c_{ik}+y^c_{ir}-\rho^c$, we have $\sum_{r=1}^{d_i-1}\delta^c_{ir+1}y^c_{i,\max\{k,r\} } \leq \sum_{r=1}^{d_i-1}\delta^c_{ir+1}(y^c_{ik}+y^c_{ir}-\rho^c)$, thus $\rho^c-\overline\rho^c(1-z_{ik})\leq y^c_{ik}$ can be implied by constraint~\eqref{eq:coordinating-1-1}.
$y^c_{ik}\leq\rho^c-\underline\rho^c(1-z_{ik})$ can be written as
\[1-z_{ik}\leq (1+\sum_{j\in N}\delta^c_{j1})(\rho^c - y^c_{ik})\]
Using the constraint~\eqref{eq:mip-1}, it can be further rewritten as follows:
\[\rho^c + \sum_{j \in N} \Bigl( \delta^c_{j1} \cdot \rho^c+ \sum_{r = 1}^{d_j-1}\delta^c_{jr+1} \cdot y^c_{jr} \Bigr) - z_{ik}\leq  (1+\sum_{j\in N}\delta^c_{j1})(\rho^c - y^c_{ik})\]
We can then organize it into
\[z_{ik} \geq \biggl(1 + \sum_{j \in N}\delta^c_{j1} \biggr) y^c_{ik} + \sum_{j \in N} \sum_{r =1}^{d_j-1} \delta^c_{jr+1}y^c_{jr}.\]
Since $\sum_{r=1}^{d_i-1}\delta^c_{ir+1}y^c_{i,\max\{k,r\} } \geq \sum_{r=1}^{d_i-1}\delta^c_{ir+1}y^c_{ir}$, $y^c_{ik}\leq\rho^c-\underline\rho^c(1-z_{ik})$ is implied by constraint~\eqref{eq:coordinating-2-2}.\Halmos
\section{Proof of Proposition~\ref{necessarycut}}
Since $\delta^c_{\ell t}\le0$ for $t\in \{2,3,\cdots,d_\ell\}$, $m^{\bar z}_{jr,\ell t}(\y^c)\ge\min\{y^c_{jr},y^c_{\ell t}\}$, and $h^{\bar z}_{jr,\ell t}(\y^c,\rho^c)\le\max\{0,y^c_{jr}+y^c_{\ell t}-\rho^c\}$, the selected inequalities relax constraints~\eqref{eq:coordinating}. Hence $\mathcal P_c\subseteq\mathcal P_c(\bar \z)$.

Let $(\bar{\z},\bar{\mu}^c,\rho^c,\y^c)$ satisfy the constraints in $\mathcal P_c(\bar{\z})$. If $\bar{z}_{jr} = 0$, we have
\begin{align*}
    0 = \bar{z}_{jr}&\geq
\left(1+\sum_{\ell\in N}\delta^c_{\ell1}\right)y_{jr}^c
+\sum_{t=1}^{d_j-1}\delta^c_{jt+1}y_{j,\max\{r,t\}}^c
+\sum_{\substack{\ell\neq j\\t\in[d_\ell-1]}}
\delta^c_{\ell t+1}m_{jr,\ell t}^{\bar z}(\y^c)\\
&\ge \left(1+\sum_{\ell\in N}\delta^c_{\ell1}\right)y_{jr}^c
+\sum_{t=1}^{d_j-1}\delta^c_{jt+1}y_{jr}^c
+\sum_{\substack{\ell\neq j\\t\in[d_\ell-1]}}
\delta^c_{\ell t+1}y^c_{jr},
\end{align*}
where the first inequality follows from \eqref{eq:zjrLowerDecomp}. The second inequality is because $\delta^c_{jt+1} \le 0$ and $y^c_{jr}-y^c_{j\max\{r,t\}}\ge 0$ and because $m_{jr,\ell t}^{\bar z}(\y^c) = y^c_{jr}$ whenever $\bar{z}_{jr}\le \bar{z}_{\ell t}$. Therefore, it follows from \eqref{eq:mip-2} that $y^c_{jr} = 0$. 
On the other hand, if $\bar{z}_{jr} = 1$, we have:
\begin{align*}
0 & \geq 1-\left(\left(1+\sum_{\ell\in N}\delta^c_{\ell1}\right)y_{jr}^c
+\sum_{t=1}^{d_j-1}\delta^c_{jt+1}y_{j,\max\{r,t\}}^c
+\sum_{\substack{\ell\neq j\\t\in[d_\ell-1]}}
\delta^c_{\ell t+1}h_{jr,\ell t}^{\bar z}(\y^c,\rho^c)\right)\\
& = \left(1+\sum_{\ell\in N}\delta^c_{\ell1}\right)(\rho^c - y_{jr}^c)
+\sum_{t=1}^{d_j-1}\delta^c_{jt+1}(y^c_{jt} - y_{j,\max\{r,t\}}^c)
+\sum_{\substack{\ell\neq j\\t\in[d_\ell-1]}}
\delta^c_{\ell t+1}(y^c_{lt} - h_{jr,\ell t}^{\bar z}(\y^c,\rho^c))\\
&\ge \left(1+\sum_{\ell\in N}\delta^c_{\ell1}\right)(\rho^c - y_{jr}^c)
+\sum_{t=1}^{d_j-1}\delta^c_{jt+1}(\rho^c - y_{jr}^c)
+\sum_{\substack{\ell\neq j\\t\in[d_\ell-1]}}
\delta^c_{\ell t+1}(\rho^c - y_{jr}^c),
\end{align*}
where the first inequality follows from~\eqref{eq:zjrUpperDecomp} and the first equality follows from~\eqref{eq:mip-1}. Due to~\eqref{eq:mip-2}, $(y^c_{jt} - y^c_{j,\max\{r,t\}}) - (\rho^c - y^c_{jr}) = y^c_{j,\min\{r,t\}}-\rho^c\leq 0$, and 
\[
\bigl(y^c_{\ell t}-h^{\bar{ z}}_{jr,\ell t}(\y^c,\rho^c)\bigr)
-(\rho^c-y^c_{jr})
=
\begin{cases}
y^c_{jr}-\rho^c, & \bar z_{\ell t}=0,\\
0, & \bar z_{\ell t}=1
\end{cases}
\;\le 0.
\]
Then, the last inequality follows from $\delta^c_{jr+1},\delta^c_{\ell t+1}\leq 0$. Therefore, it follows from \eqref{eq:mip-2} that $y^c_{jr} = \rho^c$. 

Combining both cases, $y^c_{jr} = \rho^c \bar{z}_{jr}$. It follows from \eqref{eq:mip-1} that $\rho^c = \frac{1}{1+ \sum_{\ell \in N}\left(\delta_{\ell 1} + \sum_{t\in [d_\ell - 1]}\delta_{\ell t+1}\bar{z}_{\ell t}\right)}$. Then, we obtain that:
\begin{equation*}
    \bar{\mu}^c = \sum_{\ell\in N} \left(\gamma^c_{\ell 1}\rho^c + \sum_{t\in [d_\ell - 1]}\gamma^c_{\ell t+1}y^c_{\ell t}\right) = \sum_{\ell \in N}\frac{\gamma^c_{\ell 1} + \sum_{t\in [d_\ell -1]}\gamma^c_{\ell t+1}y^c_{\ell t}}{1+\sum_{j'\in N}\left(\delta_{j'1}+\sum_{t\in [d_{j'}-1]}\delta_{j't+1}\bar{z}_{j't}\right)}.
\end{equation*}
In other words, $(\bar{\z},\bar{\mu}^c)\in \mathcal \proj_{(\z,\mu^c)}\mathcal P_c$. 

\section{Implementation Choices}\label{sec:implementation}

The algorithm implementation makes several implementation choices that we detail here.
First, we do not impose the integrality restrictions on \(\z\) at every iteration. Instead, most iterations solve the continuous relaxation of the master problem and generate cuts separating the resulting fractional solutions. The integrality restrictions are periodically restored to update the incumbent solution. This strategy improves the efficiency of cut-generation in two ways:
solving an LP is cheaper than solving the integer master problem, and separating fractional master solutions targets weaknesses of the LP relaxation. Restoring integrality, after two iterations of the LP, produces incumbent solutions that improve the lower bound. 

Second, our algorithm can be viewed as a Benders Decomposition~(\cite{benders1962partitioning,fischetti2010note}), so we also use the trust-region constraint used to stabilize Benders decomposition~(\cite{linderoth2003decomposition}). Let \(\hat \z\) denote the current candidate solution, and let $d$ be the step size. We restrict the next candidate solution by
\[
\sum_{j\in N}\sum_{r=1}^{d_j-1}
|z_{jr}-\hat z_{jr}|^2\leq d^2.
\]
The constraint prevents large oscillations of the master solution and focuses cut generation on a neighborhood of a promising incumbent. The constraint is disabled during IP iterations. We choose: 
\[d=\max\{3,0.003n_{\z}\}\]
where $n_{\z}$ is the number of coordinates of $\z$.

Third, we selectively include linking inequalities in the separation
subproblems rather than constructing them for every $z_{jr}$.
We solve an auxiliary master problem that maximizes the objective function for segment $c$.
Let $\bar x^{c}$ denote the prices implied by this solution, and let
$\bar \x$ denote those implied by the current candidate $\bar \z$.
We include the linking inequalities associated with $z_{jr}$ only if
\[
    |\bar x^{c}_j-\bar x_j| \geq 10^{-3},
\]
This heuristic rule prioritizes coordinates with a discrepancy between the segment-specific master solution and the current candidate. The selection is recomputed at each separation call, reducing the number of linking inequalities and associated multipliers in the cut-generating problem.

\section{Benefits of the Decomposition Algorithm}\label{sec:neceofdec}
\subsection{Large Scale Instances}
We still generate synthetic instances using the segment-specific utility function
\[
U_{cj}(p_j)=a_{cj}-b_{cj}p_j,
\]
where \(a_{cj}\in[-7,7]\) and \(b_{cj}\in[0.001,0.1]\), and we normalize segment weights $\omega_c$ after drawing them independently from \([0,1]\).

We conduct experiments with
\[
(|M|,|N|)\in\{(100,100),(500,20)\},
\qquad L=10,
\]
where $L$ denotes the number of admissible price levels of each product.
For each configuration, we generate 20 independent random instances.
We compare Price Decomposition presented in Algorithm~\ref{alg:pipeline} (PD) with directly solving the two MILP formulations: \eqref{CPD} and \eqref{fracprog}. All methods use a relative primal--dual gap tolerance of $0.01\%$ and a time limit of $3600$ seconds.

\begin{figure}
    \centering
    \includegraphics[width=\linewidth]{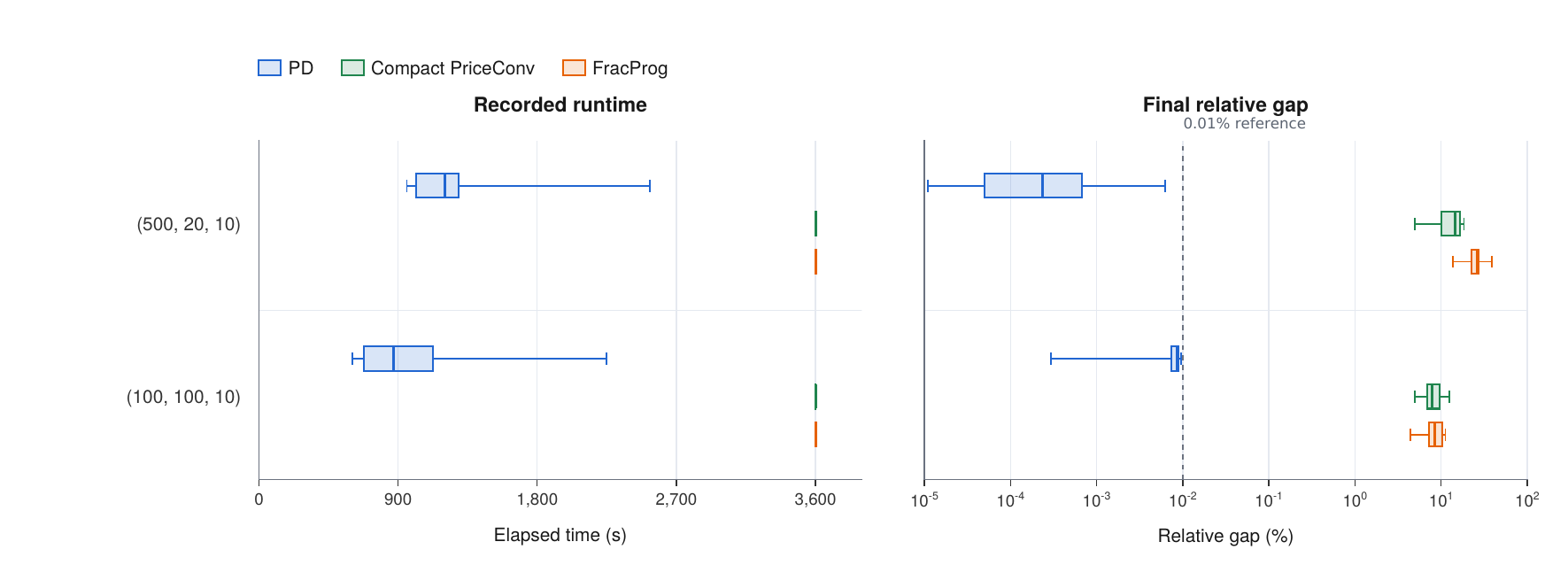}
     \caption{Runtime and relative gap of PD,~\eqref{CPD}, and~\eqref{fracprog} on synthetic instances, with 20 randomly generated instances per configuration.}
    \label{fig:unconstrainedlargebox}
\end{figure}

Figure~\ref{fig:unconstrainedlargebox} shows that PD consistently solves all tested instances to the prescribed $0.01\%$ relative optimality gap in both configurations. Its median runtime is approximately $1202$ seconds for $(M,N,L)=(500,20,10)$ and $873$ seconds for $(100,100,10)$, well below the $3600$-second time limit. In contrast, neither \eqref{CPD} nor \eqref{fracprog} solves any tested instance to this tolerance within the time limit. For $(500,20,10)$, their median final gaps remain approximately $14.5\%$ and $26.1\%$, respectively; for $(100,100,10)$, their median final gaps are approximately $7.9\%$ and $8.4\%$, respectively. Thus, PD reliably certifies global optimality to the prescribed tolerance, whereas \eqref{CPD} and \eqref{fracprog} retain substantial optimality gaps.

\begin{figure}
    \centering
    \includegraphics[width=\linewidth]{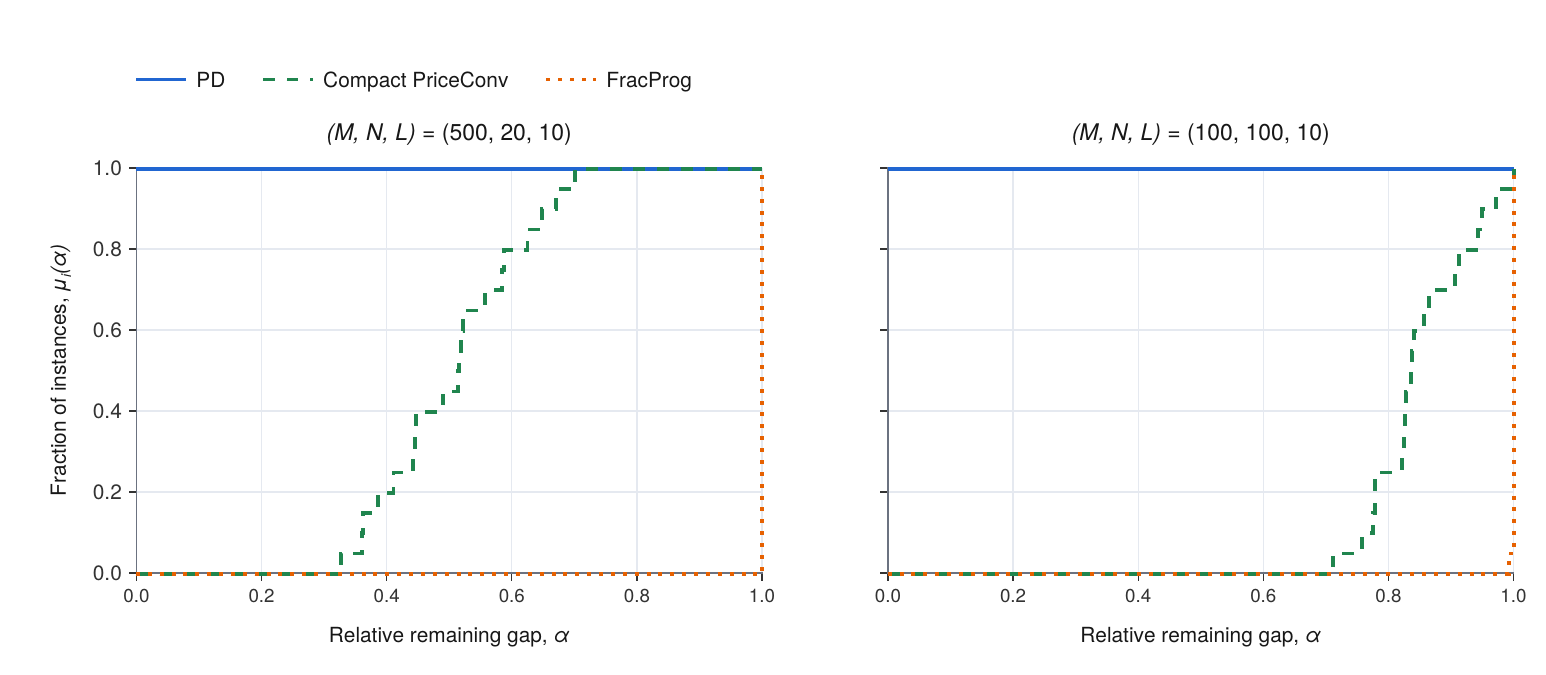}
    \caption{Relative remaining gap of PD,~\eqref{CPD}, and \eqref{fracprog} for $(M,N,L)=(500,20,10)$ and $(100,100,10)$, with 20 randomly generated instances per configuration.}
    \label{fig:unconstrainedlargetail}
\end{figure}
Figure~\ref{fig:unconstrainedlargetail} illustrates the relative remaining gaps of the three methods. 
For each instance $p\in\mathcal P$, let $v_{p,i}$ denote the final dual upper bound obtained by method $i$, and let $u_p$ be the largest primal objective value obtained by PD, \eqref{CPD}, and \eqref{fracprog}. We define the relative remaining gap as
\[
r_{p,i}=\frac{v_{p,i}-u_p}{\max_{i'}v_{p,i'}-u_p}.
\]
The relative remaining gap is set to zero whenever the corresponding relative optimality gap is below $0.01\%$.
The function $\mu_i(\alpha)$ is the empirical cumulative distribution function of the relative remaining gap for method $i$. 
For any threshold $\alpha$, $\mu_i(\alpha)$ gives the fraction of instances for which method $i$ achieves a relative remaining gap no greater than $\alpha$. 
Higher curves indicate better performance.
The figure shows that PD closes the optimality gap to within the prescribed $0.01\%$ tolerance on every tested instance, with $\mu_{\mathrm{PD}}(0)=1$ in both configurations. More informatively, Figure~\ref{fig:unconstrainedlargetail} distinguishes the performance of \eqref{CPD} and \eqref{fracprog} even though neither method reaches this tolerance. For $(M,N,L)=(500,20,10)$, \eqref{CPD} achieves relative remaining gaps of approximately $0.33$--$0.70$, whereas \eqref{fracprog} has a relative remaining gap of $1$ on every instance. Thus, \eqref{CPD} substantially reduces the gap left by \eqref{fracprog}. For $(100,100,10)$, \eqref{CPD}'s curve still lies above that of \eqref{fracprog}, although the advantage is less pronounced. These results show that, while PD reliably closes the gap, \eqref{CPD} provides meaningfully tighter final bounds than \eqref{fracprog}, particularly in the configuration with more customer segments.

\subsection{Pricing with Limited Change}
When additional constraints arising in practice are imposed, directly solving either the \eqref{CPD} or \eqref{fracprog} formulation may fail to converge.
In practice, changes to the current price vector cannot be arbitrary or substantial for two reasons. First, extensive price adjustments based on a potentially misspecified demand model entail risk. Second, implementing price changes incurs operational costs. However, it is reasonable to adjust the prices of a limited number of products. Formally, suppose $x_{j}^*$ is the current price for the $j$-th product and we can only modify up to $N_0$ products, it imposes the following constraint:
\[\sum_{j\in N}\lVert x_{j}-x_{j}^*\rVert_0\leq N_0\]

For each product \(j\), let \(\{p_{j1},p_{j2},\ldots,p_{jd_j}\}\) be the set of admissible price levels, and assume that the current price belongs to this set, with \(x_j^*=p_{jk_j}\).
Let $\z_j$ be the unary binarization of $x_j$, then $\lambda(\z_j) := (1-z_{j1},z_{j1}-z_{j2},z_{j2}-z_{j3},\cdots,z_{jd_j-2}-z_{jd_{j}-1},z_{jd_j-1})$ is the one-hot binarization of $x_j$ where $\lambda(\z_{j})_k = 1$ if and only if $x_{j} = p_{jk}$. Thus the price of $j$-th product is changed if and only if $1-\lambda(\z_j)_{k_j} = 1$. Then the price change limitation constraint can be written as follows:
\[\sum_{j\in N}(1-\lambda(\z_j)_{k_j})\leq N_0\tag{Price Change Limitation}\label{pricechangelimitation}\]
And its scaled version can be written as follows:
\[\sum_{j\in N}(\rho^c-\lambda^c(\y^c_j,\rho^c)_{k_j})\leq N_0\rho^c\tag{Scaled Price Change Limitation}\label{scaledpricechangelimitation}\]
where $\lambda^c(\y^c_j,\rho^c) := (\rho^c-y^c_{j1},y^c_{j1}-y^c_{j2},y^c_{j2}-y^c_{j3},\cdots,y^c_{jd_j-2}-y^c_{jd_{j}-1},y^c_{jd_j-1})$.

To incorporate price change limitation into~\eqref{CPD} and~\eqref{fracprog}, we add constraints~\eqref{pricechangelimitation} and~\eqref{scaledpricechangelimitation} for $c\in M$ to the formulations.
To incorporate price change limitation into algorithm~\ref{alg:pipeline}, we add constraint~\eqref{pricechangelimitation} into initial feasible region $\mathcal C$ and add constraint~\eqref{scaledpricechangelimitation} to $\mathcal P_c$ for $c\in M$.

\begin{figure}
    \centering
    \includegraphics[width=\linewidth]{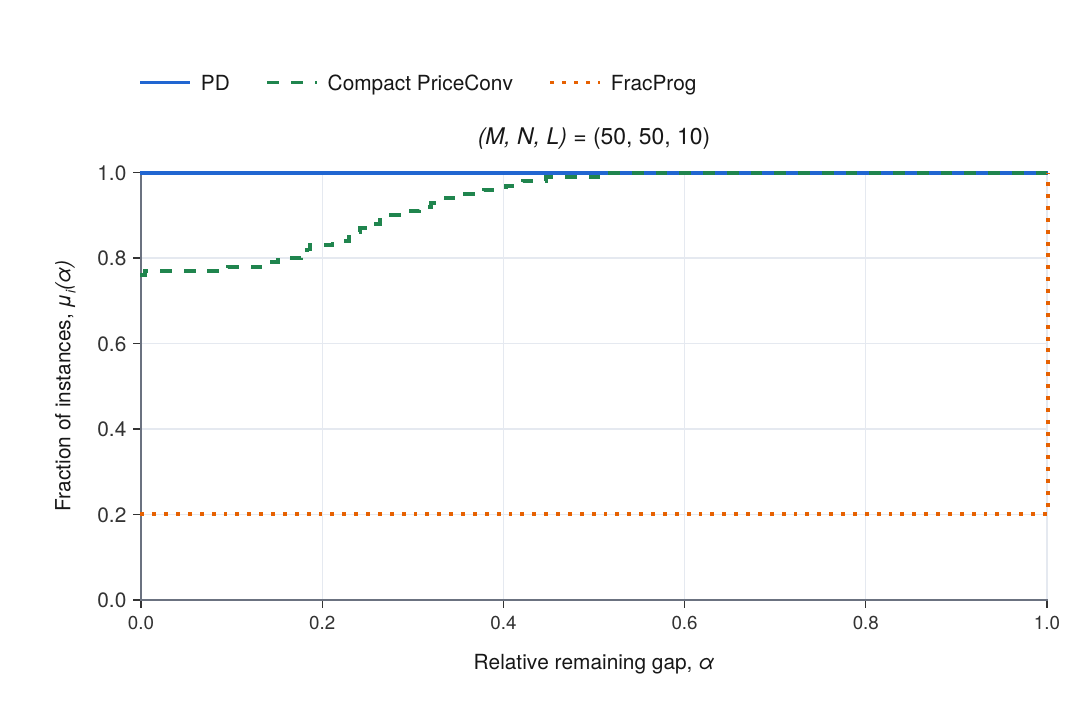}
    \caption{Relative remaining gap of PD, \eqref{CPD}, and \eqref{fracprog} on synthetic instances with price change limitation. Based on $100$ instances.}
    \label{fig:constrainedtail}
\end{figure}

We uniformly draw \(a_{cj}\in[-1,1]\) (less heterogeneous intrinsic preferences) and \(b_{cj}\in[0.001,0.1]\), and we normalize segment weights $\omega_c$ after drawing them independently from \(U[0,1]\).
We compare their performance in a relatively small instance where
\[M=50,\qquad N=50,\qquad L=10,\qquad N_0 = 20.\]
The current price of each product is randomly picked from the admissible set of price levels.
All methods use a relative primal--dual gap tolerance of $0.01\%$ and a time limit of $3600$ seconds.
\begin{figure}
    \centering
    \includegraphics[width=\linewidth]{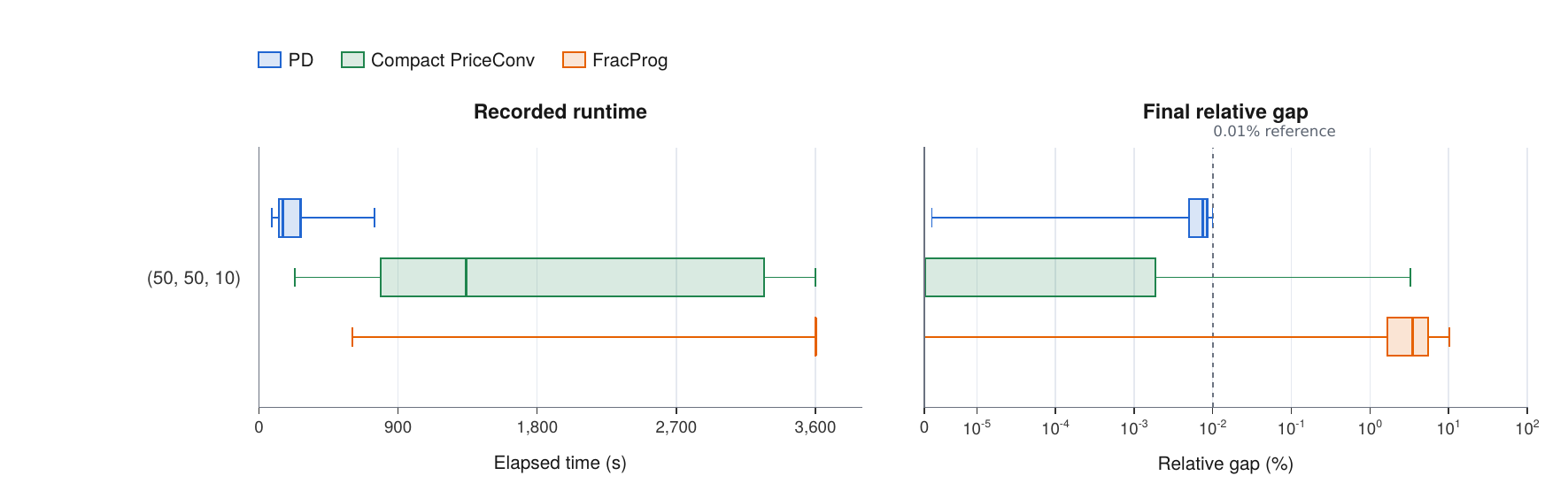}
    \caption{Runtime and relative gap of PD, \eqref{CPD}, and \eqref{fracprog} on synthetic instances with price change limitation. Based on $100$ instances.}
    \label{fig:constrainedbox}
\end{figure}
Figure~\ref{fig:constrainedbox} shows that with the price change limitation, PD consistently solves all 100 instances to the prescribed $0.01\%$ relative optimality gap, with a median runtime of approximately $156$ seconds. \eqref{CPD} also solves a majority of the instances, but at substantially greater computational cost and with much larger variation across instances: it reaches the target tolerance on 76 instances, while its median runtime is approximately $1339$ seconds. \eqref{fracprog} remains unable to solve most instances within the time limit, reaching the target tolerance on only 20 instances; its median runtime reaches the $3600$-second limit, with a median final gap of approximately $3.5\%$. Thus, PD delivers fast and consistent convergence, whereas \eqref{CPD} is considerably more expensive and less consistent, and \eqref{fracprog} leaves most instances unresolved.
Figure~\ref{fig:constrainedtail} shows that \eqref{CPD} provides substantially tighter bounds than \eqref{fracprog} even when it does not close the gap.
\section{Homogeneity is theoretically restrictive}\label{sec:homogeneityRestrictive}
Here, we show that choosing $b_{ci}$ to be independent of $c$ is also theoretically restrictive, and it reduces the flexibility of FML. In particular, such restrictive models can no longer uniformly approximate arbitrary random choice models unlike general FML models \citep{mcfadden2000mixed}. Even for a setup involving two consumer segments and two products, this approximation is not possible. We recall that, given a price vector $\x$, the probability of a consumer $c$ purchasing a product $i$ is written as follows:
\[\pi_{ci}(\x) = \frac{\exp(a_{ci}-b_{ci}x_i)}{1+\exp(a_{c1}-b_{c1}x_1)+\exp(a_{c2}-b_{c2}x_2)}\]
Then, the market share of product $i$ is given by:
\[
P_i = \sum_{c} \omega_c \pi_{ci}(\x).
\]
It follows that 
\[
\frac{\sfrac{\partial P_1}{\partial x_2}}{\sfrac{\partial P_2}{\partial x_1}} = \frac{\sum_c\omega_c b_{c2} \pi_{c1}(\x)\pi_{c2}(\x)}{\sum_c \omega_c b_{c1} \pi_{c1}(\x)\pi_{c2}(\x)}.
\]
Let $A_c = \omega_c \pi_{c1}(\x)\pi_{c2}(\x)$. Then, we have:
\begin{equation}\label{eq:ratioFractions}
\frac{\sfrac{\partial P_1}{\partial x_2}}{\sfrac{\partial P_2}{\partial x_1}} = \frac{A_1 b_{12}+A_2 b_{22}}{A_1 b_{11}+A_2 b_{21}}.
\end{equation}
Observe that this ratio depends on $(x_1,x_2)$. On the other hand, assume that one develops a model assuming that the price sensitivity $b$ is independent of $c$. Then, \eqref{eq:ratioFractions} simplifies so that $\frac{\sfrac{\partial P_1}{\partial x_2}}{\sfrac{\partial P_2}{\partial x_1}} = \frac{b_{2}}{b_{1}}$.
This quantity is independent of $(x_1,x_2)$, regardless of the number of segments in the model, making it impossible to uniformly approximate the original consumer choice model where this quantity varies with $(x_1,x_2)$. This shows that the restrictive model cannot uniformly approximate arbitrary RUM models.
\section{Upper and Lower Bounds on Optimal Prices}\label{bounds of prices}
This formulation is a straightforward extension of Proposition~1 in~\cite{van2022price}, in which the price sensitivity is assumed to be homogeneous across the customer segments.
\begin{proposition}
\label{prop:price-bounds}
Suppose that \(\x\in\mathbb R_+^n\), \(b_{cj}>0\) for all \(c\in M\) and \(j\in N\), and \(\omega_c>0\) for all \(c\in M\), with \(\sum_{c\in M}\omega_c=1\). Define
\[
B:=\max_{c\in M}\max_{\x\in\mathbb R_+^n}\sum_{i\in N}(x_i-\cost_i)\pi_{ci}(\x),\qquad \underline x_j:=\cost_j + \min_{c\in M}b_{cj}^{-1},\qquad \overline x_j:=\cost_j + \max_{c\in M}b_{cj}^{-1}+B.
\]
A global maximizer \(\x^*\) of \(\Pi(\x) = \sum_{c\in M}\omega_c\sum_{i\in N}(x_i-\cost_i) \cdot \pi_{ci}(\x)\) exists, and every global maximizer satisfies
\[
\underline x_j\leq x_j^*\leq\overline x_j,\qquad j\in N.
\]
\end{proposition}

\noindent{\bf Proof. }
For each \(c\in M\), define \(R_c(\x):=\sum_{i\in N}(x_i-\cost_i)\pi_{ci}(\x)\). 
First, if there exists product $j\in N$ such that $x_j< \cost_j$, raise all such prices to their respective costs. 
Each affected item then contributes zero rather than a negative amount to \(R_c\). Every other item with a nonnegative profit has a higher choice probability, since its own attraction value is unchanged and the sum of the attraction value is decreased. Thus this change strictly increases \(R_c\) for every \(c\in M\), and hence also increases \(\Pi\). 
We can therefore restrict attention to \(x_j\geq\cost_j\) for all \(j\in N\).
On this region,
\[
0\leq R_c(\x)\leq\sum_{i\in N}x_i\exp(a_{ci}-b_{ci}x_i)
\leq K_c:=\sum_{i\in N}\frac{\exp(a_{ci})}{e\cdot b_{ci}}.
\]
The MNL probabilities satisfy
\[
\frac{\partial \pi_{ci}(\x)}{\partial x_j}=
\begin{cases}
-b_{cj}\pi_{cj}(\x)(1-\pi_{cj}(\x)),&i=j,\\
b_{cj}\pi_{ci}(\x)\pi_{cj}(\x),&i\neq j.
\end{cases}
\]
Therefore,
\[
\frac{\partial R_c(\x)}{\partial x_j}
=b_{cj}\pi_{cj}(\x)\bigl(b_{cj}^{-1}+R_c(\x)-(x_j-\cost_j)\bigr).
\]
In particular, if \(x_j>\cost_j+b_{cj}^{-1}+K_c\), then \(\partial R_c(\x)/\partial x_j<0\). Reducing each such coordinate to its upper bound strictly increases \(R_c\). Thus the maximum of \(R_c\) is attained on the compact box \(\prod_{j\in N}[\cost_j,\cost_j+b_{cj}^{-1}+K_c]\), so \(B\) is finite and well defined.

Since \(\Pi(\x)=\sum_{c\in M}\omega_cR_c(\x)\), we have
\[
\frac{\partial\Pi(\x)}{\partial x_j}
=\sum_{c\in M}\omega_cb_{cj}\pi_{cj}(\x)
\bigl(b_{cj}^{-1}+R_c(\x)-(x_j-\cost_j)\bigr).
\]
If \(x_j<\underline x_j\), then \(b_{cj}^{-1}+R_c(\x)-(x_j-\cost_j)>0\) for every \(c\in M\), since \(R_c(\x)\geq0\). Hence \(\partial\Pi(\x)/\partial x_j>0\). Similarly, if \(x_j>\overline x_j\), then
\[
b_{cj}^{-1}+R_c(\x)-(x_j-\cost_j)
\leq \cost_j+\max_{c'\in M}b_{c'j}^{-1}+B-x_j<0
\]
for every \(c\in M\), and hence \(\partial\Pi(\x)/\partial x_j<0\).

For any price vector outside the box, move each violating coordinate, one at a time, to the nearest endpoint of its interval. By the derivative signs established above, each move strictly increases \(\Pi\). Thus every point outside the box is dominated by a point inside it. Maximizing \(\Pi\) over \(\mathbb R_+^n\) is therefore equivalent to maximizing it over this compact hyperrectangle. Since \(\Pi\) is continuous, the Weierstrass theorem guarantees the existence of a global maximizer, and every such maximizer satisfies the stated bounds. \(\Halmos\)
\section{Overview of Estimated Parameters}\label{overviewofpara}
Here, we summarize the parameters obtained when we fit an MNL model separately for each customer segment on the Expedia pricing data. For each product $j$ and parameter $\theta\in\{a,b\}$, let $(\theta_{1j},\dots,\theta_{mj})$ denote the $m$ segment-level estimates. Table~\ref{tab:expedia-parameter-summary} reports, for each center, the average over products of the per-product segment-level mean and standard deviation.

The parameters summarized here are calibrated to real choice data. The degree of parameter heterogeneity is derived from the data and the table confirms that there is significant variation in intrinsic preferences and price sensitivities that our methods can exploit.
\begin{table}[htbp]
\centering
\begin{tabular}{lrrrr}
\toprule
& \multicolumn{2}{c}{Intercept $a_{\cdot j}$}
& \multicolumn{2}{c}{Price sensitivity $b_{\cdot j}$} \\
\cmidrule(lr){2-3}\cmidrule(lr){4-5}
Center & Mean & Std & Mean & Std \\
\midrule
8467  & -0.144277 & 0.933319 & 0.012850 & \(1.82216 \times 10^{-3}\) \\
4562  & -0.126450 & 0.879466 & 0.012900 & \(1.82227 \times 10^{-3}\) \\
23246 & -0.163523 & 0.903224 & 0.012809 &  \(1.82624 \times 10^{-3}\) \\
\bottomrule
\end{tabular}
\caption{Mean and standard deviation of per-product estimated intercepts and price sensitivities within each Expedia destination center.}
\label{tab:expedia-parameter-summary}
\end{table}

\end{document}